\newcommand{\elsarticle}{}

\ifdefined \regComp

    \documentclass[12pt,final,leqno]{amsart}
\fi
\ifdefined \elsarticle

    \documentclass{elsarticle}

\fi

\RequirePackage{lmodern} 

\usepackage{amsfonts}
\usepackage{amsmath}
\usepackage{amssymb}
\usepackage{amsthm}
\usepackage{appendix}
\usepackage{changepage}
\usepackage{enumerate}
\usepackage{enumitem}
\usepackage{eucal}
\usepackage{latexsym}
\usepackage{mathrsfs}
\usepackage{scrextend}
\usepackage{textcomp}
\usepackage{todonotes}
\usepackage{verbatim}
\usepackage{pdflscape}
\usepackage{lscape}
\usepackage{chngcntr}

\usepackage{xfrac}

\relax							

\newlength{\pardefault}							
\newenvironment{psmallmatrix}
  {\left(\begin{smallmatrix}}
  {\end{smallmatrix}\right)}

\newcommand{\innerProductReg}[2]{\left\langle #1 , #2 \right\rangle}
\newcommand{\innerProductTri}[3]{\left\langle #1 , #2 \right\rangle_{#3}}

\newcommand{\norm}[1]{\| #1  \|}

\newcommand{\normTreA}[3]{ \left| \left| { #1 } \right| \right|  _{#2}^{#3} }

\newcommand{\bracketsA}[1]{\left( #1 \right)}
\newcommand{\bracketsB}[2]{\left( #1 \right)_{#2}}
\newcommand{\bracketsC}[3]{\left( #1 \right)_{#2}^{#3}}

\newcommand{\abs}[1]{\left| #1 \right|}

\newcommand{\imag}{\mathfrak{Im} ~}

\newcommand{\imagg}[1]{\mathfrak{Im} \left( #1  \right)}

\newcommand{\restrict}[1]{\raisebox{-.5ex}{$|$}_{#1}}
\newcommand\restr[2]{{
  \left.\kern-\nulldelimiterspace 
  #1 
  \vphantom{\big|} 
  \right|_{#2} 
  }}
\newcommand{\defEq}{\overset{\text{def} } {=}}

\newcommand{\refMult}[2]{(\ref{#1} - \ref{#2})}

\newcommand{\wt}[1]{\widetilde{#1}}

\newcommand{\tauBa}[1]{\tau{(#1)}}

\mathchardef\mhyphen="2D

\DeclareMathOperator*{\col}{Col}
\DeclareMathOperator*{\row}{Row}

\DeclareMathOperator*{\diag}{Diag}

\DeclareMathOperator*{\per}{Per}
\DeclareMathOperator*{\spec}{Spec}

\newcommand{\Card}[1]{\left| #1 \right|}
\newcommand{\Per}[1]{\per \left( #1 \right)}
\newcommand{\Spec}[1]{\spec \left( #1 \right)}

\newcommand{\Diag}[1]{\diag \left( #1 \right)}
\newcommand{\DiagTwo}[2]{\diag_{#1} \left( #2 \right)}

\newcommand{\btul}{{\bf u}^{\times}}

\newcommand{\forceOddPage}{
\ifodd\theCurrentPage\pagebreak\fi}

\newcommand{\forceEvenPage}{
\ifodd\value{\theCurrentPage}\newpage\else\pagebreak\fi}

\newcommand*\cleartoleftpage{%
  \clearpage
  \ifodd\value{page}\hbox{}\newpage\fi
}

\newcommand*\cleartorightpage{%
  \clearpage
  \unless\ifodd\value{page}\hbox{}\newpage\fi
}

\newcommand{\ignore}[1]{}

\newtheorem{Pa}{Paper}[section]
\newtheorem{Tm}[Pa]{{\bf Theorem}}

\newtheorem{Cy}[Pa]{{\bf Corollary}}
\newtheorem{Rk}[Pa]{{\bf Remark}}

\newtheorem{theorem}[Pa]{{\bf Theorem}}
\newtheorem{lem}[Pa]{{\bf Lemma}}
\newtheorem{corollary}[Pa]{{\bf Corollary}}
\newtheorem{prop}[Pa]{{\bf Proposition}}

\newtheorem{definition}[Pa]{{\bf Definition}}
\newtheorem{proposition}[Pa]{{\bf Proposition}}
\newtheorem{Notation}[Pa]{{\bf Notation}}
\newtheorem{Step}[]{{\bf Step}}

\newenvironment{pf}[1][\unskip]{
\par
\noindent
{\bf Proof #1:}
\noindent
}
{\hfill$\square$\\}

\def\R{\mathbb R}
\def\N{\mathbb N}
\def\C{\mathbb C}
\def\Z{\mathbb Z}

\def\P{\mathbb P}

\def\E{\mathcal E}

\def\z{\zeta}
\def\e{\varepsilon}

\usepackage{hyperref}
\hypersetup{hypertexnames=false} 

\hypersetup{
  colorlinks,
  citecolor	= violet,
  linkcolor	= blue,
  urlcolor	= red}

% ...

\usepackage{enumitem}
\numberwithin{equation}{section}
\usepackage{mathtools}
\usepackage{arydshln}

\begin{document}


\ifdefined \regComp

\title
[de Branges spaces on compact Riemann surfaces]
{de Branges spaces on compact Riemann surfaces and a Beurling-Lax type theorem}


\author[D. Alpay]{Daniel Alpay}
\address{(DA) Schmid College of Science and Technology,
Chapman University,
One University Drive
Orange, California 92866,
USA}
\email{alpay@chapman.edu}

\author[A. Pinhas]{Ariel Pinhas}
\address{(AP) Department of mathematics,
Ben-Gurion University of the Negev, P.O. Box
653, Beer-Sheva 8410501, Israel}
\email{arielp@post.bgu.ac.il}

\author[V. Vinnikov]{Victor Vinnikov}
\address{(VV) Department of mathematics,
Ben-Gurion University of the Negev, P.O. Box
653, Beer-Sheva 8410501, Israel}
\email{vinnikov@math.bgu.ac.il}

\date{}

\thanks{The first author thanks the Foster G. and Mary McGaw Professorship in Mathematical Sciences, which supported this research}
\thanks{The research of the second and the third authors was partially supported by the US--Israel Binational Science Foundation (BSF) Grant No. 2010432, 
Deutsche Forschungsgemeinschaft (DFG) Grant No. SCHW 1723/1-1, and Israel Science Foundation (ISF) Grant No. 2123/17.}
\fi

\ifdefined \elsarticle

	\title{de Branges spaces on compact Riemann surfaces and a Beurling-Lax type theorem}


	\author{Daniel Alpay\fnref{fn1}}
	\ead{alpay@chapman.edu}

	\address{Schmid College of Science and Technology,
	Chapman University,
	One University Drive
	Orange, California 92866,
	USA}

	\author{Ariel Pinhas\fnref{fn2}}
	\ead{arielp@post.bgu.ac.il}

	\address{Department of mathematics,
	Ben-Gurion University of the Negev, P.O. Box
	653, Beer-Sheva 84105, Israel}

	\author{Victor Vinnikov\fnref{fn2}}
	\ead{vinnikov@math.bgu.ac.il}

	\address{Department of mathematics,
	Ben-Gurion University of the Negev, P.O. Box
	653, Beer-Sheva 84105, Israel}

	\date{}

	\fntext[fn1]{The first author thanks the Foster G. and Mary McGaw Professorship in Mathematical Sciences, which supported this research}
	\fntext[fn2]{The research of the second and the third authors was partially supported by the US--Israel Binational Science Foundation (BSF) Grant No. 2010432, 
	Deutsche Forschungsgemeinschaft (DFG) Grant No. SCHW 1723/1-1, and Israel Science Foundation (ISF) Grant No. 2123/17.}

\fi

\begin{abstract}
Using the notion of commutative operator vessels, this work investigates de Branges-Rovnyak spaces
whose elements are sections of a line bundle of multiplicative half-order differentials on a compact real Riemann surface.
As a special case, we obtain a Beurling-Lax type theorem in the setting of the corresponding Hardy space on a finite bordered Riemann surface.
\end{abstract}

\ifdefined \regComp
\fi

\ifdefined \elsarticle


\fi

\maketitle


\setcounter{tocdepth}{1}
\tableofcontents


\section{Introduction}


\subsection{The classical case}
The theory of de Branges-Rovnyak spaces of analytic functions (see for instance \cite{alpay2017beurling,dbhsaf1,dbr1,dbr2}) 
allows us to prove Beurling-Lax type theorem in a one complex variable framework when leaving the setting of the Hardy space.
An illustrative example is given in the following theorem.
To give the statement, we must first recall that an element of $\mathbb{C}^{n \times n}$ which is
both selfadjoint and unitary is called a signature matrix.
Furthermore, for $\alpha\in\mathbb C$ and $f$ a matrix-valued function analytic in a neighborhood of $\alpha$, 
let $R_{\alpha}$ denote the resolvent operator at $\alpha$:
\begin{equation}
\label{resolventOp}
R_\alpha f(z) =
\begin{cases}
\,\,
\dfrac{f(z)-f(\alpha)}{z-\alpha},
\quad
z\not=\alpha\\
\,\,
f^\prime(\alpha),
\quad \quad \quad \quad
\,
z=\alpha.
\end{cases}
\end{equation}
It follows that the resolvent identity \cite[Theorem I]{dbhsaf1}
\[
R_{\alpha} - R_{\beta}= (\alpha - \beta)  R_{\alpha}R_{\beta}
\]
holds for any function analytic in a connected neighborhood of $\alpha$ and $\beta$.

\begin{theorem}
\label{1_1}
Let $\Omega$ be an open subset of the complex plane, symmetric with respect to the real line, and let $J$ be a signature matrix.
Let $\mathcal X$ be a reproducing kernel Hilbert space of functions analytic in $\Omega$.  Then the reproducing
kernel of $\mathcal{X}$  is of the form
\begin{equation}
\label{ktheta}
K_T(z,w) 
=
\frac{J - T(z) J T(w)^*}{-i(z-\overline{w}) },
\end{equation}
where $T$ is a $\mathbb C^{n\times n}$-valued function analytic in $\Omega$,
if and only if the following two conditions are fulfilled:
\begin{enumerate}
\item $R_\alpha \mathcal X\subset\mathcal X$ for all $\alpha\in\Omega$.
\item The structure identity
\begin{equation}
\label{eq:structure1}
\innerProductReg{R_\alpha f}{g} -
\innerProductReg{f}{R_\beta g} -
(\alpha - \overline{\beta}) \innerProductReg{R_\beta f}{R_\beta g}
=
i g(\beta)^* J f(\alpha)
\end{equation}
holds for all $\alpha,\beta \in \Omega$ and $f,g \in \mathcal X$.
\end{enumerate}
\label{thm:struct}
\end{theorem}
See \cite[Theorems III and IV]{dbhsaf1}, where $\Omega\cap\mathbb R$ is assumed non-empty; this restriction was later removed in \cite{rov-66} and in
\cite{ball-contrac} (for the corresponding theorem in the case of the disk, see \cite{ad-jfa}).
\smallskip

We remark that $T$ is $J$-contractive, that is, $T$ satisfies $T(z)J T(z)^*\le J$ for all $z$ in the intersection of $\Omega$ and the upper half-plane $\C_+$.
The positivity of the kernel \eqref{ktheta} in an open subset of the upper half-plane implies that $T$
has a meromorphic extension to $\mathbb C_+$ (and in fact, by reflection, to $\mathbb C\setminus \mathbb R$),
for which $K_T(z,w)$ is still positive definite.
\smallskip

If we further assume that $\infty\in\Omega$ and $\ker R_\alpha = \{0\}$ for an arbitrary $\alpha \in \Omega$
(and henceforth, for all $\alpha\in\Omega$), then the space  does not contain nonzero constants,
and $T$ is of the form
\begin{equation}
\label{eqCfCol}
T(z) = I - iC(zI - A) ^ {-1} C^* J
\end{equation}
where
\begin{equation}
\label{eqOpModelClass}
(A \, f)(z)=zf(z)-\lim_{w\rightarrow\infty}wf(w)\quad\text{and }\quad Cf=\lim_{w\rightarrow\infty}wf(w).
\end{equation}
Furthermore, the resolvent operator satisfies $R_{\alpha} = (A- \alpha I )^{-1}$.
\smallskip

In his fundamental paper \cite{MR0027954} Beurling introduced a characterization of 
invariant subspaces under the shift operator in the Hardy space $H^2(\mathbb{D})$.
These subspaces are characterized as the ones of the form $T H^2(\mathbb{D})$, where $T$ is an inner function. 
An inner function is, by definition, an analytic function on the unit disk such that $|T(z)| \leq 1$ for $|z|<1$ and 
with non-tangential boundary values (which exist almost everywhere since $T$ is bounded) having modulus one.
Important generalizations, the vector-valued case and the infinite dimensional case, were presented later by Lax \cite{MR0140931} and Halmos \cite{MR0152896}, respectively.
\smallskip

We note that \eqref{eq:structure1} is automatically satisfied in the Hardy space $\mathbf H^2(\mathbb C_+)$.
Therefore, applying Theorem \ref{1_1} (restricted to subspaces of $\mathbf H^2(\mathbb C_+)$)
on the orthogonal complement implies the Beurling-Lax theorem, under the hypothesis of a symmetric domain of analyticity.

\subsection{Beurling theorem on Riemann surfaces}
\label{sucSectionIntoCRS}

In the setting of finite bordered Riemann surfaces, 
several generalizations of Beurling's theorem were presented.
Sarason, in \cite{sarason2}, studied the invariant subspaces of certain multiplication operators on ${\bf L}^2$ in the case of the annulus.
Invariant subspaces in Hardy spaces over a finite bordered Riemann surface were studied by
Voichick \cite{MR0160920,voichick1966invariant}, Voichick-Zalcman \cite{MR0183883}, Hasumi \cite{MR0190790} and Forelli \cite{MR0193534} in the scalar-valued case,
while Abrahamse and Douglas \cite{AbDo} considered the vector-valued case
(see also Abrahamse \cite{Abr1} and Ball \cite{MR526287} for the closely related interpolation problems).
A dedicated investigation of invariant subspaces in Hardy space over a multiply connected 
planer domain, including corresponding operator theory and Riesz bases,
was carried out by
Fedorov and Pavlov \cite{MR1025160} and Fedorov \cite{MR1072300,MR1070485}.
In addition to lifting the problem to the universal covering,
they use extensively the representation of the domain as a
finite branched covering of the unit disk.
\smallskip

Neville \cite{MR0301206,MR0586558} and Hasumi \cite{MR0364647,MR0407283} obtained Beurling Lax type theorems 
for certain classes of infinitely connected open Riemann surfaces.
The crucial condition is that any flat complex line bundle on the open Riemann surface
admits bounded holomorphic sections; such open Riemann surfaces
were characterized completely by Widom \cite{Widom} and referred as Parreau-Widom type, see also \cite{MR0050023}.
Finally, we also mention the work of Kupin and Yuditskii \cite{MR1473264},
as they have studied interpolation problems on 
infinitely connected open Riemann surfaces of Parreau-Widom type.


\subsection{Compact real Riemann surfaces}
In the present paper we prove a counterpart of Theorem \ref{thm:struct} in the setting of compact real Riemann surface
(i.e. compact Riemann surface equipped with an anti-holomorphic involution $\tau$);
for a theorem in the finite dimensional case see \cite[Theorem 5.1]{av3}.
Taking the compact real Riemann surface to be the double of a finite bordered Riemann surface $\mathscr{S}$
we obtain, as a corollary, a Beurling Lax type characterization of invariant subspaces of a Hardy space on $\mathscr{S}$ under
a pair of certain multiplication operators.

The counterpart of the kernel $\frac{1}{-i(z-\overline{w})}$
in the compact real Riemann surface case is given by
\begin{equation}
\label{theta111}
K_{\zeta}(z,w)
=
\frac{\vartheta[\zeta](\tauBa{w}-z)}{i\vartheta[\zeta](0)E(z,\tauBa{w})}
\end{equation}
(see \cite[Section 2.4]{av3} and, up to conjugation and multiplication by a constant, \cite{av2} and \cite{MR1704479}), 
while the counterpart of the kernel $K_T(z,w)$ is given in \eqref{eqKernPos} below.
Here $\vartheta[\zeta]$ is the theta function with characteristic $\zeta$ and $E(\cdot,\cdot)$ is the prime form
(for more details see Section \ref{sub21} below or \cite{fay1}).
Note that the kernel \eqref{theta111} is not always positive definite on a subset of a given real Riemann surface. 
The case where it is positive definite ($\zeta \in T_0$ and the compact real Riemann surface is of dividing type, see below) corresponds to the Hardy space case.
\smallskip

Our approach to prove de Branges type theorem (and hence Beurling Lax type theorem) is substantially different 
compared to the approaches of the papers mentioned in Section \ref{sucSectionIntoCRS}.
Essentially, all of these papers are using the method of lifting the problem to the universal 
covering and using the classical Beurling theorem to solve the problem.
We, on the other hand, use commutative operators vessels and operator model theory \cite{MR2043236,KLMV,MR1634421} for pairs of commuting nonselfadjoint operators
in order to build explicitly the (inner) function $T$.
The function $T$ (actually a multiplicative multi-valued function, i.e. a mapping of line bundles)
is then the transfer function of a commutative two-operator vessel.
The inner space of this vessel is the given reproducing kernel Hilbert space 
(the orthogonal complement in $H^2$ to the invariant subspace in the case of Beurling Lax theorem)
and the main operators are the compressed multiplication operators by a pair of real 
meromorphic functions on the Riemann surface.
\smallskip

Furthermore, all the papers mentioned in Section \ref{sucSectionIntoCRS} consider the setting of multiplicative functions.
They use a non-canonical approach in the sense it is required to choose a measure on the boundary of the Riemann surface, 
resulting in more complicated calculations. Our approach (as in \cite{av2,av3}) is to consider the half-order multiplicative differentials setting. 
This approach is canonical in the sense that no selection of a measure is required.


\subsection{Structure of the paper}
The paper consists of six sections besides the introduction.
In Section \ref{secPre}, which consists of preliminaries, we review some basic definitions and results related to compact real Riemann surfaces.
In addition, we survey the theory of vessels associated to pairs of commutative nonselfadjoint operators.
Section \ref{secRtOm} is dedicated to present new supplementary results in realization theorem and functional models of two-operator vessels.
\smallskip

The main theorem, namely the de Branges structure theorem, the counterpart of Theorem \ref{1_1} in the compact real Riemann surfaces setting, is presented in Section \ref{secMainThm}.
In Section \ref{secH2}, we specialize our results to the setting of the Hardy space and obtain versions of Beurling's theorem on finite bordered Riemann surfaces.
Section \ref{secMultOp} is dedicated to study the compressed multiplication operators  associated to a real Riemann surface and a pair of meromorphic functions.
To ease the presentation, the proof of de Branges structure theorem is given later in Section \ref{secProofStruct}.
\smallskip

{\it The authors wish to thank the referee for carefully reading the manuscript 
and for raising valuable suggestions to improve the original draft.}


\section{Preliminaries}
\setcounter{equation}{0}
\label{secPre}

In the preliminaries section we review some necessary background and notions related to the content of this paper.
In the first part we survey the compact real Riemann surfaces and their Jacobians. 
The second part is dedicated to the theory of vessels associated to pairs of commutative operators. 
The model space associated to an expansive mapping between certain line bundles defined on a compact real Riemann surface is presented in the last part.

\subsection{Compact real Riemann surfaces}
\label{sub21}

It is a well-known fact (see \cite{livsic2,KLMV} and Section \ref{sub22} in the upcoming pages) that real algebraic curves and compact real Riemann surfaces play an important role in the theory of operators vessels.
A survey of the main needed tools (including the prime form and the Jacobian) can be found in \cite[Section 2]{av3},
the descriptions of the Jacobian variety of a real curve and the real torii is in \cite{vinnikov5}.
For general background, we refer the reader to \cite{fay1,GrHa,gunning2,mumford1,mumford2}.
\smallskip

A compact Riemann surface $X$ of genus $g$ is called {\it real}
if it comes equipped with an anti-holomorphic involution $\tau : p \mapsto \tauBa{p}$ from $X$ into itself.
Let $X_{\mathbb R}$ be the set, assumed nonempty, of real points in $X$ (that is the set of points $p \in X$ such that $p = \tauBa{p}$).
The set $X_{\mathbb R}$ consists of $k \geq 1$ disjoint connected components homeomorphic to circles and denoted by $X_0,\ldots,X_{k-1}$.
Two cases should be distinguished. The first, called the dividing case, is when $X \setminus X_{\mathbb R}$ is not connected; then it is a union of two connected components, $X_{-}$ and $X_{+}$.
The second, the non-dividing case, is when $X \setminus X_{\mathbb R}$ is connected.
In the non-dividing case we set an arbitrary orientation on each of the boundary components,
while in the dividing case we set the orientations according to the choice of $X_+$ 
(i.e. so that $X_\R = X_0 + \ldots + X_{k-1}$ is the positive oriented boundary of $X_+$). 
\smallskip

Repeating the constructions in \cite{ahlfors,av3,vinnikov5}, 
we fix points $p_j \in X_j$ where $j=0,\ldots,k-1$ 
and we choose paths from $p_0$ to $p_j$ for $j=1,\ldots,k-1$,
denoted by $C_j$, which do not contain any other real points.
Then the collection
$A_{g+1-k+j} = X_j$ and $B_{g+1-k+j}= \pm C_j \mp \tauBa{C_j}$ where $j=1,\ldots,k-1$,
can be extended to a canonical basis $A_1,\ldots,A_g$, $B_1,\ldots,B_g$ 
of the homology group $H_1(X,\mathbb{Z})$ 
such that the complex conjugation is given by
$T = \left( \begin{smallmatrix} I& H \\ 0 &-I \end{smallmatrix} \right)$,
where $H$ is a $g \times g$ symmetric matrix depending on the number $k$ of real circles
and whether we are in the dividing case or the non-dividing case, 
see \cite[Equations (2.7-2.8)]{av3}.
Let $\omega_1,\ldots,\omega_g$ be the corresponding normalized basis of the space of holomorphic differential on $X$, 
where $g$ is the genus of $X$, i.e. $\int_{A_i}\omega_j = \delta _{i,j}$.
The period matrix $\Gamma$ is defined by $\Gamma_{ij} = \int_{B_i}{\omega_j}$, and $\Gamma$ is symmetric and satisfies ${\rm Im}\, \Gamma > 0$. 
It can be shown that the Hermitian part of $\Gamma$ is equal to $\frac{1}{2}H$ and
we write $ \Gamma = \frac{1}{2}H + i Y^{-1}$.
\smallskip

One associates to $X$ the Jacobian variety $J(X) = \mathbb{C}^g/\Lambda$,
where $\Lambda$ is the lattice defined by $\Lambda = \mathbb Z^g+\Gamma\mathbb Z^g$.
The Abel-Jacobi mapping from $X$ to $J(X)$ defined by
\[
\mu : p \longrightarrow
\left(
    \begin{array}{c}
      \int_{p_0}^{p}{\omega_1} \\
     \vdots \\
      \int_{p_0}^{p}{\omega_g} \\
          \end{array}
  \right),
\]
is well-defined for an arbitrary base point $p_0 \in X$;
we take $p_0 \in X_0$ as chosen above in the construction of the canonical homology basis.
We note that $\Lambda$ is invariant under complex conjugation and hence the complex conjugation is defined on $J(X)$.
Furthermore, since we chose $p_0 \in X_\R$, it follows that $\tau$ and the complex conjugation 
are equivariant under the Abel-Jacobi map, i.e. $\mu(\tauBa{p}) = \overline{\mu(p)}$.
We extend $\mu$ by linearity to divisors on $X$.
Since by Abel-Jacobi theorem a divisor $D$ is equivalent to zero
if and only if $\mu(D)=0$, we also view $\mu$ as defined on 
line bundles on $X$.
\smallskip

The corresponding theta function is given by
\[
\vartheta(\lambda) =
\sum_{n \in \mathbb{Z}^g}{\exp\bracketsA{i \pi n ^ t \Gamma n + 2 i \pi n ^ t \lambda}},
\]
and is a quasi-periodic function with respect to the lattice $\Lambda$,
that is,
\[
\vartheta (\lambda + m) = \vartheta (\lambda)
\quad {\rm and } \quad
\vartheta (\lambda + \Gamma n) = \exp \bracketsA{ -i \pi n ^ t \Gamma n - 2 \pi i n \lambda} \vartheta (\lambda),
\]
where $n,m \in \mathbb Z ^g$.
Therefore, the theta function defines a divisor in $J(X)$.
The theta function with characteristic $a$ and $b$ in $\mathbb R^g$ is defined by
\[
\vartheta\left[ \!
\begin{array}{c}        a \\       b     \end{array} \!
\right] (\lambda) =
\sum_{n \in \mathbb{Z}^g}
{\exp{ \bracketsA{i \pi (n + a) ^ t \Gamma (n + a) + 2 i \pi (n + a) ^ t (\lambda+b)}}}.
\]
In this paper (as in \cite{av2,av3}), we consider the framework of multiplicative half-order differentials.
In order to construct and define the half-order differentials, we choose an atlas $\bracketsB{V_j , z_j}{j\in \mathcal J}$ on $X$,
for which every nonempty intersection is assumed to be simply connected.
Then there exists a family of analytic square-roots (see \cite{hawley-schiffer}), $\bracketsB{\sqrt{dz_j/dz_i}}{i,j\in \mathcal J}$, such that the following cocycle condition
\begin{equation}
\label{cocycleCond}
\sqrt{\frac{dz_i}{dz_j}}
=
\sqrt{\frac{dz_i}{dz_k}}
\sqrt{\frac{dz_k}{dz_j}}
\end{equation}
holds on $V_j \cap V_i \cap V_k $ (whenever the intersection is not empty).
Among the line bundles defined by \eqref{cocycleCond}, 
the line bundle corresponding to $-\kappa \in J(X)$ (where $\kappa$ is the Riemann constant) 
plays an important role and is denoted by $\Delta$. 
The sections of $\Delta$ are referred to as half-order differentials.
\smallskip

For $\zeta \in J(X)$ we define the corresponding multipliers over the cycles $A_j$ and $B_j$ by
\[
\chi (A_j) = \exp \bracketsA{-2 \pi i a_j},
\qquad
\chi (B_j) = \exp \bracketsA{2 \pi i b_j}
\]
where $j=1,\ldots,g$ and $\zeta = b+ \Gamma a$ where $a,b \in \R^g$.
The corresponding flat unitary line bundle is denoted by $L_\zeta$,
notice that  $\mu(L_\zeta)=\zeta$.
\smallskip

A multiplicative half-order differential corresponding to $\zeta$
is a section of $L_\zeta \otimes \Delta$.
Thus, a multiplicative half-order differential
is a family of functions $\bracketsB{f_j}{j \in \mathcal J}$ 
defined on the pre-images $\widetilde{V}_j$ of $V_j$, $j \in \mathcal J$, on the universal covering $\widetilde{X}$ of $X$, satisfying
\[
f_i(\widetilde{u}) 
= 
\sqrt{\frac{dz_j}{dz_i}}f_j(\widetilde{u}), 
\qquad 
\widetilde{u} \in \widetilde{V}_i \cap \widetilde{V}_j.
\]
and
\[
f_j(\widetilde{u}_2) 
= 
f_j(\widetilde{u}_1) \exp \bracketsA{2 \pi i (b^tm - a^t n)},
\]
where $\widetilde{u}_1,\widetilde{u}_2$ are points on the universal covering such that 
$\widetilde{\mu}(u_2)-\widetilde{\mu}(u_1) = n + \Gamma m$
($\widetilde{\mu}: \widetilde{X} \rightarrow \mathbb C ^g$ is the lifting of $\mu$ to the universal coverings and $n,m \in \Z$).
We often abuse the notation and view $f_j$ as multi-valued function on $V_j$.
We always assume that $\vartheta(\zeta) \neq 0$ which is equivalent to the fact that
$L_\zeta \otimes \Delta$ does not have global holomorphic sections
i.e. there is no non-trivial multiplicative half-order differential which is globally holomorphic.
\smallskip

The Cauchy kernel plays an important role in this framework as it is the analogue in the compact real Riemann surface to
of the kernel $\frac{1}{-i(z-\overline{w})}$ in the upper half--plane.
In the line bundle case, unlike the vector bundle case
(for explicit formulas for the Cauchy kernel in the vector valued case in genus one, see \cite{MR3584679}),
the Cauchy kernel can be described explicitly by (see \cite{av3})
\begin{equation}
\label{Cauchy}
K_{\zeta}(u,v)
=
\frac
{\vartheta[\zeta](\tauBa{v}-u)}
{i\vartheta[\zeta](0)E(u,\tauBa{v})}.
\end{equation}
Here we identify the points $u,\tauBa{v} \in X$ with the images $\mu(u)$ and $\overline{\mu(v)}$ in $J(X)$.
The prime form $E(u,v)$ is a multiplicative differential of order $-\frac{1}{2}$ in each of the variables $u$ and $v$.
It is defined by 
\[
E(u,v) = \frac{\vartheta[\delta](v-u)}{\sqrt{\xi_\delta(u)}\sqrt{\xi_\delta(v)}}.
\]
Here $\delta$ is a non-singular odd half-period and
$\xi_\delta$ is defined by 
$(\sqrt{\xi^2_\delta(u)})^2 = \sum_{j=1}^g{\frac{\partial \vartheta [ \delta ]}{\partial z_j}  (0)\omega_j(u)    } $,
for more details see \cite{fay1}.
Its main property is that $E(u,v) = 0$ if and only if $u = v$,
and thus it can be considered as the analogue for the compact Riemann surface case of the difference between two numbers in $\mathbb C$.
In terms of local coordinates, the prime form satisfies
\[
E(u,v) = (t(v)-t(u)) + o((t(v)-t(u))^2).
\]
For more details on the prime form, we refer to \cite[section 2.3]{av3} and \cite{fay1,mumford2}. \smallskip

The Cauchy kernel is Hermitian (see \cite[Proposition 2.8]{av3} and \cite{vinnikov5}) 
as long as $\zeta$ belongs to a disjoint union of the $g$-dimensional real torii, given by:
\begin{align}
    \nonumber
T_{\nu}
=
\big\{
\zeta \in J(X) : \zeta 
= &
\frac{1}{4} \text{diag}(H) +
\frac{\nu_1}{2}e_{g-k+2}
+
\ldots
+
\frac{\nu_{k-1}}{2}e_{g}
+
\\
&
i a_1 \imag \Gamma_1 + \cdots + i a_g \imag \Gamma_g
\big\},
\label{eqRealTorii}
\end{align}
where $\Gamma_1,\ldots,\Gamma_g$ are the columns of the period matrix $\Gamma$,
$a_1,\ldots,a_{g-k+1} \in \mathbb R / 2 \mathbb Z$,
$a_{g-k+2},\ldots,a_{g} \in \mathbb R / \mathbb Z$
and $\nu \in \{ 0,1\}^{k-1}$.
Furthermore, if $X$ is of dividing type and $\zeta \in T_0$ (\cite[Theorem 2.1]{av3}),
then automatically $\vartheta(\zeta) \neq 0$ and
the kernel \eqref{Cauchy} is positive on $X_+$ and negative on $X_-$
(since $\pm K_{\zeta}(\cdot,v)$ is the reproducing kernel for the Hardy spaces $H^2(L_\zeta \otimes \Delta , X_{\pm})$, 
see also Section \ref{secH2} below).
\smallskip

The Cauchy kernels satisfy an important identity, which is referred as the \textit{collection formula} (see \cite{av2} and \cite{av3}) and is used repeatedly in the sequel.
First, it is convenient to define, using the notations from \cite{av2}, the following matrices.
\begin{definition}
\label{defCollFor}
Let $y$ be a meromorphic function on $X$ of degree $n$ with simple poles $(p^{(j)})_{j=1}^{n}$ 
and residues $(-c_j)_{j=1}^{n}$ given in terms of some fixed local coordinates at the poles.
Then for $\lambda_1,\lambda_2 \in \mathbb{C}$, we set
\begin{align}
\label{formulaK1}
\mathbb{K}(\lambda_1,\lambda_2) 
& =
(\lambda_1 - \lambda_2)
\left(
\frac{1}{\sqrt{dy}(u^{(i)})}
\frac
{\vartheta[\zeta](v^{(j)}-u^{(i)})}
{\vartheta[\zeta](0)E(v^{(j)},u^{(i)})}
\frac{1}{\sqrt{dy}(v^{(j)})}
\right)_{i,j=1}^n
\\
\label{formulaK2BLA}
\mathbb{K}(\lambda_1,\infty) 
& = 
-
\left(
\frac{1}{\sqrt{dy}(u^{(i)})} 
\frac
{\vartheta[\zeta](p^{(j)}-u^{(i)})}
{\vartheta[\zeta](0)E(p^{(j)},u^{(i)})}
\frac{\sqrt{c_j}}{\sqrt{dt_j}(p^{(j)})}
\right)_{i,j=1}^n
\\
\label{formulaK3}
\mathbb{K}(\infty,\lambda_1) 
& =
\left(
\frac{\sqrt{c_i}}{\sqrt{dt_i}(p^{(i)})}
\frac
{\vartheta[\zeta](u^{(j)}-p^{(i)})}
{\vartheta[\zeta](0)E(u^{(j)},p^{(i)})}
\frac{1}{\sqrt{dy}(u^{(j)})}
\right)_{i,j=1}^{n}
\end{align}
where $(u^{(i)})_{j=1}^{n}$ and $(v^{(j)})_{j=1}^{n}$ are the $n$ pre-images in $X$ of
$\lambda_1$ and $\lambda_2$, respectively, assumed to be all distinct 
(it is possible to extend the kernel by continuity to ramified fibers,
but we will not consider this case in the sequel).
We set $\mathbb{K}(\lambda_1, \lambda_2)=I$ for $\lambda_1 = \lambda_2$ by continuity.
\end{definition}
The matrices \refMult{formulaK1}{formulaK3} satisfy the following relations, also known as the
collection formulas (\cite[Section 4]{av2}):
\begin{align}
\label{CollFormLimitVer}
\mathbb{K}(\lambda_1,\infty) \mathbb{K}(\infty,\lambda_2) & =  \mathbb{K}(\lambda_1,\lambda_2)
\\
\label{CollFormLimitVer2BLA}
\mathbb{K}(\lambda_1,\lambda_3) \mathbb{K}(\lambda_3,\lambda_2) & =  \mathbb{K}(\lambda_1,\lambda_2)
\\
\label{CollFormLimitVer3BLA}
\mathbb{K}(\lambda,\lambda) & = I
.
\end{align}
An additional version of the collection formula can be found in \cite[Lemma 4.1]{av3}.


\subsection{Commutative Vessels in Hilbert space}
\label{sub22}
It is a well-known fact that the best way to study a bounded nonselfadjoint operator is to view it 
as an element of an underlying \textit{colligation} rather than studying directly the operator itself. 
There is a deep connection
between invariant subspaces of such an operator and factorizations of the characteristic function of a colligation (see \cite{MR20:7221}).
As soon as we consider several commuting nonselfadjoint operators, the colligation does not carry enough structure to reflect the interaction between the operators.
It seems, see \cite{livsic1,KLMV,vinnikov4,MR1634421}, that the appropriate framework to study several commuting nonselfadjoint operators is via the notion of a
commutative operator vessel. In this paper we consider the case of a pair of commuting operators, 
although a generalization to the case of an $n$-tuple of commuting operators 
(and even beyond to the case of a representation of a finite dimensional real Lie algebra)
does exist, see \cite{KLMV,MR3609238,MR3906269}.
\smallskip

A commutative two-operator vessel is a collection
\begin{equation*}
\mathcal{V} = 
\bracketsA{A_1 \, , \, A_2 \, ; \, H \, , \, \Phi \, ,  \, E \, ; \,  \sigma_1 \, , \, \sigma_2 \, , \, \gamma \, , \, \widetilde{\gamma}}
\end{equation*}
where $H$ ("the inner space")  and $ E$ ("the external space") are Hilbert spaces, and ${\rm dim}\,  E < \infty$.
The operators $A_1$ and $A_2$ are bounded in $H$, commute ($A_1A_2=A_2A_1$) and with finite
non-Hermitian rank, i.e.
\begin{equation}
\label{collCond}
\frac{1}{i} \left(A_k - A_k^* \right) = \Phi^* \sigma_k \Phi,\quad k=1,2,
\end{equation}
where $\sigma_1$ and $\sigma_2$ are selfadjoint operators in $E$.
Finally, $\gamma$ and $\widetilde{\gamma}$ are bounded selfadjoint operators on $E$
which satisfy the {\it vessel conditions}:
\begin{align}
\label{VesselCond1}
\sigma_1 \Phi A_2^* - \sigma_2 \Phi A_1^* & = \gamma \Phi,
\\
\label{VesselCond2BLA} 
\sigma_1 \Phi A_2 - \sigma_2 \Phi A_1 & = \widetilde{\gamma} \Phi,
\\
\label{VesselCond3}
i(\sigma_1 \Phi \Phi^* \sigma_2 - \sigma_2 \Phi \Phi^* \sigma_1)
& =
\widetilde{\gamma} - \gamma.
\end{align}
Equations \refMult{VesselCond1}{VesselCond3} are called the input, output and linkage vessel conditions, respectively.
It is easily seen that if the colligation condition \eqref{collCond} and the output vessel condition \eqref{VesselCond2BLA} are satisfied
and we define $\gamma$ by the linkage vessel condition \eqref{VesselCond3}, 
then, the input vessel condition \eqref{VesselCond1} holds.
\smallskip

The \textit{complete characteristic function} (CCF) of a vessel is defined by \cite[Section 3.4]{KLMV},
\begin{equation}
\label{CCF}
W(\xi_1,\xi_2,z) = I -i \Phi(\xi_1 A_1 + \xi_2 A_2 -zI)^{-1} \Phi ^* (\xi_1 \sigma _1 + \xi_2 \sigma_2)
\end{equation}
and it satisfies the following metric properties ($\xi_1,\xi_2 \in \R$):
\begin{align}
\label{eqWinA}
W(\xi_1,\xi_2,z)
\left( \xi \sigma \right)
W(\xi_1,\xi_2,z)^*
 & =
\left( \xi \sigma \right)
\text{ where {\rm Im}($z$)}=0
\\
W(\xi_1,\xi_2,z)
\left( \xi \sigma \right)
W(\xi_1,\xi_2,z)^*
 & \geq 
\left( \xi \sigma \right)
\text{ where {\rm Im}($z$)}>0,
\label{eqWinB}
\end{align}
where $\left( \xi \sigma \right)$ stands for $\xi_1\sigma_1 +\xi_2 \sigma_2$.
The complete characteristic function, for fixed $\xi_1,\xi_2 \in \C$, 
is analytic for all $z$ lie outside the spectrum of $\xi_1 A_1 + \xi_2 A_2$.
\smallskip

Given a commutative vessel, we define a polynomial (that we assume to be not identically zero) of two complex variables, 
called the \textit{(affine) discriminant polynomial}, by
\[
p(\lambda_1,\lambda_2) = \det (\lambda_1 \sigma_2 - \lambda_2 \sigma_1 + \gamma).
\]
The associated real (affine) plane curve $C_0$ is called the \textit{(affine) discriminant curve} associated to $\mathcal{V}$.
Writing the discriminant polynomial in homogeneous form
\[
P(\nu_0,\nu_1,\nu_2) = \det (\nu_1 \sigma_2 - \nu_2 \sigma_1 + \nu_0 \gamma),
\]
leads to a corresponding real algebraic curve $C$ in $\mathbb P ^2$. 
We always assume that $\det{(\xi_1 \sigma_1 + \xi_2 \sigma_2)}$ is not identically zero,
so $C$ is the projective closure of $C_0$.
\smallskip

In order to continue, we set the principal subspace $\widehat{ H} \subseteq H$ of a vessel to be
\[
\widehat{ H} =
\bigvee _{m_1,m_2 = 0} ^{\infty}{A_1^{m_1}A_2^{m_2}\Phi^*(E)}=
\bigvee _{m_1,m_2 = 0} ^{\infty}{A_1^{*m_1}A_2^{*m_2}\Phi^*(E)},
\]
and we say that a vessel is {\it irreducible} whenever $\widehat{H}=H$.
\smallskip

Two vessels
\[
\bracketsA{A_1 ,A_2; H ,\Phi ,E;  \sigma_1 ,\sigma_2 , \gamma , \widetilde{\gamma}}
\quad {\rm and} \quad
\bracketsA{A_1^\prime ,A_2^\prime ; H^\prime ,\Phi^\prime , E;  \sigma_1 ,\sigma_2 , \gamma , \widetilde{\gamma}}
\]
are said to be {\it unitary equivalent} if they share common external data
(that is, both have a common external space $E$ and matrices $\sigma_1,\sigma_2, \gamma $ and $ \widetilde{\gamma}$)
and there exists a unitary operator $U$ from 
${ H}$ to ${ H} ^\prime$ such that
\begin{equation}
\label{eqVesselEq}
A_1^\prime   U =  U A_1  ,
\quad
A_2^\prime  U =  U A_2 
\quad
{\rm and }
\quad
\Phi^\prime U = \Phi .
\end{equation}

The following theorems were established by Liv\v{s}ic (see for instance \cite{livsic1,KLMV}).

\begin{theorem}[{\cite[Theorem 2]{livsic3}}]
\label{polyVanish}
The polynomial $ p(A_1,A_2)$ vanishes on the principal subspace of the colligation
$\widehat{H}$.
\end{theorem}

As a consequence, it follows through the spectral mapping theorem 
that the joint spectrum of $A_1$ and $A_2$ restricted to $\widehat{H}$ lies on the curve $C_0$.
We recall that there are several definitions for the notion of the joint spectrum of a pair of commuting operators.
Here we mention the definition of Harte \cite{harte1972spectral} 
(Harte's spectrum consists all points $(\lambda_1,\lambda_2)\in \C^2$ such that 
either the map $x \rightarrow \begin{psmallmatrix}	(A_1 - \lambda_1 I) x \\ (A_2 - \lambda_2 I)x \end{psmallmatrix}$ is not left invertible 
or the map $\begin{psmallmatrix} x_1 \\ x_2 \end{psmallmatrix} \rightarrow (A_1 - \lambda_1 I) x_1 + (A_2 - \lambda_2 I)x_2 $ is not right invertible) 
and Taylor \cite{MR42:3603}.
These various definitions coincide when the operators have compact imaginary parts, and in particular in our setting 
(see \cite{MR1320541}, where it also shown that in this case $(\lambda_1,\lambda_2)$ is in the joint spectrum if and only if
$\xi_1\lambda_1+\xi_2\lambda_2$ is in the spectrum of $\xi_1A_1+\xi_2A_2$ for all $\xi_1,\xi_2\in{\mathbb C}$).
\smallskip

The discriminant polynomial can be also described in terms of the selfadjoint matrices $\sigma_1,\sigma_2$ and $\widetilde{\gamma}$.
\begin{theorem}[{\cite[Corollary 1]{livsic2}}]
The following equality holds:
\[
\det (\lambda_1 \sigma_2 - \lambda_2 \sigma_1 + \gamma)
=
\det (\lambda_1 \sigma_2 - \lambda_2 \sigma_1 + \widetilde{\gamma}).
\]
\end{theorem}
Hence
$\nu_1 \sigma_2 - \nu_2 \sigma_1 + \nu_0\gamma$
and
$\nu_1 \sigma_2 - \nu_2 \sigma_1 + \nu_0\widetilde{\gamma}$
are (the input and output, respectively) determinantal representations of the discriminant curve $C$.
Then for each $\nu = (\nu_0,\nu_1,\nu_2)\in C$ we define non-trivial subspaces of $E$ by:
\begin{align*}
\mathcal{E}(\nu) & = \ker (\nu_1 \sigma_2 - \nu_2 \sigma_1 +\nu_0\gamma)\\
\mathcal{\widetilde{E}}(\nu) & = \ker (\nu_1 \sigma_2 - \nu_2 \sigma_1 +\nu_0\widetilde{\gamma}).
\end{align*}
The complete characteristic function is a function of two independent variables
and does not admit a satisfactory factorization theory as in the single-operator colligation case.
However, we can restrict the CCF to the discriminant curve $C$ and to the family of subspaces $\mathcal{E}(\lambda)$.
More precisely, the \textit{joint characteristic function} (JCF) \cite[Section 10.3]{KLMV},
\begin{equation}
\label{eqJcf}
S(\lambda) 
= \restr{W(\xi_1,\xi_2,\xi_1 \lambda_1 + \xi_2 \lambda_2)}{\mathcal E (\lambda)},
\end{equation}
where $\lambda = (\lambda_1,\lambda_2)\in C_0$, 
is a map between the subspaces $\mathcal{E}(\lambda)$ and $\widetilde{\mathcal{E}}(\lambda)$. 
Furthermore, it is independent of the choice of $\xi_1,\xi_2 \in \C$ when
$\xi_1 \lambda_1 + \xi_2 \lambda_2$  does not belong to the spectrum of $\xi_1 A_1+ \xi_2 A_2$.
As a system theory interpretation, the joint characteristic function determines the input-output relation 
of the corresponding (overdetermined)two-dimensional system by $y_0 = S(\lambda) u_0$, where $y_0 \in \widetilde{\mathcal{E}}(\lambda)$ 
and $u_0 \in \mathcal{E}(\lambda)$ are the amplitudes of the double periodic wave functions 
with frequency $(\lambda_1,\lambda_2) \in C_0$ (see \cite{MR97m:30051,MR1634421}).
\smallskip

In order to continue, we assume that the discriminant polynomial $P(\nu)$, defining the discriminant curve $C$, 
is an irreducible homogeneous polynomial and we denote by $X$ the desingularizing Riemann surface of $C$.
It can be shown that the $\dim \widetilde{\mathcal E}(\nu) \leq s$ where $s$ is the multiplicity of $\nu$ on $C$, see \cite[Prop. 10.5.1]{KLMV}.
Therefore $\dim \widetilde{\mathcal E}(\nu)=\dim{\mathcal E}(\nu)=1$ for all non-singular $\nu \in C$.
Therefore, in this paper, as opposed to a more general vector bundle setting, we consider the line bundle setting.
Hence, we may use the theory of theta functions and Jacobian varieties 
to express concretely the Cauchy kernels, see Section \ref{sub21} above
(for more details, see \cite{MR97m:30051,MR1704479}).
We assume furthermore that the input and the output determinantal representations are maximal and fully saturated 
as we define in the next several paragraphs allowing us to consider line bundles on the 
desingularizing Riemann surface rather than line bundles on the singular algebraic curve,
see \cite{MR2964618} in setting of determinantal representation
and \cite{MR918564,MR0024985} in general.
\smallskip

We say that the output determinantal representation 
$\widetilde{U}(\nu)= \nu_1 \sigma_2 - \nu_2 \sigma_1 + \nu_0 \widetilde{\gamma}$ 
is maximal if any point $\nu \in C$ satisfies $\dim \widetilde {\mathcal E} (\nu) = s$.
A stronger condition is as follows (see \cite{MR2962792}).
We define $\widetilde{V}(\nu) = \operatorname{adj}[\widetilde{U}(\nu)]$,
this is a matrix of homogeneous polynomial of degree $m-1$ where $m = \deg{(P)} = \dim (E)$.
Then we say that the output determinantal representation is {\it fully saturated}  
if all the entries of the adjoint matrix $\widetilde{V}$ vanish on the adjoint divisor:
$(\widetilde{V}_{i,j}) \geq D_{\rm sing}$ for all $i, j = 1, . . . , m$.
Here the adjoint divisor or the divisor of singularities is given by 
$D_{\rm sing} = (m-3)(\nu_0) - (\omega)$, 
where $\omega$ is the meromorphic differential defined by
\begin{equation}
\label{eqOmegaDef}
\omega 
= 
\frac{dy_1}{\sfrac{\partial p}{\partial \lambda_2}} 
=
-\frac{dy_2}{\sfrac{\partial p}{\partial \lambda_1}},
\end{equation}
i.e. $\omega$ has poles at the singular points of $C$ and zeros of order $m-3$ at infinity
(we denote by $\lambda_1$, $\lambda_2$ the affine coordinates
and by $y_1$, $y_2$ the corresponding meromorphic functions
on the desingularizing Riemann surface $X$).
Every fully saturated determinantal representation is maximal 
and the converse is true if at each singular point of $C$,
there are no two distinct branches with the same tangent.
As we will see below, under certain conditions, maximality (or full saturation)
of the output determinantal representation implies the maximality (or full saturation) 
of the input determinantal representation (see Lemma \ref{lemMaxIO}).
\smallskip

As follows from the theory of determinantal representations (see \cite{MR97m:30051,MR2043236,MR1634421}),
the pair of families $\mathcal E$ and $\widetilde{\mathcal{E}}$ 
of one-dimensional vector spaces $\mathcal E(\nu)$ and $\widetilde{\mathcal E}(\nu)$ 
defined on the subset of non-singular points of $C$ can be lifted to a pair of line bundles, 
also denoted by $\mathcal E$ and $\widetilde{\mathcal{E}}$ on $X$, the normalization of $C$.
We note that up to a sign, selfadjoint determinantal representations are uniquely defined up to Hermitian equivalence:
$\nu_1 \sigma_2 - \nu_2 \sigma_1 + \nu_0 \widetilde{\gamma} \mapsto
\rho^* \left(\nu_1 \sigma_2^\prime - \nu_2 \sigma_1^\prime + \nu_0 \widetilde{\gamma}^\prime \right) \rho$,
$\rho$ is non-singular, by their corresponding line bundles.
Furthermore, there exist isomorphisms (up to some twists) between the 
kernel bundles $\mathcal E$ and $\widetilde{\mathcal E}$ 
and line bundles of the form $L_{\zeta} \otimes \Delta$ and $L_{\widetilde{\zeta}} \otimes \Delta$
where $\z,\widetilde{\z}\in T_\nu$ for some $\nu$, see \eqref{eqRealTorii},
and $\vartheta(\zeta) \neq 0$ and $\vartheta(\widetilde{\zeta}) \neq 0$.
More precisely,
$\mathcal E \otimes \mathcal O (m-2)(-D_{\rm sing}) \cong L_{\zeta} \otimes \Delta$ 
and
$\widetilde{\mathcal{E}} \otimes \mathcal O (m-2)(-D_{\rm sing}) \cong L_{\widetilde{\zeta}} \otimes \Delta$
.
\smallskip

These isomorphisms can be given explicitly in terms of normalized sections of $\mathcal E$ and $\widetilde{\mathcal E}$, 
denoted by ${\bf u}^{\times}$ and $\widetilde{\bf u}^{\times} $, respectively, 
as follows (see \cite[Equation 2-28]{MR1634421}):
\begin{equation}
\label{eqIsoExplicit}
f(p) \mapsto \frac{1}{\omega(p)} f(p) \widetilde{\bf u}^\times(p),
\end{equation}
similarly for ${\mathcal{E}}$.
A normalized section $\widetilde{\bf u}^{\times}$ is the unique, up to a constant, 
meromorphic section of $\widetilde{\E} \otimes L_{-\widetilde{\zeta}} \otimes \Delta$ with poles at the points of $C$ at infinity.
A normalized section can be explicitly described by the Cauchy kernels 
(and their derivatives in the case where the poles are not simple).
More explicitly, ${\bf u}^{\times}(p)$ and $\widetilde{\bf u}^{\times}(p) $ 
are column-vectors containing the Cauchy kernels (and possibly their derivatives) 
evaluated at $p$ and the $m$ poles of the pair of coordinate meromorphic functions $y_1$ and $y_2$ on $X$,
see \eqref{normSectionCK} for the canonical determinantal representation corresponding to $\zeta$.
We refer the reader to \cite{MR1704479} for more details.

The joint characteristic function lifts to a holomorphic mapping
of the line bundles ${\mathcal E}$ and $\widetilde{\mathcal E}$ on $X$
on the complement of the points lying above the joint spectrum
of $A_1$ and $A_2$, in particular,
in a neighborhood of the points of $C$ at infinity.
Under the above isomorphisms, the joint characteristic function is translated into a mapping between flat unitary line bundles.
This mapping is called the \textit{Normalized joint characteristic function} (NJCF) 
and is denoted in the sequel by $T(p)$.
An explicit relation between the two characteristic functions is therefore given by (see \cite[Section 6]{MR1704479}):
\begin{equation}
\label{st}
S(p) {\bf u}^{\times}(p) = \widetilde{\bf u}^{\times}(p) T(p).
\end{equation}

We can consider also the left kernel families
\begin{align*}
\mathcal{E}_l(\nu) & = \ker_l (\nu_1 \sigma_2 - \nu_2 \sigma_1 +\nu_0\gamma)\\
\mathcal{\widetilde{E}}_l(\nu) & = \ker_l (\nu_1 \sigma_2 - \nu_2 \sigma_1 +\nu_0\widetilde{\gamma}),
\end{align*}
and their liftings, also denoted by $\mathcal{E}_l$ and $\mathcal{\widetilde{E}}_l$, to line bundles on $X$.
There exists an isomorphism between the dual bundles of 
$\mathcal E$ and $\widetilde{\mathcal E}$ and $\mathcal E_l$ and $\widetilde{\mathcal E}_l$
up to an appropriate twist.
More explicitly, 
$\widetilde{\mathcal{E}}_l \otimes \mathcal O (m-2)(-D_{\rm sing})$ is the Serre dual of 
$\widetilde{\mathcal{E}} \otimes \mathcal O (m-2)(-D_{\rm sing})$
(i.e. the dual tensored with the canonical bundle $K = \Delta \otimes \Delta$),
so that $\widetilde{\mathcal{E}}_l \otimes \mathcal O (m-2)(-D_{\rm sing})$
is isomorphic to $L_{-\widetilde{\zeta}} \otimes \Delta $.
The corresponding left normalized section $\widetilde{\bf u}_l^{\times}$ is the unique, up to a constant, 
meromorphic section of $\widetilde{\mathcal{E}}_l \otimes L_{\widetilde{\zeta}} \otimes \Delta$ with poles at the points of $C$ at infinity,
similarly for $\mathcal{E}_l$.
We note that the duality between the left and right kernels
can be given by the inner product with respect to $\xi\sigma = \xi_1\sigma_1 + \xi_2\sigma_2$,
so that given $\widetilde{\bf u}^{\times}$, $\widetilde{\bf u}_l^{\times}$ can be determined uniquely by
\begin{equation}
\label{eqMatNormSecRelation}
\frac{{\bf \widetilde{u}}_l^{\times} (p) (\xi \sigma) {\bf \widetilde{u}}^{\times}(q)}{ \xi y(p) - \xi y(q)}  
= 
\frac{\vartheta[\widetilde{\zeta}](q-p)}{ \vartheta[\widetilde{\zeta}](0)E(q,p)},
\quad
{\bf \widetilde{u}}_l^{\times} (p) (\xi \sigma) {\bf \widetilde{u}}^{\times}(p) = \xi dy(p),
\end{equation}
where $\xi y = \xi_1 y_1 + \xi_2 y_2$ and $\xi dy = \xi_1 dy_1 + \xi_2 dy_2$.
On the other hand, since our determinantal representations are Hermitian, 
if ${\bf \widetilde{u}}^{\times}(p)$ is a normalized section, then ${\bf \widetilde{u}}^{\times}(\tauBa{p})^*$ is a left normalized section.
It follows from \eqref{eqMatNormSecRelation} that (after multiplying ${\bf \widetilde{u}}^{\times}$ by a constant)
${\bf \widetilde{u}}^{\times}_l(p) = \e {\bf \widetilde{u}}^{\times}(\tauBa{p})^*$ where $\e = \pm 1$ is the 
sign of the determinantal representation alluded to above.
It is easily seen that the input and the output determinantal representations have the same sign
and we will always assume in the sequel, unless the contrary is explicitly stated, 
that the determinantal representations have sign $+1$.

\section{Realization theorem and functional models}
\label{secRtOm}

\subsection{Realization theorem}
The following fundamental realization question arises, whether for a given line bundle mapping $T$ there exists a commutative vessel,
such that $T$ is its NJCF. 
A solution to this question is given, in the line bundle setting, by the following realization theorem.
\begin{theorem}[{\cite[Theorem 2.3]{MR1634421}}, {\cite[Theorem 10.5.7]{KLMV}} and {\cite[Theorem 3.3]{MR97m:30051}}]
\label{realizationTh}
A multiplicative function $T(p)$ on $X$ with multipliers corresponding to $\widetilde{\zeta} - \zeta$,
i.e. a mapping between flat unitary line bundles $L_\zeta$ and $L_{\widetilde{\zeta}}$ where $\zeta,\widetilde{\zeta} \in T_\nu$,
is the normalized joint characteristic function of a
commutative two-operator vessel with discriminant polynomial $p(\lambda_1,\lambda_2)$ and maximal input and output determinantal representations corresponding to 
$\zeta,\widetilde{\zeta} \in J(X)$
if and only if
$T(p)$ is holomorphic at the points of $C$ at infinity, meromorphic on $X \setminus X_{\mathbb R}$ and satisfies
$ T(p)\overline{T(\tauBa{p})}=1$ and the kernel
\begin{equation}
\label{eqKernPos}
T(p)
\frac{\vartheta[\zeta](\tauBa{q}-p)}{i \vartheta[\zeta](0)E(p,\tauBa{q})}
\overline{T(q)}
-
\frac{\vartheta[\widetilde{\zeta}](\tauBa{q}-p)}{i \vartheta[\widetilde{\zeta}](0)E(p,\tauBa{q})}
,
\end{equation}
is positive definite on the domain of analyticity of $T$.
\end{theorem}

A multiplicative function $T$ on $X$ which is meromorphic on $X \setminus X_\R$,
satisfying $ T(p)\overline{T(\tauBa{p})}=1$ and such that the kernel \eqref{eqKernPos} is positive 
is called $(\zeta,\widetilde{\zeta})$-expansive.

\begin{proposition}
\label{propCont}
Let $X$ be a compact Riemann surface of dividing type and let $\z,\widetilde{\z} \in T_0$.
Then $T$ is $(\z,\widetilde{\z})$-expansive 
if and only if $T$ is expansive on $X_+$ (i.e $\abs{T(p)} \geq 1$ there).
\end{proposition}
\begin{pf}
Assume that the kernel \eqref{eqKernPos} is positive and let us choose $p=q\in X_+$.
As already mentioned following \eqref{eqRealTorii}, since $\zeta$ and  $\widetilde{\zeta}$ belong to $T_0$ and $X$ is of dividing type,
the Cauchy kernels are positive on $X_+$ and \eqref{eqKernPos} implies
\begin{equation}
\label{eqYtrp}
\frac{\vartheta[\zeta](0)}{\vartheta[\widetilde{\zeta}](0)}
\frac{\vartheta[\widetilde{\zeta}](\tauBa{p}-p)}{\vartheta[\zeta](\tauBa{p}-p)}
\leq
\abs{T(p)}^2,
\end{equation}
for $p \in X_+$.
Consider the limit on the left hand side as $p \rightarrow q_0 \in X_\R$. 
Assuming $q_0 \in X_k$, $\tauBa{q_0}-q_0$ is equal to the $k$-th column of $\Gamma$. 
Since $\z,\widetilde{\z} \in T_0$ and by the properties of the theta function, the limit of the left hand side is equal to one.
Hence, for $\e>0$, a point $p$ in $X_+$ sufficiently close to the boundary $X_\R$ satisfies
$\abs{T(p)}^2 \geq 1 +\e$.
By reflection, a point $p \in  X_-$ sufficiently close to $X_\R$ satisfies $\abs{T(p)}^2 \leq \frac{1}{1 +\e}$.
On the other hand, since the Cauchy kernels are negative on $X_-$, it follows similarly to \eqref{eqYtrp} that
\[
\frac{\vartheta[\zeta](0)}{\vartheta[\widetilde{\zeta}](0)}
\frac{\vartheta[\widetilde{\zeta}](\tauBa{p}-p)}{\vartheta[\zeta](\tauBa{p}-p)}
\geq
\abs{T(p)}^2,
\]
for $p \in X_-$. In particular, this implies that $T$ is analytic on $X_-$.
It then remains to apply the maximum principle to conclude that $\abs{T(p)} \leq 1$ for all $p \in X_-$
and therefore by reflection $\abs{T(p)} \geq 1$ for all $p \in X_+$.
\smallskip

Conversely, let us assume that $\abs{T(p)} \geq 1$ for $p \in X_+$. 
Then, by reflection, $\abs{T(p)} \leq 1$ for $p \in X_-$.
In particular, $T$ is holomorphic in $X_-$ and the multiplication by $T$ is a contraction 
from the Hardy space $H^2(L_\zeta \otimes \Delta , X_{-})$ to the Hardy space 
$H^2(L_{\widetilde{\zeta}} \otimes \Delta , X_{-})$.
Since $-K_\zeta(p,q)$ and $-K_{\widetilde{\zeta}}(p,q)$ are the reproducing kernels of 
$H^2(L_{{\zeta}} \otimes \Delta , X_{-})$ and $H^2(L_{\widetilde{\zeta}} \otimes \Delta , X_{-})$, respectively, 
the kernel \eqref{eqKernPos} is positive on $X_-$, see for instance \cite{Dym_CBMS} for the classical case,
and by reflection \eqref{eqKernPos} is positive on $X_+$.
If the domain $\Omega$ of analyticity of $T$ is connected then the positivity 
of \eqref{eqKernPos} on $\Omega$ follows from \cite[Theorem 1.1.4]{adrs}. 
In the general case we argue as follows.
\smallskip

We choose a pair of dividing functions $y_1$ and $y_2$ on $X$ (see Section \ref{secH2}) 
and two selfadjoint determinantal representations that correspond to $\zeta$ and to $\widetilde{\zeta}$, respectively.
Then one may build a function $W(\xi_1,\xi_2,z)$ by the restoration formula, 
see \cite[Chapter 10]{KLMV} and \cite{MR1634421},
so that 
$S(\lambda) = \restr{W(\xi_1,\xi_2,\xi_1\lambda_1+\xi_2\lambda_2)}{\mathcal{E}(\lambda)}$ 
and $S$ and $T$ are related by \eqref{st}.
We note that for any $\xi_1,\xi_2 >0$ and for any $\imagg{z}>0$
all the points that satisfy $\xi_1 y_1(p)+ \xi_2 y_2(p) = z$ are in $X_+$,
and hence the restoration formula implies that $W(\xi_1,\xi_2,z)$ 
satisfies the metric properties in \eqref{eqWinA} and \eqref{eqWinB}.
Then it follows from the classical theory that 
$\frac{W(\xi_1,\xi_2,z) (\xi \sigma) W(\xi_1,\xi_2,w)^* - (\xi \sigma)}{-i(z- \overline{w})}$ 
is a positive kernel on $\C \setminus \R$, see for instance \cite{Dym_CBMS}.
Finally, since $T$ corresponds as in \eqref{eqJcf} and \eqref{st}
to the restriction of $W$ to the curve $C$
and to the line bundle $\mathcal{E}$, we conclude that \eqref{eqKernPos} is positive on $X \setminus X_{\R}$.
\end{pf}

In the definition of $(\zeta,\widetilde{\zeta})$-expansive functions, 
we do not need to assume a priori
that $\zeta$ and $\widetilde{\zeta}$ belong to the same torus $T_\nu$.
Indeed, let us assume that $\zeta\in T_\nu$ and $\widetilde{\zeta}\in T_{\widetilde{\nu}}$ and the kernel \eqref{eqKernPos} is positive.
Arguing as in the first part of Proposition \ref{propCont}, we see that for a point $p$ sufficiently close to $q_0 \in X_k$
(and on the appropriate side of $X_k$ with respect to the chosen orientation),
$\frac{\vartheta[{\zeta}](\tauBa{p}-p)}{i\vartheta[{\zeta}] (0)E(p,\tauBa{p})}$
and 
$\frac{\vartheta[\widetilde{\zeta}](\tauBa{p}-p)}{i\vartheta[\widetilde{\zeta}] (0)E(p,\tauBa{p})}$
have signs $(-1)^{\nu_k}$ and $(-1)^{\widetilde{\nu}_k}$, respectively (cf. the proof of \cite[Theorem 2.10]{av3}).
On the other hand, \eqref{eqKernPos} implies that
\[
\abs{T(p)}^2 
\frac
{\vartheta[{\zeta}](\tauBa{p}-p)}
{i\vartheta[{\zeta}] (0)E(p,\tauBa{p})}
\geq
\frac
{\vartheta[\widetilde{\zeta}](\tauBa{p}-p)}
{i\vartheta[\widetilde{\zeta}] (0)E(p,\tauBa{p})}.
\]
Hence the signs of the LHS and the RHS are equal, i.e. $(-1)^{\nu_k} = (-1)^{\widetilde{\nu}_k}$ and hence ${\nu_k} = {\widetilde{\nu}_k}$.

\begin{Rk}
\label{rk26D}
Given $\z$ and $\widetilde{\z}$ in $J(X)$, corresponding determinantal representations
$\lambda_1 \sigma_2 - \lambda_2 \sigma_1 + \gamma$ and $\lambda_1 \sigma_2 - \lambda_2 \sigma_1 + \widetilde{\gamma}$
are only determined up to Hermitian equivalence, see \cite{vinnikov1,vinnikov2,vinnikov5}.
However, given the normalized joint characteristic function and the output determinantal representation 
$\lambda_1 \sigma_2 - \lambda_2 \sigma_1 + \widetilde{\gamma}$,
the input determinantal representation $\lambda_1 \sigma_2 - \lambda_2 \sigma_1 + \gamma$ is determined uniquely.
\smallskip

Assume again that the poles of $y_1$ and $y_2$ are all simple 
(for the case of multiple poles see Remark \ref{rkUs2} below).
Let $\lambda_1 \sigma_2^\prime - \lambda_2 \sigma_1^\prime + \gamma^\prime$ 
be some determinantal representation corresponding to $\zeta$, so that 
$$\rho^* (\lambda_1 \sigma_2 - \lambda_2 \sigma_1 + \gamma) \rho =  (\lambda_1 \sigma_2^\prime - \lambda_2 \sigma_1^\prime + \gamma^\prime),$$
with the normalized section ${\bf u}^{\prime \times}(p)$ satisfying $\rho {\bf u}^{\prime \times}(p) = {\bf u}^{\times}(p)$.
\smallskip

Since $S(p^{(j)})=I$ and using \eqref{st}, we have
$\wt{u}^\times_{j} T(p^{(j)})  =  u^\times_{j}$,
where $\wt{u}^\times_{j}$ is the residue of $\wt{\bf u}^\times(p)$ at $p^{(j)}$.
We further use the notation 
$U^\times = \row_{j=1,\ldots,n}{(u^{\times}_{j})}$
(similarly, we define $\wt{U}^{\times}$ and $U^{\prime \times}$).
Then, the relation between the input and the output determinantal representations is 
given in terms of the values of $T$ at infinity:
\begin{equation}
\label{eq26A}
U^{\times} = \widetilde{U}^{ \times}T_{\infty},
\quad
{\mathrm where}
\quad
T_{\infty} = \Diag{T(p^{(j)})}.
\end{equation}
Furthermore, using $\rho {U}^{\prime \times} = {U}^{\times}$ and $\wt{U}^\times T_{\infty}  =  U^\times$, 
we conclude that the input determinantal representations is uniquely determined, 
namely $\rho$ is uniquely determined, by the behavior of $T$ at infinity 
\begin{equation}
\label{eq26B}
\rho=\widetilde{U}^{ \times}
T_\infty
(\widetilde{U}^{ \prime \times})^{-1}
.
\end{equation}
\end{Rk}

\subsection{The model space}

We continue the discussion about commutative vessels, 
by presenting the model space associated to a commutative two--operator vessel, see also \cite{av1,av3}. 
For the analogue in the single operator case, the de Branges-Rovnyak operator model, see \cite{bc,dbr1,dbr2} and \cite{nf}.
\smallskip

Let $\zeta,\widetilde{\zeta} \in J(X)$ satisfying $\vartheta(\zeta) \neq 0$ and $\vartheta(\widetilde{\zeta}) \neq 0$ 
such that $\zeta$ and $\widetilde{\zeta}$ belong to the same real torus $T_\nu$
and let $T(p)$ be $(\zeta,\widetilde{\zeta})$-expansive function on $X$ (see Theorem \ref{realizationTh})
Then, the kernel \eqref{eqKernPos} defines a reproducing kernel Hilbert space denoted by $\mathcal{H}(T)$.
Its elements are sections of a line bundle of $L_{\widetilde{\zeta}} \otimes \Delta$ that are holomorphic on the domain of analyticity of $T$.
\smallskip

We continue with the definition of the operator model in the case of simple poles \cite[Equation 3-3]{MR1634421},
i.e. when the meromorphic coordinate functions $y_1,y_2$ have only simple poles on $X$, equivalently, 
the divisor $(\nu_0)$ is supported on $n$ distinct points on $X$
(see \cite{av3} and in Section \ref{secMultOp} below for general case).
The counterpart of the operator model \eqref{eqOpModelClass}, denoted by $M^y$, 
is defined on the sections of the line bundle $L_{\widetilde{\zeta}}\otimes \Delta$ analytic in neighborhoods of the poles of $y$. 
It is given by
\begin{equation}
M^{y}f(u)
\label{m_y}
=
y(u)f(u) + \sum_{m=1}^{n}{c_m f(p^{(m)})
\frac
{\vartheta[\widetilde{\zeta}](p^{(m)}-u)}
{\vartheta[\widetilde{\zeta}] (0)E(p^{(m)},u)}
}.
\end{equation}
where $y$ is a meromorphic function on $X$ of degree $n$ with distinct simple poles $p^{(1)},\ldots,p^{(n)}$ 
and $c_1,\ldots,c_n$ are the negatives of the corresponding residues.
The operator $M^y$ is in fact independent of the choice of the local coordinates at the poles of $y$ used to compute the residues,
this point is discussed in Section \ref{secMultOp}.
Furthermore, for any pair of meromorphic functions $y_1,y_2 \in \mathcal{M} (X)$, the operators $M^{y_1}$ and $M^{y_2}$ commute
(see \cite{av3} and Section \ref{secMultOp} below for the theorem in the case of non-simple poles).
\smallskip

The counterpart of the resolvent operator \eqref{resolventOp}, denoted by $R_{\alpha}^{y}$, satisfies 
(for $\alpha$ in the neighborhood of infinity) 
$R_{\alpha}^{y} = \left(M^{y} - \alpha I \right)^{-1}$ and is defined, for an $\alpha$ with $n$ distinct pre-images, 
by (see \cite[Equation 3-4]{MR1634421})
\begin{equation}
\label{Ralpha}
R_{\alpha}^{y}f(u) =
\frac{f(u)}{y(u) - \alpha} -
\sum_{j=1}^{n}{\frac{f(u^{(j)})}{dy(u^{(j)})}\frac{\vartheta[\widetilde{\zeta}](u^{(j)}-u)}{\vartheta[\widetilde{\zeta}](0)E(u^{(j)},u)}},
\end{equation}
where $y(u^{(j)})=\alpha$ for $j=1,\ldots,n$. 
$R_{\alpha}^{y}$ is indeed the resolvent of $M^y$ by Theorem \ref{rkModelOpAlg} below, since $R_{\alpha}^{y} = M^{\frac{1}{y-\alpha}}$
(in the case that $y$ has only simple poles a direct calculation can be done).
Furthermore, the resolvent operators satisfy the resolvent identity (see \cite[Theorem 4.2]{av3})
\[
R^y_\alpha -R^y_\beta
=
(\alpha -\beta) R^y_\alpha R^y_\beta.
\]
The model (commutative two-operator) vessel is the collection (see \cite[Theorem 3.1]{MR1634421}):
\begin{equation}
\label{eqModelVess}
\mathcal{V}_{T} = 
\bracketsA{M^{y_1} \, , \, M^{y_2} \, ; \, \mathcal{H}(T) \, , \, \Phi_{\mathrm{Mod}} \, ,  \, \mathbb C ^n \, ; \,  \sigma_1 \, , \, \sigma_2 \, , \, \gamma \, , \, \widetilde{\gamma}},
\end{equation}
where $\Phi_{\mathrm{Mod}}$ is the evaluation operator at the poles of $y_1$ and $y_2$ 
in the case where these poles are all simple and the output determinantal representation 
is the canonical determinantal representation corresponding to $\wt{\z}$
(in general, $\Phi_{\mathrm{Mod}}$ involves derivatives and the adjustment for the choice of the output determinantal representation,
see the proof of Theorem \ref{vesselAreEq} in Section \ref{secMultOp} below for details).
To define the mapping between the inner space $H$ of a vessel to its model space,
we first consider the mapping, see \cite[Equation 3-5]{MR1634421},
\begin{align}
\label{ModelMap0}
h \mapsto  &
\frac{\xi_1 dy_1(p)+\xi_2 dy_2(p)}{\omega(p)}
\widetilde{P}(\xi_1,\xi_2,p)
\Phi 
\\ & 
(\xi_1A_1+\xi_2A_2 - \xi_1y_1(p) -\xi_2y_2(p))^{-1}h
\nonumber,
\end{align}
which defines a section of $\widetilde{\mathcal{E}}$
(more precisely of $\widetilde{\mathcal{E}}\otimes \mathcal{O}(m-2)(-D_{\rm sing})$)
and is independent of $(\xi_1,\xi_2) \in \C^2$ as long as the resolvent exists (see Proposition \ref{lemmModelSpaceMapProp} below).
In \eqref{ModelMap0}, $h \in H$, $p \in X$
and $\omega(p)$ is the meromorphic differential \eqref{eqOmegaDef} with zeros of order $(m-3)$ at the points of $C$ at infinity
and poles on the divisor of singularities.
$\widetilde{P}(\xi_1,\xi_2,p)$ is the projection onto $\widetilde{\mathcal{E}}(p)$ along
$\widetilde{\mathcal{E}}(p^{(2)})\dotplus\cdots\dotplus\widetilde{\mathcal{E}}(p^{(m)})$ where
$p=p^{(1)},p^{(2)},\ldots,p^{(m)}$ are the $m$ distinct points of the intersections of the 
line $\xi_1\lambda_1+\xi_2\lambda_2 = \xi_1 y_1(p) + \xi_2 y_2(p)$ with the curve $C$.
The projection $\widetilde{P}$ is given by (see \cite[Equations 2.29]{MR97m:30051}):
\[
\widetilde{P}(\xi_1,\xi_2,p) 
= 
\widetilde{\bf u}^{\times}(p) \widetilde{\bf u}^{\times}_l(p)
\frac{\xi_1 \sigma_1 + \xi_2 \sigma_2}{\xi_1 d y_1(p) + \xi_2 d y_2(p)}
,
\]
where $\widetilde{\bf u}^{\times} (p)$ and $ \widetilde{\bf u}^{\times} _l (p)$
are the normalized sections and the left normalized sections of $\widetilde{\mathcal E}$ and $\widetilde{\mathcal E}_{l}$, respectively.
The mapping in \eqref{ModelMap0} then defines a section of $L_{\widetilde{\z}}\otimes \Delta$ by pulling back according to the isomorphism \eqref{eqIsoExplicit}
and hence \eqref{ModelMap0} becomes:
\begin{align}
\label{ModelMap1}
h \mapsto \varphi_h(p) 
\defEq &
{\bf \widetilde{u}}_l^\times(p)
(\xi_1 \sigma_1 + \xi_2 \sigma_2)
\Phi 
\\ & 
(\xi_1A_1+\xi_2A_2 - \xi_1y_1(p) -\xi_2y_2(p))^{-1}h.
\nonumber
\end{align}

\begin{proposition}
\label{lemmModelSpaceMapProp}
The mapping \eqref{ModelMap1} to the model space is well defined 
and does not depend on the choice of $\xi_1$ and $\xi_2$.
\end{proposition}

\begin{proposition}
\label{lemmModelSpaceMapProp123}
The mapping to the model space \eqref{ModelMap1} is injective if and only if the vessel is irreducible.
\end{proposition}

\begin{Tm}
\label{vesselAreEq}
Let $\mathcal{V}$ be an irreducible commutative two-operator vessel.
Then $\mathcal{V}$ is unitary equivalent to the model vessel $\mathcal{V}_{T}$ \eqref{eqModelVess}.
\end{Tm}

Since the proofs of those statements in \cite{MR1634421} are rather sketchy,
we present the comprehensive proofs of Proposition \ref{lemmModelSpaceMapProp} and Proposition \ref{lemmModelSpaceMapProp123}
in this section and the proof of Theorem \ref{vesselAreEq} in Section \ref{secMultOp} below.
Notice that given a $(\zeta, \widetilde{\zeta})$-expansive function $T$, 
we can construct the model vessel in \eqref{eqModelVess},
by taking $\lambda_1 \sigma_2 - \lambda_2 \sigma_1 +\widetilde{\gamma}$ any
determinantal representation corresponding to $\wt{\z}$
(in particular, for the case of simple poles we can take the canonical determinantal representation constructed explicitly in \cite{MR1704479},
see Steps \ref{step1} \& \ref{step2} of the proof of Theorem \ref{preTh} in Section \ref{secProofStruct}),
and then computing the corresponding input determinantal representation  
$\lambda_1 \sigma_2 - \lambda_2 \sigma_1 +\gamma$ as in Remark \ref{rk26D} and Remark \ref{rkUs2}.
\smallskip

We note that two vessels sharing the same "external data" and the same characteristic functions
are unitarily equivalent to the same model vessel $\mathcal{V}_{T}$ \eqref{eqModelVess}.
Thus, we may conclude the following result.
\begin{Cy}
\label{coro12345}
Any two irreducible vessels with the same $E, \sigma_1, \sigma_2 , \gamma$ and $\widetilde{\gamma}$
are unitarily equivalent if and only if the pair of the associated normalized joint characteristic functions are equal.
\end{Cy}

Corollary \ref{coro12345} has been proved in \cite{KLMV} by a substantially different method using the restoration formula.
\smallskip


\begin{pf}[of Propostion \ref{lemmModelSpaceMapProp}]
We have to show that \eqref{ModelMap0} is independent of the choice of $\xi_1,\xi_2$ 
as long as the resolvent exists. That is, we show that
\begin{multline*}
\btul(p) (\xi_1\sigma_1+\xi_2\sigma_2) \Phi (\xi_1A_1+\xi_2A_2-\xi_1y_1(p)-\xi_2y_2(p))^{-1} \\ =
\btul(p) (\eta_1\sigma_1+\eta_2\sigma_2) \Phi (\eta_1A_1+\eta_2A_2-\eta_1y_1(p)-\eta_2y_2(p))^{-1}.
\end{multline*}
Let us multiply on the right by the invertible operator
$$
(\xi_1A_1+\xi_2A_2-\xi_1y_1(p)-\xi_2y_2(p)) (\eta_1A_1+\eta_2A_2-\eta_1y_1(p)-\eta_2y_2(p)),
$$
we obtain that the equality to be verified is 
\begin{multline*}
\btul(p) (\xi_1\sigma_1+\xi_2\sigma_2) \Phi (\eta_1A_1+\eta_2A_2-\eta_1y_1(p)-\eta_2y_2(p)) \\ =
\btul(p) (\eta_1\sigma_1+\eta_2\sigma_2) \Phi (\xi_1A_1+\xi_2A_2-\xi_1y_1(p)-\xi_2y_2(p)).
\end{multline*}
After opening the parentheses and cancelling equal terms this becomes
\begin{multline*}
\btul(p) \left(\xi_1\eta_2 \sigma_1 \Phi A_2 + \xi_2\eta_1 \sigma_2 \Phi A_1 
- \xi_1\eta_2 y_2(p) \sigma_1 \Phi - \xi_2 \eta_1 y_1(p) \sigma_2 \Phi\right) \\ =
\btul(p) \left(\eta_1\xi_2 \sigma_1 \Phi A_2 + \eta_2\xi_1  \sigma_2 \Phi A_1
- \eta_1\xi_2 y_2(p) \sigma_1 \Phi - \eta_2\xi_1 y_1(p) \sigma_2 \Phi\right).
\end{multline*}
Taking everything to the left hand side we obtain that we have to verify the equality
$$
(\xi_1\eta_2-\xi_2\eta_1) \btul(p) \left(\sigma_1 \Phi A_2 - \sigma_2 \Phi A_1 
- y_2(p) \sigma_1 \Phi + y_1(p) \sigma_2 \Phi\right) = 0.
$$
Using the output vessel condition the left hand side equals
\begin{multline*}
(\xi_1\eta_2-\xi_2\eta_1) \btul(p) \left(\widetilde\gamma \Phi - y_2(p) \sigma_1 \Phi + y_1(p) \sigma_2 \Phi\right) \\ =
(\xi_1\eta_2-\xi_2\eta_1) \btul(p) \left(\widetilde\gamma - y_2(p) \sigma_1 + y_1(p) \sigma_2\right) \Phi,
\end{multline*}
and this is of course zero by the definition of $\btul(p)$.
\end{pf}

\begin{pf}[of Propostion \ref{lemmModelSpaceMapProp123}]
If $h \perp \widehat H$ then clearly 
$\Phi(\xi_1 A_1 + \xi_2 A_2 - z I)^{-1} h = 0$ for all $z \not\in \operatorname{spec}(\xi_1 A_1 + \xi_2 A_2)$,
e.g., by expanding Taylor series in $z$ around $\infty$. Therefore $\varphi_h=0$.
\smallskip

Assume conversely that $\varphi_h=0$, i.e.,
$$
\widetilde{P}(\xi_1,\xi_2,p) \Phi ((\xi_1A_1+\xi_2A_2-\xi_1y_1(p)-\xi_2y_2(p))^{-1} h = 0.
$$
We proceed now as with the restoration formula for the complete characteristic function
\cite[Section 10]{KLMV}.
For any $\xi_1,\xi_2,z$ such that the line $\xi_1\lambda_1+\xi_2\lambda_2=z$ intersects $C$ in $n$ distinct affine points, 
we have 
$$
\sum_{p \in X \colon \xi_1 y_1(p) + \xi_2 y_2(p) = z} \widetilde{P}(\xi_1,\xi_2,p) = I,
$$
and we conclude that if $z \not\in \operatorname{spec}(\xi_1 A_1 + \xi_2 A_2)$ then
\begin{align*}
\Phi(\xi_1 A_1 + \xi_2 A_2 - z I)^{-1} h 
= & 
\sum_{p \in X \colon \xi_1 y_1(p) + \xi_2 y_2(p) = z} \widetilde{P}(\xi_1,\xi_2,p) \Phi \times
\\ & 
((\xi_1A_1+\xi_2A_2-\xi_1y_1(p)-\xi_2y_2(p))^{-1} h = 0.
\end{align*}
So if $\det(\xi_1\sigma_1+\xi_2\sigma_2) \neq 0$ 
(this guarantees that the line $\xi_1 y_1 + \xi_2 y_2 = z$ does not intersect $C$ at infinity),
then $\Phi(\xi_1 A_1 + \xi_2 A_2 - z I)^{-1} h = 0$ for $z$ in a neighborhood of $\infty$,
and from the power series expansion $\Phi (\xi_1 A_1 + \xi_2 A_2)^m h = 0$ for all $n$, i.e.,
$h \perp \bigvee_{m=0}^\infty ((\xi_1 A_1 + \xi_2 A_2)^*)^m \Phi^*(E)$.
But
$$
\widehat H = \bigvee_{m=0}^\infty ((\xi_1 A_1 + \xi_2 A_2)^*)^m \Phi^*(E),
$$
since $\det(\xi_1\sigma_1+\xi_2\sigma_2) \neq 0$.
This is a corollary of \cite[Proposition 10.4.1]{KLMV}.
So $h \perp \widehat H$.
\smallskip

To summarize, $\varphi_h = 0$ if and only if $h \perp \widehat H$. 
Therefore $h \mapsto \varphi_h$ is injective
if and only if the vessel is irreducible.
\end{pf}

\section{Statement of the structure theorem}
\label{secMainThm}

Before stating the main result of this section, the counterpart of Theorem \ref{1_1}, 
we fix some notations and conventions which are used in the upcoming sections.

\begin{Notation}
\label{not4A}
Let $y_1$ and $y_2$ be real (i.e. satisfying $y(p) = \overline{y(\tauBa{p})}$) 
meromorphic functions of degrees $n_1$ and $n_2$, respectively, 
generating $\mathcal{M}(X)$ (the field of meromorphic functions on $X$).
Then
\begin{enumerate}[label=(A\arabic*)]
\item
For $\alpha$ and $\beta$ in $\hat{\mathbb{C}}$ and $k=1,2$, 
we denote the $n_k$ pre-images on $X$ under $y_k$ of $\alpha$ and $\beta$,
assumed to be all distinct, by
$\bracketsC{ w_k^{(l)} }{l= 1}{n_k}$ and $\bracketsC{ v_k^{(t)} }{t=1}{n_k} $.
\item
In the case where $y_1$ and $y_2$ have only simple poles,
we denote by $\bracketsC{ p_k^{(j)} }{j=1}{n_k}$ the $n_k$ poles of $y_k$ where $k=1,2$.
\item
\label{not4A_3}
We denote by $\pi: p \mapsto (y_1(p),y_2(p))$ the corresponding birational embedding of $X$
as an irreducible affine algebraic curve $C_0\subset \C^2$.
Its projective closure is denoted by $C \subset \mathbb P^2$.
Notice that  $\pi: X \rightarrow C$ is the normalization of $C$.
\end{enumerate}
\end{Notation}
An important property of the model operator is given in the following lemma.
\begin{lem}
\label{lemMyBounded}
Let $X$ be a compact real Riemann surface and let $y$ be a meromorphic function on $X$.
Let ${\mathcal X}$ be a reproducing kernel Hilbert space of sections of 
$L_{\widetilde{\zeta}} \otimes \Delta$ (where $\widetilde{\zeta} \in T_{\nu}$) 
analytic in an open and connected set $\Omega$ containing the poles of $y$ which is invariant under $M^y$.
Then the model operator $M^y $ is bounded on ${\mathcal X}$.
\end{lem}
\begin{pf}
    Since in a reproducing kernel Hilbert space, strong (or even weak)
    convergence implies pointwise convergence,
    we see that the operator $M^y$ is closed.
 Hence, since $M^y$ is closed and everywhere defined,
then, by the closed graph theorem, $M^y$ is bounded.
\end{pf}

We notice the following lemma is presented in \cite{av3} in the case where $\mathcal X$ is finite dimensional.
This result determines under which conditions $R^y_\alpha$ is the resolvent of the model operator $M^y$.
Since the model operator is bounded (Lemma \ref{lemMyBounded}) it follows that the resolvent exists in a neighborhood of infinity
and hence we may state the following result.

\begin{lem}[{\cite[Section 4]{av3}}]
\label{mainThLem}
Let $\Omega \subseteq X $ be an open and connected set containing the poles of a meromorphic function $y$.
Let ${\mathcal X}$ be a reproducing kernel Hilbert space of sections of $L_{\widetilde{\zeta}} \otimes \Delta$ analytic in $\Omega$.
Then, for $\alpha$ in a neighborhood of infinity and $f \in {\mathcal X}$, 
$R^y_\alpha f$ and $M^y f$ are well-defined analytic sections in $\Omega$.
Furthermore,
if $\mathcal X$ is invariant under $M^y$ then $\mathcal X$
is invariant under $R^y_\alpha$ and under this condition
\begin{equation*}
R^y_\alpha = (M^y-\alpha I)^{-1}.
\end{equation*}
In particular, the kernel of the resolvent operator is trivial, i.e. $\ker R_\alpha^{y} = \{0\}$.
\end{lem}

We notice also that a simple calculation (see the proof of Step \ref{stepMax} of Theorem \ref{preTh}
in Section \ref{secProofStruct}) shows that if $\mathcal X$ is invariant under $R^y_{\alpha^0}$
then it is unvariant under $R^y_{\alpha}$ for $\alpha$ in a neighborhood of $\alpha^0$.

The main result is presented in Theorem \ref{MainTh} below.
To prove Theorem \ref{MainTh}, we first state a simpler result in Theorem \ref{preTh}.
It contains the counterpart of the "if" part of Theorem \ref{1_1} under the assumption that the meromorphic functions 
$y_1(\cdot)$ and $y_2(\cdot)$ have only simple poles. This assumption is dropped later in Theorem \ref{MainTh}.
In these statements, recall again that $T_\nu$ are the real torii given in \eqref{eqRealTorii}.

\begin{theorem}
\label{preTh}
Let $X$ be a compact real Riemann surface and let ${\mathcal X}$ 
be a reproducing kernel Hilbert space of sections of 
$L_{\widetilde{\zeta}} \otimes \Delta$ ( where $\widetilde{\zeta} \in T_{\nu}$) 
analytic in an open and connected set $\Omega$.
We pick two meromorphic functions, $y_1$ and $y_2$, with simple poles generating $\mathcal{M}(X)$,
such that $\Omega$ contains the points above the singular points of $C$ and the poles of $y_1$ and $y_2$.
Furthermore, assume that for some $\alpha,\beta \in \mathbb C$ in the neighborhood of infinity 
such that $\alpha \neq \beta$ and such that their $n$ pre-images lie within $\Omega$, 
the following conditions hold:
\begin{enumerate}[label=(\roman*)]
\item
${\mathcal X}$ is invariant under $M^{y_1}$ and $M^{y_2}$.
\item
\label{assump3}
For every choice of $f,g \in {\mathcal X}$ it holds that
\begin{align}
\label{StructIdent4}
\innerProductReg{R_\alpha^{y_k}f}{g} & -
\innerProductReg{f}{R_\beta^{y_k}g} -
(\alpha - \overline{\beta})\innerProductReg{R_\alpha^{y_k}f}{R_\beta^{y_k}g}
=
\nonumber
\\ &
i(\alpha - \overline{\beta}) \sum_{l,t=1}^{n_k}
{
\frac{f(w^{(l)}) \overline{g(v^{(t)})}}
{dy_k(w^{(l)})\overline{dy_k(v^{(t)})}}
\frac
{\vartheta[\widetilde{\zeta}](w^{(l)}-\tauBa{v^{(t)}})}
{\vartheta[\widetilde{\zeta}] (0)E(w^{(l)},\tauBa{v^{(t)}})}
}.
\end{align}
\end{enumerate}
Then the reproducing kernel of $\mathcal X$ is of the form
\begin{equation}
\label{RKstruct}
K_{\mathcal X}(p,q) =
T(p) K_{\zeta}(p,q) T(q)^* - K_{\widetilde{\zeta}}(p,q)  
\end{equation}
for some $\zeta \in T_{\nu}$ and where $T(\cdot)$ is a $(\zeta, \widetilde{\zeta})$-expansive line bundles mapping.
\end{theorem}

\begin{Rk}
\label{LocCoor}
It easily seen that the structure identity does not depend on the choice of the local coordinates:
each summand on the right hand side of \eqref{StructIdent4} is well defined since $\wt{\z} \in T_\nu$ and hence
$\wt{\z} + \overline{\wt{\z}} = \overline{\kappa} + \kappa$
and since the transition functions of the line bundle
defining $\Delta$ were chosen to be symmetric.
\end{Rk}

\begin{Rk}
\label{rkBlaBla}
One may further consider the structure identity \eqref{StructIdent4} when $\alpha = \overline{\beta}$.
Then, the right hand side of \eqref{StructIdent4} is interpreted as the limit.
Unless specifically stated otherwise, this case is not considered below.
\end{Rk}

The comprehensive counterpart (removing the simple poles assumption and adding the converse statement) 
of Theorem \ref{1_1}, in the setting of compact real Riemann surfaces, is given below.
\begin{theorem}
\label{MainTh}
Let $X$ be a compact real Riemann surface and 
let ${\mathcal X}$ be a reproducing kernel Hilbert space of sections 
of $L_{\widetilde{\zeta}} \otimes \Delta$ 
(where $\widetilde{\zeta} \in T_{\nu}$).
\begin{enumerate}[label=(\alph*)]
\item
\label{stateA}
Assume that the elements of $\mathcal X$ are analytic in an open and connected set $\Omega$.
Let $y_1$ and $y_2$ be two meromorphic functions generating $\mathcal{M}(X)$,
such that $\Omega$ contains the points above the singular points of $C$ and contains the 
pre-images under $(y_1,y_2)$of some 
$(\beta_1^0,\beta_2^0) \in \R^2$.
Furthermore, we assume that 
$\mathcal X$ is invariant under $R^{y_1}_{\beta_1^0}$ and $R^{y_2}_{\beta_2^0}$,
and for some $\alpha_k,\beta_k \in {\mathbb C}$ in a neighborhood of $\beta_k^0$,  
such that $\alpha_k,\beta_k \neq \beta_k^0$, $\alpha_k \neq \beta_k$, $k = 1,2$, 
the structure identity \eqref{StructIdent4} holds.
Then, the reproducing kernel of $\mathcal X$ is of the form
\begin{equation}
\label{rkhsKernel}
K_{\mathcal X}(p,q) =
T(p) K_{\zeta}(p,q) T(q)^*-  K_{\widetilde{\zeta}}(p,q)
\end{equation}
for some $\zeta \in T_{\nu}$ and where $T(\cdot)$ is a $(\zeta, \widetilde{\zeta})$-expansive line bundles mapping.
\item
\label{stateB}
Conversely,
let $T$ be a $(\zeta,\widetilde{\zeta})$-expansive mapping and assume that
$\mathcal X$ has a reproducing kernel of the form \eqref{rkhsKernel}.
Then all the elements of $\mathcal X$ are analytic on $\Omega_T$ (the region of regularity of $T$).
Furthermore, for any $y(\cdot)$, a real meromorphic function on $X$ 
such that all its poles are contained in $\Omega_T$,
$\mathcal X$ is $M^{y}$--invariant and the structure identity \eqref{StructIdent4} holds.
\footnote{\label{fn1} A sign difference between the characteristic functions \eqref{eqCfCol} and \eqref{CCF}, causes
a difference in the statements of Theorem \ref{1_1} and Theorem \ref{MainTh}.}
\end{enumerate}
\end{theorem}

Note that Theorem \ref{MainTh} shows in fact that if the kernel \eqref{rkhsKernel} is positive,
then $T$ admits a meromorphic extension to $X \setminus X_\R$, see Section \ref{secProofStruct}.
\smallskip

In the finite dimensional case, it has been proved \cite[Section 3]{av3} that $T(\cdot)$ 
has the form of a finite Blaschke product on a compact Riemann surface, that is, 
a finite product of the Blaschke factors
\[
b_a(u) = \frac{E(u,a)}{E(u,\overline{a})}\exp \left(-2\pi (a -\overline{a})^t Y u\right).
\]
The proofs of Theorems \ref{preTh} and \ref{MainTh} are presented in Section \ref{secProofStruct}.
We now give the outline of the "only if" part of the proof of Theorem \ref{preTh}.
We start with the observation that by Lemma \ref{mainThLem}
the kernels of the operators $R_\alpha^{y_1}$ and $R_\beta^{y_2}$ are trivial and we have
\begin{equation}
\label{StructIdent5}
R_\alpha^{y_1} = (M^{y_1} - \alpha I) ^{-1}, \qquad R_\beta^{y_2} = (M^{y_2} - \beta I) ^{-1}.
\end{equation}
Then one proceeds as follows:
we start by presenting and constructing a natural two-operator vessel embedding the operators $M^{y_1}$ and $M^{y_2}$ 
and then we follow the next steps:
\begin{enumerate}[label=({\bf  Step {\arabic*}})]
\item
Prove that the colligation conditions for $M^{y_1}$ and $M^{y_2}$ are equivalent to the structure identities for $y_1$ and $y_2$, respectively.
\item
Show that the output vessel condition holds.
\item
Construct the matrix $\gamma$ such that the input vessel condition holds.
\item
Prove that, in our setting, the mapping \eqref{ModelMap1}, between the inner space of the vessel to the model space, is the identity mapping.
\item
Present the reproducing kernel in terms of the joint characteristic function.
\item
Show that the input and output determinantal representations are fully saturated.
\item
Conclude, by the reproducing kernel Hilbert space properties, that the reproducing kernel has the desired structure \eqref{RKstruct}.
\end{enumerate}
\section{Subspaces of \texorpdfstring{$H^2_{\widetilde{\zeta}}$}{ a } and a version of Beurling's Theorem}
\label{secH2}

In this section, we present three versions of Beurling's theorem on finite bordered Riemann surfaces. 
For a short survey and related previous results, see Section \ref{sucSectionIntoCRS}.
\smallskip

Let $\mathscr{S}$ be an open Riemann surface so that $\mathscr{S} \cup \partial \mathscr{S}$ is a finite bordered Riemann surface of genus $g_{\mathscr{S}}$ 
whose boundary consists of $k \geq 1 $ connected components, denoted by $X_0, \ldots, X_{k-1}$.
The double of $\mathscr{S}$ is a compact Riemann surface $X$ with a natural antiholomorphic involution $\tau$, 
turning $X$ into a compact real Riemann surface of genus $g = 2g_{\mathscr{S}}	+k-1$.
The boundary $\partial \mathscr{S}$ coincides with the set of fixed points of $\tau$ on $X$ (denoted by $X_{\mathbb R}$).
Furthermore, $X$ is a compact real Riemann surface of dividing type since 
$X \setminus X_{\mathbb R}$ contains two connected components $X_- = \mathscr{S}$ and $X_+$.
\footnote{The usual convention is $X_+=\mathscr{S}$, 
but $X_-=\mathscr{S}$ is more convenient for us since with our choice of notation the characteristic function
is contractive in the lower half plane, cf. Footnote \ref{fn1}.}
\smallskip

As discussed in \cite[Section 2]{av2}, in a more general context,
a flat unitary line bundle on $\mathscr{S}\cup\partial \mathscr{S}$ can be uniquely extended to a 
flat unitary line bundle on $X$  such that certain symmetry properties are fulfilled,
more precisely, to a line bundle $L_{\widetilde{\zeta}}$ where $\widetilde{\zeta} \in T_0$.
\smallskip

Let now $L_{\widetilde{\zeta}}$ where $\widetilde{\zeta} \in T_0$ be a flat unitary line bundle on $X$ 
and let $\Delta$ be a square root of the canonical line bundle as in Section \ref{sub21}.
Then, the corresponding Hardy space consists essentially of sections of the line bundle $L_{\widetilde{\zeta}} \otimes \Delta$ analytic on $\mathscr{S}$ and satisfying
(in the sense of non tangential boundary values)
\[
\sum_{j=0}^{k-1}{\int_{X_j}{f(p)^*f(p)}} < \infty,
\]
it becomes a Hilbert space equipped with the inner product
\[
\innerProductReg{f}{g} =  2 \pi \sum_{j=0}^{k-1}{\int_{X_j}{g(p)^*f(p)}}.
\]
See \cite{av2} for the precise definitions.
Following the classical modification of the inner product in $H^2$, 
also here, for the sake of simplicity, 
the inner product is multiplied by $2\pi$ and hence the reproducing kernel is 
$
-
\frac
{\vartheta[\widetilde{\zeta}](\tauBa{w} - u^{(j)})}
{ i \vartheta[\widetilde{\zeta}] (0)E(u^{(j)},\tauBa{w})}
$
and not multiplied by $\frac{1}{2\pi}$ (as in \cite{av2}).
To simplify notations, we set $H^2_{\widetilde{\zeta}} = H^2 ( L_{\widetilde{\zeta}}\otimes \Delta , X_-)$.
\smallskip

Before turning to the main theorem, let us recall the definition of dividing functions on a compact real Riemann surface.
\begin{definition}
A real meromorphic function $y$ on a compact real Riemann surface $X$ is {\it dividing} 
if $u \in X_{\mathbb R}$ if and only if $y(u) \in {\mathbb R}$.
\end{definition}
It is easy to see that if $y$ is a dividing function 
on a compact real Riemann surface $X$,
then $X$ is of dividing type and 
$y$ maps $X_+$ onto either the upper or the lower half-plane;
we will always assume that $y(X_+)={\mathbb C}_+$.
A known result regarding dividing functions is presented below.
However, we note that only the first part is used in the sequel.
\begin{proposition}[{\cite{MR0021108,ahlfors}, \cite{MR1072300,MR1070485} and \cite[Proposition 5.2]{shamovich2014livsic}}]
\label{elisLemma}
Let $y$ be a dividing function on $X$.
Then $y$ has only real simple poles and simple
zeros and its residues at the poles, with respect to a real local coordinate with positive orientation, are negative.
Conversely, if $X$ is of dividing type and $y$ is a real meromorphic function on $X$ with simple real poles and negative residues with respect to positive real
local coordinate, then $y$ is dividing.
\end{proposition}

For more information about dividing functions, see also \cite{gabard2012ahlfors}.
In particular, for any compact Riemann surface $X$ of dividing type,
there always exists a pair of dividing functions that generates $\mathcal{M}(X)$,
see also \cite{kummer2015real} for a far reaching high dimensional generalization.
\smallskip

We now turn to state the first version of Beurling's theorem on finite bordered Riemann surfaces.
\begin{theorem}[Beurling's Theorem for Finite bordered Riemann surfaces: version I]
\label{BLth}
Let $\mathscr{S}$ be a finite bordered Riemann surface and let $X$ be its double.
Let $\widetilde{\zeta} \in {T}_0$, 
let $H^2_{\widetilde{\zeta}}$ be the corresponding Hardy space on $X_-$ and 
let $y_1$ and $y_2$ be dividing functions on $X$ generating $\mathcal{M} (X)$.
Furthermore, assume that for 
$\mathcal{H} \subseteq H^2_{\widetilde{\zeta}}$ 
the following conditions hold:
\begin{enumerate}
\item $\mathcal{H}$ is a closed subspace of $H^2_{\widetilde{\zeta}}$ and is invariant under the multiplication operators 
$\mathcal{M}_{\frac{1}{y_1(\cdot)-\overline{\alpha}}}$ and $\mathcal{M}_{\frac{1}{y_2(\cdot)-\overline{\alpha}}}$ 
for every $\alpha$ in the lower half-plane $\mathbb C_{-}$.
\item The elements of $\mathcal{H}^{\perp}$ 
(the orthogonal complement of $\mathcal{H}$) 
have analytic extensions with bounded point evaluations in a connected neighborhood of the poles of $y_1$ and $y_2$ 
and of the pre-images of the singular points of $C$ (see Notation \ref{not4A} \ref{not4A_3}).
\end{enumerate}
Then $\mathcal{H}$ is of the form
\[
\mathcal{H} = T H^2_{{\zeta}},
\]
where $T$ is a $(\zeta, \widetilde{\zeta})$-inner function for some $\zeta \in T_0$.
\end{theorem}

Notice that by Proposition \ref{propCont}, $\zeta, \widetilde{\zeta} \in T_0$ implies that if 
$T$ is $(\zeta,\widetilde{\zeta})$-expansive on $X$ then it is contractive on $X_{-} = \mathscr{S}$ 
and thus admits non-tangential boundary values almost everywhere on $X_{\R} = \partial \mathscr{S}$.
We say that $T$ is $(\zeta,\widetilde{\zeta})$-{\it inner} if the non-tangential boundary 
value is of absolute value one almost everywhere on $\partial \mathscr{S}$.
\smallskip

We note that the operator $\mathcal{M}_f$ denotes the conventional multiplication operator by a function $f$ 
while $M^y$, as before, denotes the model operator.
Furthermore, for $\alpha \in \mathbb{C}_-$ and $y$ a dividing function, the operator $\mathcal{M}_{\frac{1}{y(u)-\overline{\alpha}}}$ is well-defined on $H^2_{\widetilde{\zeta}}$.
This follows since $\overline{\alpha} \in \C_+$ and therefore the function $\frac{1}{y(u)-\overline{\alpha}}$ is meromorphic on $X$ and analytic in $\mathscr{S} \cup \partial \mathscr{S}$.
Hence, $\mathcal{M}_{\frac{1}{y(u)-\overline{\alpha}}}$ is bounded and sends $H^2_{\widetilde{\zeta}}$ to $H^2_{\widetilde{\zeta}}$.
\smallskip

Before heading to the proof, we present several preliminary results.
Our first goal is to show that the structure identity \eqref{StructIdent4} holds for all elements in $H^2_{\widetilde{\zeta}}$. 
We start with the following result and we recall the proof for the sake of completeness.
\begin{lem}[{\cite[Theorem 4.3]{av3}}]
\label{R_alpha_evs}
Let $\alpha \in \mathbb C$ have $n$ distinct pre-images with respect to a real meromorphic function $y$
and let $\widetilde{\zeta} \in J(X)$ such that $ \theta(\widetilde{\zeta}) \neq 0$.
Then the Cauchy kernels $K_{\widetilde\zeta}(\cdot,w)$ are the eigenvectors of the resolvent operator 
$R^y_\alpha$ with eigenvalues $\frac{1}{\overline{y(w)}-\alpha}$.
\end{lem}
\begin{pf}
We apply the resolvent operator $R^y_\alpha$ \eqref{Ralpha} on $K_{\widetilde{\zeta}}(\cdot,w)$.
Then, a direct computation, using the collection formula \cite[Lemma 4.1]{av3}, yields the following identity
\begin{align*}
\big( R_\alpha^y \, & K_{\widetilde{\zeta}}(\cdot,w) \big) (v)
=
\frac{K_{\widetilde{\zeta}}(v,w)}{y(v)-\alpha}
-
\sum_{j=1}^{n}
{
\frac{K_{\widetilde{\zeta}}(u^{(j)},w)}{dy(u^{(j)})}
\frac
{\vartheta[\widetilde{\zeta}](u^{(j)}-v)}
{ \vartheta[\widetilde{\zeta}] (0)E(u^{(j)},v)}
}
\\
=
&
\frac{K_{\widetilde{\zeta}}(v,w)}{y(v)-\alpha}
-
\frac{-1}{i}
\sum_{j=1}^{n}
{
\frac{1}{dy(u^{(j)})}
\frac
{\vartheta[\widetilde{\zeta}](\tauBa{w} - u^{(j)})}
{ \vartheta[\widetilde{\zeta}] (0)E(\tauBa{w},u^{(j)})}
\frac
{\vartheta[\widetilde{\zeta}](u^{(j)}-v)}
{ \vartheta[\widetilde{\zeta}] (0)E(u^{(j)},v)}
}
\\
=
&
\frac{K_{\widetilde{\zeta}}(v,w)}{y(v)-\alpha}
-
\frac
{- \vartheta[\widetilde{\zeta}](\tauBa{w} - v)}
{ i \vartheta[\widetilde{\zeta}] (0)E(\tauBa{w},v)}
\left(
\frac{1}{y(v)-\alpha}
-
\frac{1}{\overline{y(w)}-\alpha}
\right)
\\
=
&
\frac{K_{\widetilde{\zeta}}(v,w)}{y(v)-\alpha}
-
K_{\widetilde{\zeta}}(v,w)
\left(
\frac{1}{y(v)-\alpha}
-
\frac{1}{\overline{y(w)}-\alpha}
\right)
\\
=
&
\frac{1}{\overline{y(w)}-\alpha}K_{\widetilde{\zeta}}(v,w)
.
\end{align*}
\end{pf}
Using the preceding lemma, we may conclude and prove that the structure identity holds on a dense subset of $H^2_{\widetilde{\zeta}}$.
\begin{lem}
\label{stIdH2}
Let $\alpha,\beta \in \mathbb C_{-}$. Then the structure identity \eqref{StructIdent4} holds on the linear span of Cauchy kernels inside $H^2_{\widetilde{\zeta}}$.
\end{lem}
\begin{pf}
Since $y$ is dividing, it maps $X_+$ to $\mathbb C_+$ and $X_-$ to $\mathbb C_-$. 
Thus, for $\alpha,\beta \in \mathbb C _-$, their pre-images are in $X_-$.
Hence, for all $f \in H^2_{\widetilde{\zeta}}$, $R_\alpha^y f$ and $R_\beta^y f$ 
are well-defined analytic sections of $L_{\widetilde{\zeta}}\otimes \Delta$ on $X_{-}$.
\smallskip

It is enough to verify that \eqref{StructIdent4} holds on a pair of (minus) kernels functions
$f(u) = K_{\widetilde{\zeta}}(u,v)$ and $g(u) = K_{\widetilde{\zeta}}(u,v)$, where $w,v \in X_{-}$.
Then, starting with the left hand side of \eqref{StructIdent4} and using
Lemma \ref{R_alpha_evs}, we compute separately the three components.
The first two components, using Lemma \ref{R_alpha_evs}, are given by:
\begin{align}
\label{firstcomp}
\innerProductReg{R_\alpha^y K_{\widetilde{\zeta}}(u,w)}{K_{\widetilde{\zeta}}(u,v)}
& =
\frac{-1}{\overline{y(w)}-\alpha}
\innerProductReg{K_{\widetilde{\zeta}}(u,w)}{-K_{\widetilde{\zeta}}(u,v)}
\\ & =
\nonumber
\frac{-1}{\overline{y(w)}-\alpha}K_{\widetilde{\zeta}}(v,w),
\end{align}
and similarly,
\begin{equation}
\label{secondcomp}
\innerProductReg{K_{\widetilde{\zeta}}(u,w)}{R_\beta^y K_{\widetilde{\zeta}}(u,v)}
 =
\frac{-1}{y(v)-\overline{\beta}}K_{\widetilde{\zeta}}(v,w).
\end{equation}
    The third element, applying Lemma \ref{R_alpha_evs} once again, is 
\begin{align}
\label{thirdcomp}
\left( \alpha - \overline{\beta} \right)
\innerProductReg{R_\alpha^y K_{\widetilde{\zeta}}(u,w) }{R_\beta^y K_{\widetilde{\zeta}}(u,v)}
= &
\frac{\alpha - \overline{\beta} }{\overline{y(w)}-\alpha}
\innerProductReg{K_{\widetilde{\zeta}}(u,w)}{R_\beta^y K_{\widetilde{\zeta}}(u,v)}
\\
= & -
\frac
{(\alpha - \overline{\beta})K_{\widetilde{\zeta}}(v,w)}
{(\overline{y(w)}-\alpha)(y(v)-\overline{\beta})}.
\nonumber
\end{align}
Summing all three components together, \eqref{firstcomp}, \eqref{secondcomp} and \eqref{thirdcomp}, one may conclude:
\begin{align}
\label{StructIdent123}
\innerProductReg{R_\alpha^y f}{g} - &
\innerProductReg{f}{R_\beta^y g} -
(\alpha - \overline{\beta})
\innerProductReg{R_\alpha^y f}{R_\beta^y g}
=
\\ = &
\nonumber
-
K_{\widetilde{\zeta}}(v,w)
\left(
\frac{1}{\overline{y(w)}-\alpha} -
\frac{1}{y(v)-\overline{\beta}} -
\frac
{\alpha - \overline{\beta}}
{(\overline{y(w)}-\alpha)(y(v)-\overline{\beta})}
\right)
\\
\nonumber
= &
-
K_{\widetilde{\zeta}}(v,w)
\frac{y(v) - \overline{y(w)}}
{(\overline{y(w)}-\alpha)(y(v)-\overline{\beta})}.
\end{align}
On the other hand, the right hand side (we use the notation $RHS$) of \eqref{StructIdent4}
can be simplified by using \cite[Lemma 4.1]{av3} twice.
The first part of the calculation is:
\begin{align}
\nonumber
RHS
=
&
i(\alpha - \overline{\beta})
\sum_{l,t=1}^{n}
{
\frac
{K_{\widetilde{\zeta}}(w^{(l)},w) \overline{K_{\widetilde{\zeta}}(v^{(t)},v)}}
{dy(w^{(l)})\overline{dy(v^{(t)})}}
\frac
{\vartheta[\widetilde{\zeta}](w^{(l)}-\tauBa{v^{(t)}})}
{\vartheta[\widetilde{\zeta}] (0)E(w^{(l)},\tauBa{v^{(t)}})}
}
\\
\nonumber
=
&
-
\sum_{t=1}^{n}
\frac{\overline{K_{\widetilde{\zeta}}(v^{(t)},v)}}{\overline{dy(v^{(t)})}}
\sum_{l=1}^{n}
\frac{\alpha - \overline{\beta}}{dy(w^{(l)})}
\frac
{\vartheta[\widetilde{\zeta}](\tauBa{w} - w^{(l)})}
{ \vartheta[\widetilde{\zeta}] (0)E(\tauBa{w},w^{(l)})}
\frac
{\vartheta[\widetilde{\zeta}](w^{(l)}-\tauBa{v^{(t)}})}
{\vartheta[\widetilde{\zeta}] (0)E(w^{(l)},\tauBa{v^{(t)}})}
\\
\nonumber
=
&
i(\alpha - \overline{\beta})
\left(
\frac{1}{\overline{\beta}-\alpha}
-
\frac{1}{\overline{y(w)}-\alpha}
\right)
\sum_{t=1}^{n}
{
\overline{ \left( \frac{K_{\widetilde{\zeta}}(v^{(t)},v)}{dy(v^{(t)})} \right) }
\overline{K_{\widetilde{\zeta}}(w,\tauBa{v^{(t)}})}
}
\\
\label{Req2}
=
&
i
\left(
1
+
\frac{\alpha - \overline{\beta}}{\overline{y(w)}-\alpha}
\right)
\overline{
\sum_{t=1}^{n}
{
\frac{1}{dy(v^{(t)})}
\frac
{\vartheta[\widetilde{\zeta}](\tauBa{v}-v^{(t)})}
{ \vartheta[\widetilde{\zeta}] (0)E(\tauBa{v},v^{(t)})}
\frac
{\vartheta[\widetilde{\zeta}]( v^{(t)} - w )}
{ \vartheta[\widetilde{\zeta}] (0)E(v^{(t)},w)}
}
}
.
\end{align}
Using the collection formula once again, we have the following equality
\begin{align}
\nonumber
&
-
\overline{
\sum_{t=1}^{n}
\frac{1}{dy(v^{(t)})}
\frac
{\vartheta[\widetilde{\zeta}](\tauBa{v}-v^{(t)})}
{ \vartheta[\widetilde{\zeta}] (0)E(\tauBa{v},v^{(t)})}
\frac
{\vartheta[\widetilde{\zeta}]( v^{(t)} - w )}
{ \vartheta[\widetilde{\zeta}] (0)E(v^{(t)},w)}
}
\\
\nonumber
=
&
-
\overline{
(-i)
\frac
{\vartheta[\widetilde{\zeta}](\tauBa{v}-w)}
{i\vartheta[\widetilde{\zeta}] (0)E(w,\tauBa{v})}
}
\left(
\frac{1}{\overline{y(w)}-\overline{\beta}}
-
\frac{1}{y(v)-\overline{\beta}}
\right)
\\
\label{Req1}
=
&
-
i
\left(
\frac{1}{\overline{y(w)}-\overline{\beta}}
-
\frac{1}{y(v)-\overline{\beta}}
\right)
K_{\widetilde{\zeta}}(v,w)
.
\end{align}
Substituting \eqref{Req1} in \eqref{Req2} leads to
\begin{align*}
RHS = &
i
\left(
1
+
\frac{\alpha - \overline{\beta}}{\overline{y(w)}-\alpha}
\right)
i
\left(
\frac{1}{\overline{y(w)}-\overline{\beta}}
-
\frac{1}{y(v)-\overline{\beta}}
\right)
K_{\widetilde{\zeta}}(v,w)
\\
=
&
-
\frac{y(v)-\overline{y(w)}}
{(y(v)-\overline{\beta})(\overline{y(w)}-\alpha)}
K_{\widetilde{\zeta}}(v,w),
\end{align*}
as in \eqref{StructIdent123}.
\end{pf}

Moreover, using Lemma \ref{stIdH2}, we show below that the operator $R_{\alpha}^{y}$ is a bounded operator on $H^2_{\widetilde{\zeta}}$.

\begin{lem}
\label{R_alpha_evs1}
Let $\alpha \in \mathbb C _ {-}$ and let $y$ be a dividing function.
Then the resolvent operator $R_{\alpha} ^{y}$ is a well-defined bounded operator 
on $H^2_{\widetilde{\zeta}}$.
\end{lem}
\begin{pf}
Using Lemma \ref{stIdH2}, the structure identity \eqref{StructIdent4}
holds for any linear combination of Cauchy kernels.
We use the structure identity in order to prove the boundedness of $R_{\alpha} ^{y}$.
In \eqref{StructIdent4}, we choose $f=g$ to be a linear combination of Cauchy kernels and we set $\beta = \alpha \in \mathbb{C}_{-}$.
Then we have:
\begin{align}
\label{StructIdent16}
\imag \innerProductReg{R_\alpha^{y}f}{f} & -
\imagg{\alpha} \norm{ R_\alpha^{y}f } ^2
=
\\ &
\nonumber
i \imagg{\alpha} \sum_{l,t=1}^{n}
{
\frac{f(v^{(l)})
\overline{f(v^{(t)})}}
{dy(v^{(l)})\overline{dy(v^{(t)})}}
\frac
{\vartheta[\widetilde{\zeta}](v^{(l)}-\tauBa{v^{(t)}})}
{\vartheta[\widetilde{\zeta}] (0)E(v^{(l)},\tauBa{v^{(t)}})}
}
.
\end{align}
The right hand side of Equation \eqref{StructIdent16}, in view of Remark \ref{LocCoor},
does not depend on the local coordinates choice.
As a result, we deduce from \eqref{StructIdent16} the following inequality:
\begin{equation}
\label{StructIdent977979}
\norm{ R_\alpha^{y}f } ^2
\leq
C_{\alpha} \ \norm{ R_\alpha^{y}f } \
\norm{f}
+
D_{\alpha}\norm{f}^2
,
\end{equation}
for some constants $C_{\alpha}$ and $D_{\alpha}$ (depending only on $\alpha$);
we use here the fact that the point evaluations on $H^2_{\widetilde{\zeta}}$ are bounded.
The inequality \eqref{StructIdent977979} is true for any $f$ in a dense subset of $H^2_{\widetilde{\zeta}}$.
Hence, dividing by $\norm{f}^2$ and taking the supremum over $f$ implies that the operator $R^{y}_{\alpha}$ is bounded
on this dense subspace and therefore extends to a bounded operator on $H^2_{\widetilde{\zeta}}$.
We denote the extension by $T$ and we see that $Tf(p) = R^y_{\alpha}f(p)$ for all $f\in H^2_{\widetilde{\zeta}}$ and all $p\in X_{-}$
again by the boundedness of the point evaluations.

\end{pf}

Combining the last two results, Lemma \ref{stIdH2} and Lemma \ref{R_alpha_evs1},
we conclude that since $R_\alpha^y$ is bounded and since the structure identity holds
on a dense subset, the structure identity holds in $H^2_{\widetilde{\zeta}}$.

\begin{corollary}
\label{stIdH3}
Let $\alpha,\beta \in \mathbb C_{-}$ and let $y$ be a dividing function on a compact real Riemann surface $X$. 
Then the structure identity \eqref{StructIdent4} holds in $H^2_{\widetilde{\zeta}}$.
\end{corollary}

The link between $R^y_\alpha$-invariant subspaces and subspaces which are
invariant under multiplication operators is illustrated in the following lemma.

\begin{lem}
\label{MM}
    Let $\mathcal{H}$ be a subspace of $H^2_{\widetilde{\zeta}}$ and $\alpha \in \mathbb C _{-}$.
    Then $\mathcal{H}$ is $R^y_\alpha$-invariant
    if and only if the orthogonal complement $\mathcal{H}^{\perp}$
    is invariant under the multiplication by $\frac{1}{y(\cdot)-\overline{\alpha}}$.
\end{lem}
\begin{pf}
Since $R_\alpha^y$ is bounded by Lemma \ref{R_alpha_evs1}, it is sufficient to prove that 
$(R_\alpha^y)^* f = {\mathcal M}_{\frac{1}{y(\cdot)-\overline{\alpha}}} f$ for $f$ in $H^2_{\widetilde\zeta}$,
where $\mathcal{M}$ is the multiplication operator. 
One may obtain the following
\begin{align}
\label{eval1}
\innerProductReg{ K_{\widetilde{\zeta}}(\cdot,v)  }{ R^y_\alpha K_{\widetilde{\zeta}}(\cdot,w)}
& =
-\innerProductReg{ \left(R^y_\alpha \right)^* K_{\widetilde{\zeta}}(\cdot,v) }{ -K_{\widetilde{\zeta}}(\cdot,w) }
\\
& =
\nonumber
-\left(\left( R^y_\alpha \right)^* K_{\widetilde{\zeta}}(\cdot,v)\right) (w).
\end{align}
On the other hand, using Lemma \ref{R_alpha_evs}, we have
\begin{align}
\label{eval2}
\nonumber
\innerProductReg{ K_{\widetilde{\zeta}}(\cdot,v) }{ R^y_\alpha K_{\widetilde{\zeta}}(\cdot,w)  }
& =
- \innerProductReg{ K_{\widetilde{\zeta}}(\cdot,v) }{ - \frac{1}{\overline{y(w)}-\alpha}  K_{\widetilde{\zeta}}(\cdot,w)  } \\
& =
- 
\frac{1}{y(w)-\overline{\alpha}}  K_{\widetilde{\zeta}}(w,v).
\end{align}
Hence, combining \eqref{eval1} and \eqref{eval2} 
yields the desired result
for $f = K_{\widetilde{\zeta}}(\cdot,v)$, a (minus) kernel function,
and hence for all $f$ in $H^2_{\widetilde{\zeta}}$.
\end{pf}

To apply Theorem \ref{preTh},
one needs to show that the structure identity holds also on $\partial \mathscr{S}$ in a neighborhood of the poles of $y_1$ and $y_2$.
This point is proved in Lemma \ref{MM1} below. 
\begin{lem}
\label{MM1}
Let $X$ be a compact real Riemann surface of dividing type, 
let $y(\cdot)$ be a dividing function on $X$ 
and let $H^2_{\widetilde{\zeta}}$ be the Hardy space corresponding to 
$\widetilde{\zeta} \in T_0$.
Furthermore, 
we assume that $\mathcal{H}$ is a closed subspace which is invariant under $\mathcal{M}_{\frac{1}{y(u)-\overline{\alpha}}}$
where $\alpha \in \C_{-}$
and the elements of $\mathcal{H}^{\perp}$ have analytic extensions 
in a neighborhood of the poles of $y(\cdot)$.
Then:  
\begin{enumerate}
\item 
The subspace $\mathcal{H}^{\perp}$ is $R_{\alpha_0}^y$-invariant for any 
$\alpha_0\in\R$ in a neighborhood of infinity.
\item
\label{lemAlpha0SI}
The structure identity may be extended to $\alpha_0 \neq \beta_0$ 
in $\mathbb R$ in a neighborhood of infinity
(see Remark \ref{rkBlaBla}, for the case where $\alpha_0=\beta_0\in \mathbb R$). 
\end{enumerate}
\end{lem}
\begin{pf}
We fix $\alpha_0 \in \R$ and let $f$ be an element of $H^2_{\widetilde{\zeta}}$ 
that has an analytic continuation in the neighborhood of the poles of $y$ and of the fiber of $y$ above $\alpha_0$,
so that $R^y_{\alpha_0} f$ is a well defined analytic section of $L_{\widetilde{\zeta}} \otimes \Delta$ on $X_-$.
First, we show explicitly that $R_{\alpha_0} ^y \, f$ belongs to $H^2_{\widetilde{\zeta}}$.
The strategy is to divide the integration path $X_{\R}(\e)$ (a contour approximating $X_\R$, see \cite{av2} for the precise details) into two parts.
The first is over the set of arcs near the pre-images of $\alpha_0$ (where $f$ can be continued analytically).
We denote the first integration path by $X^{\alpha_0}_{\R}(\e)$.
The second path, $X^{\alpha_0}_{\R}(\e)^c \defEq X_{\R}(\e) \backslash X^{\alpha_0}_{\R}(\e)$, 
is contained in a compact set where $\frac{1}{y(u)-\alpha_0}$ can be bounded.
Namely, we examine
\begin{equation}
\label{eqNormFdivide}
\sup_{\e > 0} \int_{X^{\alpha_0}_{\R}(\e)} \abs{R^y_{\alpha_0} f} ^2
+
\sup_{\e > 0} \int_{X^{\alpha_0}_{\R}(\e)^c}  \abs{R^y_{\alpha_0} f} ^2
.
\end{equation}
We note that since $f$ has analytic continuation in a neighborhood of the fiber above $\alpha_0$,
$R^y_{\alpha_0} f$ is holomorphic on this neighborhood and hence, in particular, continuous there. 
Thus the first summand is bounded.
\smallskip

On $X^{\alpha_0}_{\R}(\e)^c$, $\frac{1}{y(u)-\alpha_0}$ is bounded from above by some constant $0<M$ 
and the second summand of \eqref{eqNormFdivide} becomes
\begin{align}
\nonumber
\sup_{\e > 0} & \int_{X^{\alpha_0}_{\R}(\e)^c} \abs{R^y_{\alpha_0} f} ^2
= 
\sup_{\e > 0} 
\int_{X^{\alpha_0}_{\R}(\e)^c}
\abs{
\frac{f(u)}{y(u) -\alpha} - 
\sum_{j=1}^n 
h_{j}
\frac{\vartheta[\zeta](u^{(j)}_0-u)}{\vartheta[\zeta](0)E(u^{(j)}_0,u)}
}^2
\\ \leq & \nonumber
\sup_{\e > 0} 
\left(
\sqrt{\int_{X^{\alpha_0}_{\R}(\e)^c} \left|\frac{f(u)}{y(u)-\alpha_0}\right|^2}
+ 
\sqrt{\int_{X^{\alpha_0}_{\R}(\e)^c} \left|\sum_{j=1}^n h_{j}\frac{\vartheta[\zeta](u^{(j)}_0-u)}{\vartheta[\zeta](0)E(u^{(j)}_0,u)}\right|^2}
\right)^2
\\ \leq &
\sup_{\e > 0} 
\left(
M
\sqrt{\int_{X^{\alpha_0}_{\R}(\e)^c} \left|f(u)\right|^2}
+ 
\sqrt{\int_{X^{\alpha_0}_{\R}(\e)^c} \left|\sum_{j=1}^n h_{j} \frac{\vartheta[\zeta](u^{(j)}_0-u)}{\vartheta[\zeta](0)E(u^{(j)}_0,u)}\right|^2}
\right)^2,
\label{eqNormRf2}
\end{align}
where the coefficients $\left( h_{j} \right)_{j=1}^n$ depend on the pre-images of $\alpha_0$.
The statement follows since $f$ belongs to $H^2_{\widetilde{\zeta}}$ and furthermore 
$\frac{\vartheta[\zeta](u^{(j)}_0-u)}{\vartheta[\zeta](0)E(u^{(j)}_0,u)}$ is bounded on a compact set containing $X^{\alpha_0}_{\R}(\e)^c$.
Thus, $f \in \mathcal{H} ^\perp$ implies that $R^y_{\alpha_0} f$ belongs to $H^2_{\widetilde{\zeta}}$.
\smallskip

In order to continue and conclude that $\mathcal{H} ^\perp$ is invariant under $R^y_{\alpha_0}$, 
we fix a sequence $\left( \alpha_j \right)_{j=1}^{\infty} \subseteq \C_{-}$ converging to $\alpha_0 \in \R$ in a neighborhood of infinity
and we show that $R^y_{\alpha_j} f$ converges to $R^y_{\alpha_0} f$ in norm.
Also here, we divide the integration path, now on $X_\R$, into two parts:
\begin{align}
\lim_{i \rightarrow \infty}
\normTreA{R^y_{\alpha_0} f  - R^y_{\alpha_i} f }{H^2_{\widetilde{\zeta}}}{2}
= &
\lim_{i \rightarrow \infty}
\int_{X^{\alpha_0}_{\R}(0)} 
\abs{ R^y_{\alpha_0} f(u) - R^y_{\alpha_i} f(u)}^2
+ \nonumber \\ &
\lim_{i \rightarrow \infty}
\int_{X^{\alpha_0}_{\R}(0)^c} 
\abs{R^y_{\alpha_0} f(u) - R^y_{\alpha_i} f(u) }^2.
\label{eqNormRf3}
\end{align}
The second summand tends to zero due to a similar analysis as in \eqref{eqNormRf2}.
To show that the first summand also vanishes, we note that in general,
for any section $f$ analytic on a neighborhood $U$ of the fiber of $y$ over $\alpha_0$ (assuming the fiber is unramified),
$(R^y_{\alpha} f)(p)$ is jointly analytic in
$(\alpha,p) \in \{|\alpha-\alpha_0|<\e\} \times U$ (for sufficiently small $\e$).
The analyticity in $p$ follows immediately from the definition of $R^y_\alpha$.
The analyticity in $\alpha$ follows by the following argument.
We examine $(R^y_{\alpha} f)(p)$ as a function of $\alpha$, where $p = p_0$ is fixed.
The first term in $(R^y_{\alpha} f)(p_0)$, that is $\frac{f(p_0)}{y(p_0)- \alpha}$, is meromorphic with a simple pole at $\alpha = y(p_0)$.
On the other hand, also the second term, containing the summation of the Cauchy kernels, has a simple pole at $\alpha = y(p_0)$.
Furthermore, both terms have the same residue at the simple pole $y(p_0)$ 
and hence their subtraction yields a holomorphic function of $\alpha$.
It then remains to apply the Hartog's Theorem, see for instance \cite[Theorem 2.2.8]{hormander1973introduction}, 
to conclude the local joint-analyticity property of $(R^y_{\alpha} f)(p)$.
\smallskip

As a consequence, $R^y_{\alpha_i} f(p)$ converges to $R^y_{\alpha_0} f(p)$ uniformly in $p$ on a neighborhood
of each $u^{(j)}_0$ in the fiber over $\alpha_0$, in particular on a small arc of the boundary near $u^{(j)}_0$ and so 
the second limit in \eqref{eqNormRf3} is also zero.
Hence $\|R^y_{\alpha_i} f - R^y_{\alpha_0} f\|^2_{H^2_{\widetilde{\zeta}}}$ converges to zero,
i.e., $R^y_{\alpha_i} f$ converges to $R^y_{\alpha_0} f$ in $H^2_{\widetilde{\zeta}}$.
Since we know that for any $f \in \mathcal{H}^{\perp}$, $R^y_{\alpha_i} f \in \mathcal{H}^{\perp}$ for all $i \in \N$,
we conclude that $R^y_{\alpha_0} f \in \mathcal{H}^{\perp}$.
\smallskip

To show that \eqref{lemAlpha0SI} holds, 
we choose a sequence $(\beta _j)_{j \in \mathbb N} \subset \mathbb C _ {-}$  
converging to $\beta_0 \in \mathbb R$ in a neighborhood of infinity such that $\alpha_0 \neq \beta_0$.
According to Lemma \ref{stIdH3}, the structure identity \eqref{StructIdent4} holds for any pair of elements of the sequences
$(\alpha _j)_{j \in \mathbb N}$ and $(\beta _j)_{j \in \mathbb N}$. 
Considering $j \rightarrow \infty$, we obtain the following identity:
\begin{align}
\nonumber
\lim_{j\rightarrow \infty}
\big(
&
\innerProductReg{ R^{y}_{\alpha_j} f}{ g} -
\innerProductReg{ f}{  R^{y}_{\beta_j} g} -
(\alpha_j - \overline{\beta_j})\innerProductReg{  R^{y}_{\alpha_j} f }{  R^{y}_{\beta_j} g}
\big)
=
\\
\label{1234}
&
\lim_{j\rightarrow \infty}
\left(
i(\alpha_j - \overline{\beta_j}) \sum_{l,t=1}^{n}
{
\frac{f(w_j^{(l)}) \overline{g(v_j^{(t)})}}
{dy_k(w_j^{(l)})\overline{dy_k(v_j^{(t)})}}
\frac
{\vartheta[\widetilde{\zeta}](w_j^{(l)}-\tauBa{v_j^{(t)}})}
{\vartheta[\widetilde{\zeta}] (0)E(w_j^{(l)},\tauBa{v_j^{(t)}})}
}
\right)
.
\end{align}
The limit on the right hand side exists due to the analytic continuation to neighborhoods of the fibers over $\alpha_0$ and $\beta_0$.
Since $R^y_{\alpha_j}$ converges to $R^y_{\alpha_0}$ and $R^y_{\beta_j}$ converges to $R^y_{\beta_0}$ in the operator norm, 
the limit on the left hand side coincides with the corresponding value for the operators $R^y_{\alpha_0}$ and $R^y_{\beta_0}$.
\smallskip

To complete the proof, by continuity, 
the limit on the right hand side of \eqref{1234} 
exists also for $\alpha_0 = \beta_0 \in \mathbb R$,
and we take the expression on the right hand side as its definition, 
see Remark \ref{rkBlaBla}.
\end{pf}

As we have gathered all the required preliminary results, 
we may present the proof of the first version of Beurling's theorem.
\\

\begin{pf}[of Theorem \ref{BLth}]
Let $y_1(\cdot)$ and $y_2(\cdot)$, as in the statement, be dividing functions on $X$.
Thus, using Proposition \ref{elisLemma}, the poles of $y_1$ and $y_2$ are real and simple. 
\smallskip

Applying Lemma \ref{MM}, the assumption that $\mathcal{H}$ is invariant under the multiplication operators 
$\mathcal{M}_{\frac{1}{y_1(u)-\overline{\alpha}}}$ and $\mathcal{M}_{\frac{1}{y_2(u)-\overline{\beta}}}$ 
is translated to: $\mathcal{H}^{\perp}$ is invariant under the operators $R_{\alpha}^{y_1}$ and $R_{\beta}^{y_2}$ where $\alpha,\beta \in \mathbb C_{-}$.
By Corollary \ref{stIdH3}, the structure identity automatically holds in $H^2_{\widetilde{\zeta}}$ for all $\alpha,\beta \in \mathbb C _{-}$.
\smallskip

By assumption, the elements of $\mathcal{H}^{\perp}$ 
have analytic extensions with bounded point evaluations to a (connected) neigborhood of the poles of $y_1$ and $y_2$ 
and of the preimages of the singular points of $C$.
Then, using Lemma \ref{MM1}, $\mathcal{H} ^\perp$ is invariant under the bounded operator $R_{\alpha}^{y_k}$ where $\alpha \in \mathbb R$ in a
neighborhood of infinity. 
Furthermore, the structure identity can be extended to $\alpha,\beta \in \mathbb R$ in a neighborhood of infinity.
\smallskip

Combining all the observations above, we can apply Theorem \ref{MainTh} to the orthogonal complement $\mathcal{H}^{\perp}$.
Thus, $\mathcal{H}^{\perp}$ is a reproducing kernel Hilbert space with reproducing kernel of the form
\[
K_{\mathcal{H}^{\perp}}(p,q) = T(p) K_{\zeta}(p,q) T(q)^* - K_{\widetilde{\zeta}}(p,q),
\]
where $T$ is $(\zeta,\widetilde\zeta)$-expansive function for some $\zeta \in T_0$.
By Proposition 3.2, $T$ is contractive on $X_-$, and since $-K_\zeta$
and $-K_{\widetilde\zeta}$ are the reproducing kernels of $H^2_\zeta$
and of $H^2_{\widetilde\zeta}$, respectively,
it follows by the general theory of reproducing kernel spaces, see \cite{dbr2,sarason94},
that ${\mathcal H}^\perp$ is contractively included in $H^2_{\widetilde\zeta}$
and its generalized orthogonal complement is $TH^2_\zeta$
which is likewise contractively included in $H^2_{\widetilde\zeta}$ 
with respect to the range norm.
On the other hand, we know that ${\mathcal H}^\perp$ is isometrically included in $H^2_{\widetilde\zeta}$ therefore so is its generalized orthogonal complement, 
and hence the multiplication by $T$ is an isometry from $H^2_\zeta$ to $H^2_{\widetilde\zeta}$ and so $T$ is inner.
\end{pf}

The natural condition for a Beurling Lax theorem using de Branges structure theorem, see \cite{sarason94}, 
is to assume that the invariant subspace is a contractively included subspace of the Hardy space.
However, in our approach, we note that the structure identity does not automatically hold in a contructively included subspace of the Hardy space
and Lemma \ref{MM} does not hold either when replacing $\mathcal{H}^{\perp}$ by the Brangesian generalized orthogonal complement.
Thus, we add these assumptions in order to state a Beurling's theorem when $T$ is contractive instead of inner.

\begin{theorem}[Beurling's Theorem for Finite bordered Riemann surfaces: version II]
Let $\mathscr{S}$ be a finite bordered Riemann surface where $X$ is its double and
let $y_1$ and $y_2$ be dividing real meromorphic functions generating $\mathcal{M} (X)$.
Let $\mathcal{H}$ be a contractively included subspace of the Hardy space $H^2_{\widetilde{\zeta}}$ on $\mathscr{S}=X_{-}$
such that its generalized orthogonal complement $\mathcal{H}^{[\perp]}$ 
is invariant under $R^{y_1}_\alpha$ and $R^{y_2}_\alpha$ for all $\imagg{\alpha}<0$ 
and the structure identity holds in $\mathcal{H}^{[\perp]}$.
Furthermore, we assume that the elements of ${\mathcal H}^{[\perp]}$ have analytic extensions with bounded point evaluations
to a connected neigborhood of the poles of $y_1$ and $y_2$ and of the preimages of the singular points of $C$.
Then $\mathcal{H}$ is of the form
\[
\mathcal{H} = T H^2_{\zeta},
\]
for some $(\zeta, \widetilde{\zeta})$-contractive mapping $T$.
\end{theorem}

The Beurling's theorems in \cite{MR1070485} and \cite{MR0183883}
assume that $\mathcal{H}$ is invariant under all multiplication operators belonging to some algebra of functions 
(in \cite{MR0183883}, the collection of functions analytic inside a multiply connected domain $R$ and continuous in $\overline {R}$).
Hence, we are motivated to formulate the following version of Beurling's theorem.

\begin{Cy}[Beurling's Theorem for Finite bordered Riemann surfaces: version III]
Let $\mathscr{S}$ be a finite bordered Riemann surface and let $X$ be its double. 
Let $\mathcal{H}$ be a closed subspace of the Hardy space $H^2_{\widetilde{\zeta}}$ on $\mathscr{S}=X_{-}$
corresponding to $\widetilde{\zeta} \in T_0$ such that the following conditions hold:
\begin{enumerate}
\item
$\mathcal{H}$ is invariant under the multiplication operators 
by all functions in the algebra of functions analytic in $\mathscr{S}$ and continuous on ${\mathscr S} \cup \partial{\mathscr S}$.
\item
The elements of $\mathcal{H}^{\perp}$ have analytic extensions with bounded point evaluations in an open neigborhood of a relatively open set
$U \subseteq \partial \mathscr{S}$.
\item 
There exists a pair of dividing functions $y_1$ and $y_2$ such that their poles and the 
pre-images of the singular points of the corresponding algebraic curve belong to $U$.
\end{enumerate}
Then $\mathcal{H}$ is of the form
\[
\mathcal{H} = T H^2_{\zeta},
\]
where $T$ is a $(\zeta, \widetilde{\zeta})$-inner function for some $\zeta \in T_0$.
\end{Cy}
\begin{pf}
The functions
$f_1(u) = {\frac{1}{y_1(u)-\overline{\alpha}}}$
and
$f_2(u) = {\frac{1}{y_2(u)-\overline{\beta}}}$
(where $\alpha,\beta \in \C_{-}$)
are analytic in $\mathscr{S}$ and continuous on ${\mathscr S} \cup \partial{\mathscr S}$.
Therefore, by assumption, $\mathcal{H}$ is invariant under the multiplications by $f_1$ and $f_2$.
It remains to apply Theorem \ref{BLth}.
\end{pf}

We note that every dividing function has poles on every boundary component
and hence a necessary condition is that $U \cap X_j \neq \emptyset$ for $j=0,\ldots,k$.
It is plausible to assume that it is a sufficient condition for 
the existence of the required pair of dividing functions,
see \cite{MR1072300} for related results in the case of planar domains.

\section{Compressed multiplication operators on compact real Riemann surfaces}
\label{secMultOp}

We begin by giving an alternative definition of the compressed multiplication operator $M^y$ 
which is applicable also when $y$ has arbitrary poles.
Let $\mathcal{S}$ be a finite set of points in $X$.
We denote by $\mathcal{M}_{\mathcal{S}}(L_{\widetilde{\zeta}} \otimes \Delta)$ the vector space of 
germs of meromorphic sections of $L_{\widetilde{\zeta}} \otimes \Delta$ 
in a neighborhood of $\mathcal{S}$ and with poles only in $\mathcal{S}$.
Then, we consider the decomposition
\[
\mathcal{M}_{\mathcal{S}}(L_{\widetilde{\zeta}} \otimes \Delta) 
=
\mathcal{M}_{\mathcal{S},+}(L_{\widetilde{\zeta}} \otimes \Delta)
\oplus 
\mathcal{M}_{\mathcal{S},-}(L_{\widetilde{\zeta}} \otimes \Delta),
\]
where $\mathcal{M}_{\mathcal{S},+}(L_{\widetilde{\zeta}} \otimes \Delta)$ denotes the set of germs of holomorphic sections of $L_{\widetilde{\zeta}} \otimes \Delta$ in a neighborhood of $\mathcal{S}$
while $\mathcal{M}_{\mathcal{S},-}(L_{\widetilde{\zeta}} \otimes \Delta)$ is the set of global meromorphic sections of $L_{\widetilde{\zeta}} \otimes \Delta$ with poles in $\mathcal{S}$.
This decomposition is indeed a direct sum decomposition, since by assumption
$L_{\widetilde{\zeta}} \otimes \Delta$ has no global non-zero holomorphic section and therefore
$$\mathcal{M}_{\mathcal{S},+}(L_{\widetilde{\zeta}} \otimes \Delta) \cap \mathcal{M}_{\mathcal{S},-}(L_{\widetilde{\zeta}} \otimes \Delta)=\{0\}.$$
The sum of the two subspaces equals all of ${\mathcal M}_{\mathcal S}(L_{\widetilde\zeta} \otimes \Delta)$ 
since there exist global meromorphic sections with any prescribed principal parts, see \eqref{eqK2Var} below.
For $f \in \mathcal{M}_{\mathcal{S}}(L_{\widetilde{\zeta}} \otimes \Delta)$ we write the corresponding decomposition as $f_{+} + f_{-}$.
We fix a real meromorphic function $y$ such that all the poles of $y$ belong to $\mathcal{S}$ and we define
\begin{align*}
M^y:  \mathcal{M}_{\mathcal{S},+} & \mapsto \mathcal{M}_{\mathcal{S},+} \\
f(p)   & \mapsto  (y(p)f(p))_+.
\end{align*}

In other words, $M^y$ is given by
$$(M^y f )(u)= y(u)f(u)+c_f(u).$$
Here
$f$ is a holomorphic section of $L_{\widetilde{\zeta}} \otimes \Delta $ in a neighborhood of $\mathcal{S}$
and $c_f$ is the unique global meromorphic section with 
divisor of poles contained in the divisor of poles of $y$, 
such that $y(u)f(u)+c_f(u)$ is analytic at the poles of $y$.
\smallskip

It follows immediately from the properties of the Cauchy kernel that in the case when $y$ has only simple poles
the new definition of $M^y$ coincides with the definition of the model operator in \eqref{m_y}.
It also follows that we can write the model operator when $y$ has also non-simple poles as
\begin{align}
\nonumber
M^{y}f(u)
= &
y(u)f(u) +
\\ &
i
\sum_{m=1}^{n}
\sum_{l=1}^{s^{m}}
\sum_{j=l}^{s^{m}}
a_{m,-j}
\frac{ f^{(j-l)}(p^{(m)})}{ (j-l)!(l-1)!}
K_{\widetilde{\zeta}}^{(0,l-1)}(u,\tauBa{p^{(m)}})
,
\label{MyMatrix}
\end{align}
where the set of the poles of $y$ and their orders are given by $\left( p^{(m)} \right)_{m=1}^{n}$ and $\left( s^{m} \right)_{m=1}^{n}$, respectively, 
and where $a_{m,-j}$ is the $-j$-th Laurent coefficient of $y$ at $p^{(m)}$ (see also \cite[Equation 4.21]{av3}).
Here (and in the rest of this section) we use the notation
\begin{equation}
\label{eqK2Var}
K^{(j,l)}_{\widetilde{\zeta}}(u,w) 
\defEq
\frac{\partial^{j+l}}{\partial u ^{j} \partial \overline{w} ^{l}}K^{(j,l)}_{\widetilde{\zeta}}(u,w),
\end{equation}
where the derivatives are computed with respect to some local coordinates 
(in coordinate free terms, this is a higher order connection for the corresponding line bundles).
In particular, $K^{(0,l)}_{\widetilde\zeta}(\cdot,w)$ is a global meromorphic section of $L_{\widetilde\zeta}\otimes \Delta$ which has
a pole of order $l+1$ at $\tauBa{w}$ with the principal part,
in terms of a local coordinate $t$ centered at $\tauBa{w}$ used to compute the derivative in \eqref{eqK2Var}, 
given by $-\frac{l!}{i t^{l+1}}$.
\smallskip

In order to prove Theorem \ref{lemCollCodI} below, we first present some preliminary results.
The next important result appeared partially in \cite[Theorem 4.10]{av3}.
Here, we present a more general statement with a comprehensive and an alternative proof.

\begin{theorem}
\label{rkModelOpAlg}
Let $\mathcal{S}$ be a finite set of points on $X$ and let $\mathcal A$ be the algebra of meromorphic functions on $X$ whose poles are contained in $\mathcal{S}$.
The mapping $y \mapsto M^{y}$, where $y$ is a meromorphic function on $X$ whose poles are contained in $\mathcal{S}$,
is an algebra homomorphism from $\mathcal{A}$ into the algebra of linear operators on the vector space $\mathcal{M}_{\mathcal{S},+}(L_{\widetilde{\zeta}} \otimes \Delta)$.
\smallskip

In other words, for a two-variables polynomial $g(x_1,x_2)$ and pair of meromorphic functions
$y_1$ and $y_2$ with poles in $\mathcal{S}$, the model operator satisfies 
$$M^{g(y_1,y_2)} = g(M^{y_1},M^{y_2}).$$
In particular, $M^{y_1 \, y_2} = M^{y_1}M^{y_2}$ and $M^{y_1 + y_2} = M^{y_1} + M^{y_2}$.
\end{theorem}
\begin{pf}[of Theorem \ref{rkModelOpAlg}]
Let $f=f_{+}$ be an element in $\mathcal{M}_{\mathcal{S},+}(L_{\widetilde{\zeta}} \otimes \Delta)$, then, by definition, $M^{y_1 + y_2}f = M^{y_1}f+M^{ y_2}f$.
We set  $M^{y_1}M^{y_2} f = h_+ = (h_+ + h_-)_+ $ and we show that $h_+ = M^{y_1 \, y_2} f$.
We also set $y_2 \cdot f_+ = g_+ + g_-$ and then we have
\begin{align*}
h_+(p) = & y_1 (p) g_+(p) - h_-(p) 
\\ = & 
y_1(p)\big(y_2(p) f(p) - g_-(p)\big)- h_-(p),
\end{align*}
and hence
\[
y_1(p)y_2(p)f(p) = h_+(p) + \big( h_-(p)+y_1(p)g_-(p) \big).
\]
However, $h_-(p)+y_1(p)g_-(p)$ is a global meromorphic section with poles in $\mathcal{S}$,
so $\big(y_1(p)y_2(p)f(p)\big)_+  = h_+(p)$ and $M^{y_1 y_2}f = h_+$ follows.
\end{pf}

The main result in this section is given in the next theorem.

\begin{theorem}
\label{lemCollCodI}
Let $y_1$ and $y_2$ be a pair of meromorphic functions on a compact real Riemann surface $X$.
Let us assume that the structure identity, given in \eqref{StructIdent4}, 
holds in a reproducing kernel Hilbert space $\mathcal X$
for a pair of meromorphic functions $y_1$ and $y_2$.
Then, the structure identity holds for all functions in the 
algebra of meromorphic functions generated by $y_1$ and $y_2$.
\end{theorem}

We start with the following technical result required later in the proof of Proposition \ref{myGeneral1}.

\begin{lem}
\label{lemLimitZero}
Let $a,b,d \in \mathbb N$ such that $0\leq a,b < d$ and let $c_0,\ldots,c_d$ 
be a sequence of real numbers such that $c_0 \neq 0$. 
Then, the equality
\begin{equation}
\label{eqLimitDelta}
\lim_{x,y \rightarrow 0}
\sum_{j=0}^{d-1-a}
\frac{
c_{d-(j+a+1)} }{ b! j !}
\frac{\partial^{j+b}}{\partial x ^j \partial y ^b} 
\frac
{
\sum_{q=0}^{d-1} c_{q}
\sum_{t=1}^{d-q} x^{q+t-1} y^{d-t}
} 
{\sum_{p=0}^{d}c_{p}x^{p}    \sum_{p=0}^{d}c_{p}y^{p}}
= 
\delta_{b,a}
\end{equation}
holds, where $\delta$ stands for the Kronecker delta.
\end{lem}
\begin{pf}
Along the proof, we use the notations $f(x) \defEq \sum_{p=0}^{d}c_{p}x^{p}$, 
\[
g[n] \defEq \frac{1}{n!} \restr{ \frac{d^n}{dx^n} \frac{1}{f(x)} }{x=0}
    \qquad
    {\rm and}
    \qquad
h[m] \defEq \sum_{n=0}^{m} g[n] g[m-n].
\]
According to Leibniz product rule, the following identity holds
\begin{align}
\nonumber
\lim_{x,y \rightarrow 0}
\frac{\partial^{j+b}}{\partial  x ^j\partial y^b} 
\frac
{
\sum_{q=0}^{d-1} c_{q}
\sum_{t=1}^{d-q} x^{q+t-1} y^{d-t}
} 
{f(x) \, f(y)}
& 
= 
\sum_{k=0}^{j}  \sum_{l=0}^{b} 
\frac{j! g[k]}{(j-k)!}
\frac{b! g[l]}{(b-l)!}
\sum_{q=0}^{d-1} c_{q}
\times \\ 
&
\lim_{x,y \rightarrow 0} 
\frac{\partial ^{(j-k)+(b-l)}}{\partial  x ^ {j-k} \partial  y ^{b-l}}
\sum_{t=1}^{d-q} x^{q+t-1} y^{d-t}
\label{eqLeib}
.
\end{align}
The two-variable polynomial in the numerator on the LHS of \eqref{eqLeib}
has zero coefficients for all monomials of combined degree less than $d-1$,
and the same is therefore true after multiplying it by 
the Taylor series of $\frac{1}{f(x)f(y)}$.
Hence, all mixed derivatives of orders $j$ and $b$ such that $j+b<d-1$ are zero.
Since the outer sum on the LHS of \eqref{eqLimitDelta} is summed up to $d-1-a$, 
the LHS of \eqref{eqLimitDelta} is zero whenever $b<a$.
\smallskip

We move to consider the case where $a \leq b$. 
First, we note that
\begin{align*}
\nonumber
\sum_{q=0}^{d-1} c_{q}
&
\lim_{x,y \rightarrow 0}
\frac{\partial ^{(j-k)+(b-l)}}{\partial x ^{j-k} \partial  y ^{b-l}}
\sum_{t=1}^{d-q} x^{q+t-1} y^{d-t}
= \\ &
\begin{cases}
(j-k)! (b-l)! c_{(j-k)+(b-l)-(d-1)}
&
d-1 \leq (j-k)+(b-l),
\\
\quad
0
&
\it{otherwise}.
\end{cases}
\end{align*}
Here we repeat the same argument (as used in the $b<a$ case),
we note that if $j<d-b-1$, the mixed derivatives on the LHS of \eqref{eqLeib} are again zero.
Whenever $d-1 \leq b+j$, the RHS of \eqref{eqLeib} becomes
\begin{align*}
RHS
= &
\lim_{x,y \rightarrow 0}
\sum_{k = 0}^{j}
\sum_{l = 0}^{b}
\frac{j! b! g[k] g[l]}{(j-k)!(b-l)!}
\frac{\partial ^{(j-k)+(b-l)}}{\partial  x ^ {j-k} \partial  y ^{b-l}}
\sum_{q=0}^{d-1} c_{q}
\sum_{t=1}^{d-q} x^{q+t-1} y^{d-t}
\\ = &
\sum_{s = 0}^{j+b-d+1}
\sum_{l = 0}^{s}
j! g[s-l]
b! g[l]
c_{j+b-d+1 - s}
.
\end{align*}
The last equality is the result of setting $s \defEq k+l$.
It follows from the last observation that
the LHS of \eqref{eqLimitDelta} can be rewritten as follows
(we just set $j'=d-1-a-j$, $b^\prime = b-a$ and $m=b'-j'$):
\begin{align}
LHS =
&
\nonumber
\sum_{j=d-1-b}^{d-1-a}
c_{d-j-a-1} 
\sum_{s = 0}^{j+b-d+1}
c_{j+b-d+1 - s}
\sum_{l = 0}^{s}
g[s-l] g[l]
\\ = &
\nonumber
\sum_{j'=0}^{b-a}
\sum_{s = 0}^{b-a-j'}
c_{b-a-j'-s}
c_{j'}
h[s]
\\ = &
\sum_{s=0}^{b^\prime}
h[s]
\sum_{m=s}^{b^\prime}
c_{m-s}
c_{b^\prime-m}
\label{eqZero}
.
\end{align}
We note that $h[s]$ is just the $s$-Taylor coefficient at zero of $\frac{1}{f(x)}\frac{1}{f(x)}$ (recall that $c_0 \neq 0$).
On the other hand, $\sum_{m=s}^{b^\prime} c_{m-s} c_{b^\prime-m}$ is the $(b^\prime-s)$-Taylor coefficient at zero of $f(x) f(x)$.
Hence, \eqref{eqZero} is the $b^\prime$-Taylor coefficient at zero of the multiplication of $f(x) f(x)$ and $\frac{1}{f(x)}\frac{1}{f(x)}$.
The $b^\prime$-Taylor coefficient of the function identically equal to one is clearly equal to zero whenever $b^\prime>0$, 
or equivalently, when $b>a$ and equal to $1$ if $b^{\prime}=0$ or $b=a$.
Hence the claim follows.
\end{pf}

In the next proposition we will show that the structure identity is equivalent to an appropriate colligation condition for the model operator.
To identify exactly the matrix appearing on the right hand side of the colligation condition, we will need the following two technical remarks.

\begin{Rk}
\label{T0}
For each component $X_i$ where $0 \leq i \leq k-1$, one may attach a sign $(-1)^{\mu_i}$ 
with $\mu_i = 0$ or $\mu_i = 1$ as a sign of an appropriate differential
(if in addition $X$ is of dividing type and $X_\R$ is given the natural orientation, then $\mu_i=0$ for all $i$).
Then, for any $f \in \L_{\widetilde{\zeta}} \otimes \Delta$, where $\widetilde{\zeta} \in T_{\nu}$ 
and for any $p \in X_i$ we have $ f(\tauBa{p} = (-1)^{\mu_i+\nu_i} f(p)$ (we set $\nu_0=0$).
Furthermore, ${\rm sign} \, K_{\widetilde{\zeta}}(q,q) = (-1)^{\mu_i+\nu_i}$ where $q$ is 
close to $X_i$ on the positive side with respect to the chosen orientation.
For more details we refer to \cite[Section 2]{av3}.
\end{Rk}

\begin{Rk}
Let $y(p)$ be a real meromorphic function with $n$ poles of order $s^1,\ldots,s^n$.
Then, the poles of $y$ are either real or appear as conjugate pairs 
(recall that our compact real Riemann surface $X$ is not necessarily of dividing type and $y$ is not necessarily dividing). 
We denote by $\left(p^{(r)}\right)_{r=1}^{\mathbf{r}}$ the real poles, 
and by $\left(p^{(r)}\right)_{r=\mathbf{r}+1}^{\mathbf{r}+2\mathbf{m}}$ the $\mathbf{m}$ pairs of non-real poles ($n = \mathbf{r} + 2\mathbf{m}$).
\end{Rk}

\begin{prop}
\label{myGeneral1}
Let $y$ be a real meromorphic function on a compact real Riemann surface $X$.
Let $\mathcal X$ be a reproducing kernel Hilbert space of section of 
$L_{\widetilde{\zeta}}\otimes\Delta$ which is invariant under the operators
$R^y_\alpha$, $R^y_\beta$ and $M^y$ for fixed $\alpha,\beta \in \C$.
Then the structure identity \eqref{StructIdent4} is equivalent to the colligation condition \eqref{collCond} for the operator $M^y$.
Here $\Phi = \Phi_y$ is the evaluation operator that maps a section $f$ holomorphic at the poles of $y$
to the values of $f$ and its Taylor coefficients of order up to $s^r-1$ at $p^{(r)}$ where $r=1,\ldots,n$. 
$\sigma=\sigma_y$ is defined by
\begin{equation}
\label{eqDefSigmaY}
\sigma_y =
 \begin{psmallmatrix}
  \mathbf{P}^1		& 				&  				&		&			&		$		$\\
  			&	\ddots 		&  				&  		&			&		$		$\\
  			&			 	& 	\mathbf{P}^{\mathbf{r}}		    &  		&			&		$		$\\
  			&			 	& 				& 0		& \mathbf{P}^{\mathbf{r}+2}	&		&		&\\
  			&		 		&  				& \mathbf{P}^{\mathbf{r}+1}& 0		&		$		$\\
  			&		 		&  				&  		&			&\ddots	$		$\\
  			&		 		&  				&  		&			&		&0	& \mathbf{P}^{\mathbf{r}+2\mathbf{m}} \\
  			&		 		&  				&  		&			&		& \mathbf{P}^{\mathbf{r}+2\mathbf{m}-1}& 0\\
 \end{psmallmatrix},
\end{equation}
where $\mathbf{P}^r$ are upper-skew-triangular Hankel matrices 
(Hankel matrices with zero entries below the main skew-diagonal)
of sizes $s^r$ 
where the $(\gamma,\delta)$-entry (where $\gamma,\delta = 0,\ldots,s^r-1$ and $\gamma + \delta < s^r$) 
is equal to minus the corresponding Laurent coefficients 
$a_{r,-(\gamma + \delta + 1)}$ of $y$ at the pole $p^{(r)}$, see \eqref{eqYexpan},
with $\mathbf{P}^1,\ldots,\mathbf{P}^{\mathrm r}$ further multiplied by $(-1)^{\mu_i+\nu_i}$,
see Remark \ref{T0}, depending on the component $X_i$ containing the corresponding real pole.
\end{prop}
\begin{pf}
First, using the notation in \eqref{eqK2Var}, let us consider the following version of the collection formula 
\footnote{
We note that this generalization of the collection formula 
is written in terms of the Hermitian Cauchy kernel and 
not in terms of the non-Hermitian kernel as in Definition \ref{defCollFor}
}
\begin{align}
\nonumber
i \big( y(v) -  y(\tauBa{w}) \big)& K_{\widetilde{\zeta}}(v,w) 
= 
\\ &
\sum_{r=1}^{n}
\sum_{\delta,\gamma = 0} ^{s^r-1} 
\frac{K_{\widetilde{\zeta}}^{(0,\gamma)}(v,\tauBa{p^{(r)}})}{\gamma!}
\mathbf{A}_{r;\gamma,\delta} 
\frac{K_{\widetilde{\zeta}}^{(\delta,0)}(p^{(r)},w)}{\delta!}
.
\label{eqCollFor}
\end{align}
Here ${\mathbf A} = \diag_{r=1,\ldots,n} \left[{\mathbf A}_{r;\gamma,\delta}\right]_{\gamma,\delta=0,\ldots,s^r-1}$
is a block diagonal matrix (we will also write ${\mathbf A}_y$ to indicate the dependence on the meromorphic
function $y$ whenever needed for clarity) with entries defined by \eqref{lblAmatrix} below. \smallskip

To see that \eqref{eqCollFor} holds, we fix $w=w_0 \in X$, such that $w_0 \notin \mathcal{P}(y)$ and $\overline{y(w_0)}$ has $n$ distinct pre-images.
Then, both sides of \eqref{eqCollFor} are meromorphic sections of $L_{\widetilde{\zeta}}\otimes\Delta$ in $v$.
Furthermore, both sides of \eqref{eqCollFor} have poles at $(p^{(r)})_{r=1}^n$ of orders $(s^{r})_{r=1}^n$, respectively,
and the additional singularity on the left hand side at $v=\tauBa{w_0}$ is removable.
We set the entries of $\mathbf{A}$ to be the negative Laurent coefficients of $y$
at $(p^{(r)})_{r=1}^n$ such that both sides of \eqref{eqCollFor} share the same principal parts.
More precisely, we use a local parameter $t_r$ centered at $p^{(r)}$.
We assume for later purposes that $t_r$ is real on $X_i$ and compatible with the chosen orientation 
if $p^{(r)} \in X_i$ is a real pole,
and that the local coordinates centered at complex conjugate poles are complex conjugate to each other.
Then, the expansion of $y(\cdot)$ at $p^{(r)}$ is given by
\begin{equation}
\label{eqYexpan}
y(v) = a_{r,-s^r}t_r(v)^{-s^r} +\ldots+ a_{r,-1}t_r(v)^{-1}+\ldots,
\end{equation}
while the expansion of the Cauchy kernel $K_{\widetilde{\zeta}}(p^{(r)},w_0)$ is
\[
K_{\widetilde{\zeta}}(v,w_0)
=
\sum_{j=0}^{\infty} 
\frac{1}{j!} K^{(j,0)}_{\widetilde{\zeta}}(p^{(r)},w_0) t_r(v)^j.
\]
Hence, the Laurent coefficient of $t_r(v)^{-l}$ (where $0 < l \leq s^r$) of the LHS of \eqref{eqCollFor} is
\begin{equation}
\label{eqLRes}
i \sum_{j=0}^{s^r-l} a_{r,-(l + j)} \frac{1}{j!} K^{(j,0)}_{\widetilde{\zeta}}(p^{(r)},w_0)
.
\end{equation}

To get the coefficient of $t_r(v)^{-l}$ on the RHS of \eqref{eqCollFor}, 
we recall that $K_{\widetilde\zeta}^{(0,\gamma)}$
has a pole of order $\gamma+1$ at $\tauBa{p^{(r)}}$ with principal part
$- \frac{\gamma!}{i t_r(v)^{\gamma+1}}$ (see remark above).
We therefore set $\gamma = l-1$ and obtain that 
the coefficient of $t_r(v)^{-l}$ on the RHS of \eqref{eqCollFor} is
$$
i \sum_{\delta=0}^{s^r-1} {\mathbf{A}}_{r;\gamma,\delta}
\frac{K_{\widetilde\zeta}^{(\delta,0)}(p^{(r)},w_0)}{\delta!}.
$$
To match this to \eqref{eqLRes}, we take $j=\delta$. We see that we need $ {\mathbf{A}}_{r;\gamma,\delta} = 0$ for $j>s^r-l$, i.e., for $\gamma+\delta>s^r-1$, i.e., the $s^r \times s^r$ matrix $\left[{\mathbf{A}}_{r;\gamma,\delta}\right]_{\gamma,\delta=0}^{s^r-1}$ is skew-upper-triangular 
where for the elements on or above the main skew diagonal we obtain 
$i \frac{1}{\delta!} a_{r-(\gamma+\delta+1)} =  \frac{i}{\delta!} {\mathbf{A}}_{r;\gamma,\delta}$ 
yielding finally
\begin{equation}
\label{lblAmatrix}
{\mathbf{A}}_{r;\gamma,\delta} 
=
\begin{cases}
a_{r,-(\gamma+\delta+1)},
& 
\gamma+\delta < s^r,
\\
0, & {\text otherwise}.
\end{cases}
\end{equation}
Since, by assumption, there are no non-zero global holomorphic sections, 
we may conclude that \eqref{eqCollFor} holds for any $v \in X$.
\smallskip

It is convenient to use the following notations (see Lemma \ref{mainThLem})
\begin{equation*}
F = (M^{y} - \alpha I) ^{-1} f
\qquad {\rm and} \qquad
G  = (M^{y} - \beta I) ^{-1} g.
\end{equation*}
Recall that using \eqref{StructIdent5}, $F$ and $G$ are the images of resolvent operators acting on $f$ and $g$, respectively.
Thus, the left hand side of the structure identity \eqref{StructIdent4} can then be rewritten as
\begin{align}
\label{StructIdent7}
\innerProductReg{(M^{y}-  \alpha I)^{-1} f}{g} & -
\innerProductReg{f}{(M^{y}-\beta I)^{-1}g} -
\nonumber
\\ &
(\alpha - \overline{\beta})
\innerProductReg{ (M^{y}-\alpha I)^{-1}f}{(M^{y}- \beta I)^{-1}g}
\nonumber
\\ = &
\innerProductReg{
(M^{y}-\beta I)^{-*}
\left( (M^{y})^* - M^{y} \right)
(M^{y}-\alpha I)^{-1}f}{g}
\nonumber
\\ = &
\innerProductReg{
\left( (M^{y})^* - M^{y} \right) F}{G}.
\end{align}
We then substitute \eqref{eqCollFor} (with $v=\tauBa{v^{(t)}}$ and $w=\tauBa{w^{(l)}}$) and \eqref{StructIdent7} 
on the LHS and RHS of structure identity \eqref{StructIdent4} (multiplied by $(-1)$), respectively,
and we get the expression
\begin{align}
\nonumber
\innerProductReg{(M^{y} - (M^{y})^*) F}{G}
= & 
- (\alpha - \overline{\beta})
\sum_{l,t=1}^{n}
\frac{f(w^{(l)}) }{dy(w^{(l)})}
K_{\widetilde{\zeta}}(\tauBa{v^{(t)}},\tauBa{w^{(l)}})
\frac{\overline{g(v^{(t)})}}{\overline{dy(v^{(t)})}}
\\ = &
\nonumber
- i
\sum_{r=1}^{n}
\sum_{\delta,\gamma = 0} ^{s^r-1}
\left(
\sum_{l=1}^{n}
\frac{f(w^{(l)})}{dy(w^{(l)})}
K_{\widetilde{\zeta}}^{(\delta,0)}(p^{(r)},\tauBa{w^{(l)}})
\right)
\\ 
&
\times 
\frac{\mathbf{A}_{r;\gamma,\delta}}{\gamma! \delta!}
\left(
\sum_{t=1}^{n}
K_{\widetilde{\zeta}}^{(0,\gamma)}(\tauBa{v^{(t)}},\tauBa{p^{(r)}})
\frac{\overline{g(v^{(t)})}}
{\overline{dy(v^{(t)})}}
\right)
\label{eq95}
.
\end{align}
Noting that $y(w ^{(l)})=\alpha$, the evaluation of $f$ at $w ^{(l)}$ in terms of $F$
(an equivalent formula ties $G$ and $g$) may be written as (see \eqref{MyMatrix}):
\begin{align*}
f(w^{(l)}) = & (M^y - \alpha) F(w^{(l)})
\\ = & i\sum_{m=1}^{n}\sum_{k=1}^{s^{m}} \sum_{j=k}^{s^{m}} 
a_{m,-j} \frac{ F^{(j-k)}(p^{(m)})}{ (j-k)!(k-1)!} K_{\widetilde{\zeta}}^{(0,k-1)}(w^{(l)},\tauBa{p^{(m)}}).
\end{align*}
We write the first inner sum on the right hand side of \eqref{eq95}, as follows:
\begin{align}
\nonumber
\sum_{l=1}^{n}
\frac{f(w^{(l)})}{dy(w^{(l)})}
K_{\widetilde{\zeta}}^{(\delta,0)}(p^{(r)},\tauBa{w^{(l)}})
& = i
\sum_{m=1}^{n}
\sum_{k=1}^{s^{m}}
\sum_{j=k}^{s^{m}}
a_{m,-j}
\frac{ F^{(j-k)}(p^{(m)})}{ (j-k)!(k-1)!}
\\ & \times
\sum_{l=1}^{n}
\frac{
K_{\widetilde{\zeta}}^{(\delta,0)}(p^{(r)},\tauBa{w^{(l)}})
K_{\widetilde{\zeta}}^{(0,k-1)}(w^{(l)},\tauBa{p^{(m)}})
}{dy(w^{(l)})}
\label{eq34}
.
\end{align}
Differentiating the collection formula for the function 
$\frac{1}{y(\cdot)-\alpha}$, with simple poles $w^{(l)}$, $\delta$ times with respect to $v$
and $k-1$ times with respect to $w$ leads to
\begin{align}
\nonumber
\sum_{l=1}^{n}
&
\frac
{K_{\widetilde{\zeta}}^{(\delta,0)}(v,\tauBa{w^{(l)}})K_{\widetilde{\zeta}}^{(0,k-1)}(w^{(l)},\tauBa{w})}
{dy(w^{(l)})}
=  
\\ & =
i
\frac{\partial ^{\delta+k-1}}{\partial  v ^{\delta} \partial w ^{k-1}} 
\left(  \frac{1}{y(v) - \alpha} - \frac{1}{y(w) - \alpha} \right) 
K_{\widetilde{\zeta}}( v , \tauBa{w})
\label{eqRepad45}
.
\end{align}
Therefore, the inner sum in \eqref{eq34} vanishes whenever $p^{(r)} \neq p^{(m)}$.
Furthermore, we change the summation indices to be $j^\prime = j-k$ and $k^\prime=k-1$ and the RHS of \eqref{eq34} becomes
\begin{align}
\label{eqRepad}
i \sum_{j^\prime = 0}^{s^{r}-1}
\frac{F^{(j^\prime)}(p^{(r)})}{j^\prime!}
&
\sum_{k^\prime=0}^{s^{r}-1-j^\prime}
\frac{a_{r,-(j^\prime+1+k^\prime)} }{ k^\prime!}
\\ & \nonumber
\sum_{l=1}^{n}
\frac{
K_{\widetilde{\zeta}}^{(\delta,0)}(p^{(r)},\tauBa{w^{(l)}})
K_{\widetilde{\zeta}}^{(0,k^\prime)}(w^{(l)},\tauBa{p^{(r)}})}
{dy(w^{(l)})}
.
\end{align}
Evaluating the RHS of \eqref{eqRepad45} at $\tauBa{w}=v=p^{(r)}$, 
using the local parameter $t_r$ centered at $p^{(r)}$ (notice that since $\gamma,k^\prime \leq s^r-1$ 
we can omit higher order terms in $t(v)$ and $t(w)$), leads to
\begin{align*}
\frac{\partial ^{\delta+k^{\prime}}}{\partial t(v) ^{\delta}\partial t(w) ^{k^{\prime}}} 
\frac{1}{t_r(w)-t_r(v)}
\bigg( 
& 
\frac{t_r(v)^{s^{r}}} {a_{r,-s^{r}}+\ldots+(a_{r,0}- \alpha) t_r(v)^{s^{r}}} 
-
\\ &
\frac{t_r(w)^{s^{r}}}{a_{r,-s^{r}}+\ldots+ (a_{r,0}- \alpha) t_r(w)^{s^{r}}} 
\bigg)
_{\substack{t(v)=0\\t(w)=0}}.
\end{align*}
Thus, \eqref{eqRepad45} becomes 
\begin{align}
\nonumber
\sum_{l=1}^{n}
&
\frac{
K_{\widetilde{\zeta}}^{(\delta,0)}(p^{(r)},\tauBa{w^{(l)}})
K_{\widetilde{\zeta}}^{(0,k^\prime)}(w^{(l)},\tauBa{p^{(r)}})}
{dy(w^{(l)})}
= 
- i
\bigg(
\frac{\partial ^{\delta+k^\prime}}{\partial t(v) ^{\delta} \partial t(w)^{k^\prime}} 
\\ &
\frac
{
\sum_{\beta=1}^{s^{r}} a_{r,-\beta}
\sum_{\eta=0}^{\beta-1} t_r(w)^{s^{r}-\beta+\eta} t_r(v)^{s^{r}-1-\eta} } 
{(a_{r,-s^{r}}+\ldots+(a_{r,0}- \alpha) t_r(w)^{s^{r}})(a_{r,-s^{r}}+\ldots+(a_{r,0}- \alpha) t_r(v)^{s^{r}})}
\bigg)
_{\substack{t(v)=0\\t(w)=0}}
.
\label{eqCollCondQ}
\end{align}
Substituting \eqref{eqCollCondQ} in the RHS of \eqref{eqRepad} leads to
\begin{align}
\nonumber
\sum_{j^\prime=0}^{s^{r}-1} &
\frac{F^{(j^\prime)}(p^{(r)})}{j^\prime !}
\sum_{k^\prime=0}^{s^{r}-1-j^\prime}
\frac{a_{r,-(j^\prime+k^\prime+1)} }{k^\prime!}
\bigg(
\frac{\partial ^{\delta+k^\prime}}{\partial t(v) ^{\delta} \partial t(w)^{k^\prime}} 
\\ &
\frac
{
\sum_{\beta=1}^{s^{r}} a_{r,-\beta}
\sum_{\eta=0}^{\beta-1} t_r(w)^{s^{r}-\beta+\eta} t_r(v)^{s^{r}-1-\eta} } 
{(a_{r,-s^{r}}+\ldots+(a_{r,0}- \alpha) t_r(w)^{s^{r}})(a_{r,-s^{r}}+\ldots+(a_{r,0}- \alpha) t_r(v)^{s^{r}})}
\bigg)
_{\substack{t(v)=0\\t(w)=0}}
,
\label{eqRepadQ}
\end{align}
and, by applying Lemma \ref{lemLimitZero} to the inner sum in \eqref{eqRepadQ} (with $d=s^r$, $a=j'$, $b=\delta$, $c_\cdot = a_{r,\cdot-s^r}$),
we conclude that \eqref{eqRepad} is just equal to $F^{(\delta)}(p^{(r)})$.
Similarly, for the second inner sum in \eqref{eq95}, we have
\begin{align*}
\sum_{t=1}^{n}
K_{\widetilde\zeta}^{(0,\gamma)}(\tauBa{v^{(t)}},\tauBa{p^{(r)}})
\frac{\overline{g(v^{(t)})}}{\overline{dy(v^{(t)})}}
= &
\overline{
\sum_{t=1}^{n}
K_{\widetilde\zeta}^{(\gamma,0)}(\tauBa{p^{(r)}},\tauBa{v^{(t)}})
\frac{g(v^{(t)})}{dy(v^{(t)})}
}
\\ = &
\overline{G^{(\gamma)}(\tauBa{p^{(r)}})}
.
\end{align*}
Finally, we move back to evaluate \eqref{eq95} and we summarize:
\begin{align*}
\innerProductReg{
(M^{y}-(M^{y})^*) F}{G}
= &
-
i
\sum_{r =1}^{n}
\sum_{\gamma = 0}^{s^{r}}
\sum_{\delta = 0}^{s^{r}}
\frac{\overline{G^{(\gamma)}(\tauBa{p^{(r)}})}}{\gamma!}
\mathbf{A}_{r; \gamma, \delta} 
\frac{F^{(\delta)}(p^{(r)})}{\delta!}
\\ = &
-
i
\Phi_y (G) ^*\, \Per{y} \,  \mathbf{A}_y \, \Phi_y (F)
\\ = &
i\Phi_y (G)^* \, \sigma_y \, \Phi_y (F)
.
\end{align*}
Here $\Phi_y$ returns for each pole $p^{(r)}$ of $y$ the first $s^r-1$ Taylor coefficients,
more precisely, 
\begin{equation}
\label{defPhiYOp}
\Phi_y \colon f \mapsto \col_{r=1,\ldots,n} \col_{\delta=0,\ldots,s^r-1} \frac{f^{(\delta)}(p^{(r)})}{\delta!}.
\end{equation}
Moreover, the matrix $\sigma_y$ is equal to $-\mathbf{A}_y$ (defined in \eqref{lblAmatrix}) up to the blocks 
permutation corresponding to conjugate poles of a real meromorphic function $y$ 
and up to the signs associated to each boundary component, see Remark \ref{T0}. 
We denote the corresponding signed permutation matrix by
\begin{equation}
\label{eqPerY} 
\Per{y}\defEq \left( \begin{smallmatrix} R_y & o \\ 0 & C_y \end{smallmatrix} \right),
\end{equation}
where
\begin{equation*}
R_y
=
\diag_{j=1,\ldots,\mathbf{r}}((-1)^{\mu_{i_j}+\nu_{i_j}}I_{s^j}), 
\quad
C_y=
\begin{psmallmatrix} 
 {0}                    & I_{s^{\mathbf{r}+2}}& & & \\
 {I_{s^{\mathbf{r}+1}}} & {0}   &		$		$\\
 				        &  		&			&\ddots	$		$\\
  				        &  		&			&		&0	& I_{s^{\mathbf{r}+2\mathbf{m}}} \\
  				        &  		&			&		& I_{s^{\mathbf{r}+2\mathbf{m}-1}}& 0\\
\end{psmallmatrix} 
,
\end{equation*}
and where $i_j$ is defined by $p^{(j)} \in X_{i_j}$.
\end{pf}

We note that the matrix $\sigma_y$ is selfadjoint. It is an immediate consequence of the assumption
that $y$ is real (and hence the Laurent coefficients of $y$, 
with appropriate choice of local coordinates,
are real at real points and conjugate to each other at complex conjugate points).
\smallskip

In the upcoming proofs we extensively use the following notation.

\begin{Notation}
\label{notatMatrix1}
We denote by 
$\mathcal P(y_k)$ the set of poles of $y_k$,
we set $n_k = \Card{\mathcal{P}(y_k)}$ for $k=1,2$ and 
$n = \Card{\mathcal{P}(y_1) \cup \mathcal{P}(y_2)}$.
The set of poles is denoted by $\left(p^{(m)}\right)_{m=1}^{n}$
and we denote by $s_{k}^{m}$ the order of $y_k$ at a pole $p^{(m)}$
(we set $s_{k}^{m}=0$ when $p^{(m)} \notin \mathcal{P}(y_k)$), for $k=1,2$ and $m=1,\ldots,n$.
Furthermore, we denote the cumulative order of a pole by $s^{m}=s_1^{m}+s_2^{m}$ and the maximum by $\hat{s}^{m} = \max(s_1^{m},s_2^{m})$.
\smallskip

The operators $\Phi_{11},\Phi_{22},\Phi_{12}$ and $\widehat{\Phi}_{12}$ 
return the Taylor coefficients of a section $f$ of $L_{\widetilde{\zeta}}\otimes \Delta$ at the poles of 
$\mathcal P (y_1) \cap \mathcal P (y_2)^c$, $\mathcal P (y_2) \cap \mathcal P (y_1)^c,\mathcal P (y_1) \cap \mathcal P (y_2)$ and $\mathcal P (y_1) \cap \mathcal P (y_2)$, 
up to orders of $s^{m}_1-1$, $s^{m}_2-1$, $s^{m}-1$ and $\widehat{s}^{m}-1$, respectively.
$\Phi_1^{p^{(m)}}$, $\Phi_2^{p^{(m)}}$ and $\Phi_{12}^{p^{(m)}}$ return the first $s_1^m$, $s_2^m$ and $s^m$ Taylor coefficients of $f$ at $p^{(m)}$, respectively.
\smallskip

The block matrix $\mathbf{K}_{\widetilde{\zeta}}(P_{1},P_2)$ is defined by 
\[
\left[\left[
\frac{K_{\widetilde{\zeta}}^{(k_1,k_2)}(p^{(m_1)}_1,\tauBa{p^{(m_2)}_2})}{k_1!k_2!}
\right]_{k_1,k_2}\right]_{m_1,m_2},
\quad
p^{(m_1)}_1 \in \mathcal P (y_1)  \, \, {\rm  and  } \,\,  p^{(m_2)}_2 \in \mathcal P (y_2)
,
\]
where $m_1=1,\ldots,n_1$, $k_1=0,\ldots,s_1^{m_1}-1$, $m_2=1,\ldots,n_2$, $k=0,\ldots,s_2^{m_2}-1$.
Similarly, we define the row vector ${\mathbf K}_{\widetilde\zeta}(u,P_\ell)$:
$$
{\mathbf K}_{\widetilde{\zeta}}(u,P_\ell) = 
\row_{m=1,\ldots,n_\ell} \row_{k=0,\ldots,s_\ell^{m}-1}
\frac{K_{\widetilde{\zeta}}^{(0,k)}(u,\tauBa{p^{(m)}_\ell})}{k!}, \quad  p^{(m)}_\ell \in \mathcal P (y_\ell)
.
$$
For the case $p_1^{(m_1)} = p_2^{(m_2)}$, we replace $\frac{K_{\widetilde\zeta}^{(k_1,k_2)}(p_1^{(m_1)},\tauBa{p_2^{(m_2)}})}{k_1!k_2!}$
in ${\mathbf K}_{\widetilde\zeta}(P_1,P_2)$ by the coefficient of $t_{1,m_1}(u)^{k_1}$, where $t_{1,m_1}(u)$ is the corresponding local coordinate, 
in the Laurent series of $\frac{K_{\widetilde\zeta}^{(0,k_2)}(u,\tauBa{p_2^{(m_2)}})}{k_2!}$ at $p_1^{(m_1)}$.
\end{Notation}

\begin{Rk}
Using Notation \ref{notatMatrix1}, the model operator \eqref{MyMatrix} can be written in the following form
\begin{equation}
\label{mymyColl2}
M^{y_k}f(u)
=
y_k(u)f(u) + i \mathbf{K}_{\widetilde{\zeta}}(u,P_k) \, \mathbf{A}_{y_k} \, \Phi_{y_k}(f),
\end{equation}
where $\mathbf{A}_{y_k}$ 
is the block diagonal matrix with 
upper-skew-triangular Hankel blocks given by \eqref{lblAmatrix}.
\end{Rk}

\begin{pf}[of Theorem \ref{lemCollCodI}]
Using Proposition \ref{myGeneral1}, it is sufficient to prove that for 
an arbitrary two-variables polynomial $g(x_1,x_2)$ 
with real coefficients and the meromorphic function $z=g(y_1,y_2)$,
the operator $M^z$ satisfies the appropriate colligation condition.
It is enough to show that if the colligation condition holds for $M^{y_1}$ and $M^{y_2}$
then it holds for $M^{y_1}+M^{y_2}$ and $M^{y_1}M^{y_2}$.
We let $\Phi_k = \Phi_{y_k}$, see \eqref{defPhiYOp}, and $\sigma_k = \sigma_{y_k}$, see \eqref{eqDefSigmaY}.
\smallskip

First, by summing the colligation conditions for $M^{y_1}$ and $M^{y_2}$, we simply have
\begin{equation*}
M^{y_1} + M^{y_2} - M^{{y_1}*} - M^{{y_2}*}
=
i \Phi_1 ^* \sigma_1 \Phi_1 + i \Phi_2 ^* \sigma_2 \Phi_2
=
i \Phi ^* \sigma \Phi,
\end{equation*}
where $\Phi = \Phi_{11} \oplus \widehat{\Phi}_{12} \oplus \Phi_{22}$, see Notation \ref{notatMatrix1},
and $\sigma = \sigma_{11} \oplus \widehat{\sigma}_{12} \oplus \sigma_{22}$.
The block diagonal matrices $\sigma_{11}$ and $\sigma_{22}$ contain the blocks corresponding to the Laurent coefficients of $y_1$ and $y_2$ 
at the poles in $\mathcal P (y_1) \cap \mathcal P (y_2)^{c}$ and $\mathcal P (y_2) \cap \mathcal P (y_1)^{c}$, respectively.
The matrix $\widehat{\sigma}_{12}$ contains the summation of the Laurent coefficients of $y_1$ and $y_2$ at the joint poles.
As for the case where a joint pole does not belong to $\mathcal{P}(y_1 + y_2)$ or its order is less than $\widehat{s}^m$ for $y_1+y_2$, 
then the matrix $\widehat{\sigma}_{12}$ will contain zeros at the corresponding entries.
Thus, the colligation condition for $M^{y_1+y_2}$ follows.
\smallskip

Our next aim is to show that the colligation condition for $M^{y_1y_2}$ holds with $\sigma_{y_1y_2}$ as constructed in \eqref{eqDefSigmaY}
and so we first examine the block entries $\sigma_{y_1y_2}$, more precisely the Hankel matrix ${\mathbf A}_{y_1y_2,r}$ corresponding
to the pole $p^{(r)}$, see \eqref{lblAmatrix}. Its $(\gamma,\delta)$-entry is
\begin{equation}
\label{eqPy1y2}
{\mathbf A}_{y_1y_2,r;\gamma,\delta} = 
\sum_{j=0}^{k+s^r} a_{1,r,-s^r_1+j} a_{2,r,s^r_1+k-j},
\quad \gamma + \delta < s^r,
\end{equation}
where $k=-(\gamma + \delta + 1)$ and $a_{1,r,\cdot}$ and $a_{2,r,\cdot}$ are the Laurent coefficients of $y_1$ and $y_2$ respectively at $p^{(r)}$.
Then one may show that \eqref{eqPy1y2} becomes
\begin{align}
\label{eqSigma12}
{\mathbf A}_{y_1y_2,r}
=
\begin{pmatrix}
{\mathbf A}_{y_1,r}Y^{r}_2 \\ \hdashline[2pt/2pt]
0_{s^{r}_2 \times s^r}
\end{pmatrix}
+
\begin{pmatrix}
{\mathbf A}_{y_2,r} Y^{r}_1 \\ \hdashline[2pt/2pt]
0_{s^{r}_1 \times s^r}
\end{pmatrix}
^T
,
\end{align}
where $Y^r_1$ and $Y^r_2$ are Toeplitz matrices of sizes $s^r_2 \times s^r$ and $s^r_1 \times s^r$, respectively, given by
\begin{equation}
\label{defY2r}
[Y_\ell^r]_{\scriptsize\begin{aligned}\gamma &=1,\ldots,s^r-s_\ell^r  \\[-4pt] \delta & =1,\ldots,s^r \end{aligned}}
=
\begin{cases}
a_{\ell,r,{\gamma-\delta}} &  s^r_\ell \geq \delta - \gamma\\
0 & s^r_\ell < \delta - \gamma
\end{cases}
\end{equation}
(the entries of the first and the second matrix products on the right hand side of \eqref{eqSigma12}
correspond to the index $j \leq s^r_1 - \gamma - 1$ and to $j \geq s^r_1 - \gamma$ in \eqref{eqPy1y2}, respectively).

Moving on to examine $M^{y_1y_2}=M^{y_1}M^{y_2}$, one can show, using the commutativity of $M^{y_1}$ and $M^{y_2}$ 
and the colligation conditions associated to $y_1y_2$ (Proposition \ref{myGeneral1}), the following
\begin{align}
\nonumber
\frac{1}{i}
&
\big(
\innerProductReg{M^{y_1}M^{y_2}f}{g} - \innerProductReg{f}{M^{y_1}M^{y_2}g}
\big)
\\
\nonumber
&
=
\innerProductReg{\sigma_1 \Phi_1 M^{y_2} f }{\Phi_1 g} +
\frac{1}{i}\big(\innerProductReg{M^{y_2}f}{M^{y_1}g}- \innerProductReg{f}{M^{y_1}M^{y_2}g}\big)
\\
\nonumber
&
=
\innerProductReg{\sigma_1 \Phi_1 M^{y_2} f }{\Phi_1 g} +
\frac{1}{i}\big(\innerProductReg{M^{y_2}f}{M^{y_1}g}- \innerProductReg{f}{M^{y_2}M^{y_1}g}\big)
\\
\label{mymyColl}
&
=
\innerProductReg{\sigma_1 \Phi_1 M^{y_2} f }{\Phi_1 g} +
\innerProductReg{\sigma_2 \Phi_2 f }{\Phi_2 M^{y_1} g}
.
\end{align}
For the sake of simplicity, we first assume that $\mathcal P (y_1)\cap \mathcal P (y_2)=\emptyset$.
It follows that the expression in \eqref{mymyColl}, using the matrix representations as presented in \eqref{mymyColl2}, becomes
\begin{align}
\nonumber
\frac{1}{i}
\big(
\innerProductReg{M^{y_1}M^{y_2}f}{g} -
&
\innerProductReg{f}{M^{y_1}M^{y_2}g}
\big)
= \\ =  &
\nonumber
\Phi_{1}(g)^*
\sigma_1
\left(
\Phi_{1}(y_2 \cdot f) 
+
i
\mathbf{K}_{\widetilde{\zeta}}(P_1,P_2)
\mathbf{A}_{y_2}
\Phi_{2}(f)
\right)
+
\\
&
\left(
\Phi_{2}(y_1 \cdot g)^*
-
i 
\Phi_{1}(g)^*
\mathbf{A}_{y_1}^*
\mathbf{K}_{\widetilde{\zeta}}(P_2,P_1)^*
\right)
\sigma_2
\Phi_{2}(f)
\label{eqCollMultE}
,
\end{align}
where ${\mathbf K}_{\widetilde\zeta}(P_2,P_1)$ is defined similarly
to ${\mathbf K}_{\widetilde\zeta}(P_1,P_2)$.

We now show that two of the terms in \eqref{eqCollMultE} vanish.
Using the hermitian structure of the Cauchy kernels,
we have ${\mathbf K}_{\widetilde\zeta}(P_2,P_1)^*= {\mathbf K}_{\widetilde\zeta}(\tauBa{P_1},\tauBa{P_2})$,
where the last matrix is again similarly defined.
Furthermore, using the signed permutation matrices corresponding to real and complex conjugate poles, see \eqref{eqPerY}, we have:
\begin{align}
\nonumber
&
\sigma_1\mathbf{K}_{\widetilde{\zeta}}(P_1,P_2)\mathbf{A}_{y_2}
-
\mathbf{A}_{y_1}^*\mathbf{K}_{\widetilde{\zeta}}(\tauBa{P_1},\tauBa{P_2})\sigma_2
= \\ 
&
\sigma_1\mathbf{K}_{\widetilde{\zeta}}(P_1,P_2)\mathbf{A}_{y_2}
-
{\mathbf A}_{y_1}^*\Per{y_1} {\mathbf K}_{\widetilde\zeta}(P_1,P_2) \Per{y_2}  \sigma_2
=0,
\label{eqK0}
\end{align}
since $\sigma_k = - \Per{y_k} {\mathbf A}_{y_k} = - {\mathbf A}_{y_k}^* \Per{y_k}$,
see Proposition \ref{myGeneral1}.

For the remaining two terms in \eqref{eqCollMultE}, since we assumed that ${\mathcal P}(y_1) \cap {\mathcal P}(y_2) = \emptyset$, it follows
from \eqref{eqSigma12} that for all $p^{(r)} \in {\mathcal P}(y_1)$, ${\mathbf A}_{y_1y_2,r} = {\mathbf A}_{y_1,r} Y_2^r$ 
and the matrix $Y^r_2$ is a square lower-triangular Toeplitz matrix consisting of $a_{1,r,0},\ldots,a_{1,r,s^r_1-1}$.
The $k$-th Taylor coefficient ($0 \leq k < s^{r}=s^{r}_1$) of $y_2 f$ at $p^{(r)}$ is equal to
$\sum_{j=0}^{k} a_{2,r,k-j} \frac{f^{(j)}(p{(r)})}{j!}$. Then we may write
$$
\Phi_1(y_2 f) = \DiagTwo{p^{(r)} \in \mathcal{P}(y_1)}{Y_2^r} \Phi_1(f) 
,$$
and similarly for $\Phi_2(y_1f)$, and \eqref{eqCollMultE} becomes 
\begin{align*}
\frac{1}{i}
\big( &
\innerProductReg{M^{y_1}M^{y_2}f}{g} -
\innerProductReg{f}{M^{y_1}M^{y_2}g}
\big)
\\ = & 
\Phi_{1}(g)^*
\sigma_1
\diag_{p^{(r)} \in \mathcal{P}(y_1)} Y_2^r 
\Phi_1(f)
+
\Phi_{2}(g)^* 
(\diag_{p^{(r)} \in \mathcal{P}(y_2)} Y_1^r)^*
\sigma_2
\Phi_{2}(f)	
\\ = & 
(\Phi_{1} \oplus \Phi_{2})(g)^*
\begin{psmallmatrix} 
\sigma_1 \DiagTwo{p^{(r)} \in \mathcal{P}(y_1)}{Y_2^r}    & 0  \\
0  &     \sigma_2 \DiagTwo{p^{(r)} \in \mathcal{P}(y_2)}{ Y_1^r}
\end{psmallmatrix} 
(\Phi_{1} \oplus \Phi_{2})(f) 
\\ = & 
\Phi_{y_1y_2}(g)^*
\sigma_{y_1y_2}
\Phi_{y_1y_2}(f) 
.
\end{align*}
Notice that when passing from the second line to the third line, we have used the block diagonal structure of the matrices $\sigma$,
with the selfadjoint diagonal blocks at the real poles and the selfadjoint $2 \times 2$  off diagonal block matrices
at the complex conjugate poles, see \eqref{eqDefSigmaY}.
Notice also that in the case where $y_1$ has a pole at $p^{(r)}$ and $y_2$ has a zero at $p^{(r)}$ (or vice versa), 
so that the order of the pole of $y_1y_2$ at $p^{(r)}$ is less than $s^r$,
the main skew diagonal (and possibly some higher skew diagonals) of the Hankel matrix
${\mathbf A}_{y_1,r} Y_2^r$ will be zero since in the Toeplitz matrix $Y_2^r$
the main diagonal (and possibly some lower diagonals) will be zero.
It follows that we can omit these zero entries and the corresponding entries $\frac{f^{(s^r)}(p^{(r)})}{s^r!}$, $\frac{g^{(s^r)}(p^{(r)})}{s^r!}$  
(and possibly some lower derivatives) in $(\Phi_1 \oplus \Phi_2)(f)$, $(\Phi_1 \oplus \Phi_2)g)$,
so that the equality between the third and the fourth lines in the above calculation still holds.
\smallskip

For the case of a pair of meromorphic functions with common poles, it follows from the block diagonal structure of the matrices $\sigma$ that 
it is enough to verify the colligation condition for a single pole (or a pair of complex conjugate poles) and that the previous argument still holds
for $p^{(r)} \in {\mathcal P}(y_1) \cap {\mathcal P}(y_2)^c$ and for $p^{(r)} \in {\mathcal P}(y_2) \cap {\mathcal P}(y_1)^c$.
In order to examine $\Phi_1^{p^{(r)}} M^{y_2} f$ for $p^{(r)}\in \mathcal{P}(y_1)\cap\mathcal{P}(y_1)$,
we first evaluate the $k$-th Taylor coefficient of $M^{y_2}f(u)$ ($k < s_1^{r}$) at $p^{(r)}$.
We consider the expansions of $y_2(u)$, $f(u)$, and $K_{\widetilde{\zeta}}^{(0,l-1)}(u,\tauBa{p^{(r)}})$
using the local coordinate $t_r(u)$ centered at $p^{(r)}$ and then
\begin{align*}
&
\frac{d^k }{du ^k} 
\frac{1}{k!}
\left(
y_2(t_r(u))f(t_r(u)) + i \sum_{l=1}^{s^{r}_2} \sum_{j=l}^{s^{r}_2}
\frac{a_{2,r,-j} f^{(j-l)}(t_r(u))}{ (j-l)!(l-1)!} K_{\widetilde{\zeta}}^{(0,l-1)}(t_r(u),\tauBa{t_r(u)})
\right)
\\ = &
\sum_{j=0}^{s^r_2+k} a_{2,r,-(j-k)}  \frac{f^{(j)}(p^{(r)})}{j!}
+ i \sum_{l=1}^{s^{r}_2} \sum_{j=l}^{s^{r}_2} \frac{a_{2,r,-j} f^{(j-l)}(t_r(u))}{ (j-l)!(l-1)!}  \left(K_{\widetilde\zeta}^{(0,l-1)}(u,\tauBa{p^{(r)}})\right)_k 
\\ = &
[ a_{2,r,k} ,\ldots , a_{2,r,-s^r_2} ,  \underbrace{ 0 ,\dots,0 }_{s^r_1 -1 -k \text{ times}}  ] \Phi^{p^{(r)}}_{12}(f)
+
\\ &
i \row_{l=1,\ldots,s_2^r} \frac{1}{(l-1)!} \left(K_{\widetilde\zeta}^{(0,l-1)}(u,\tauBa{p^{(r)}})\right)_k {\mathbf A}_{y_2,r} \Phi_2^{p^{(r)}}(f)
,
\end{align*}
where $\left(K_{\widetilde{\zeta}}^{(0,l-1)}(u,\tauBa{p^{(r)}})\right)_k$ denotes the $k$-th Laurent coefficient and $\Phi_{12}^{p^{(r)}}(f)$
consists of the Taylor coefficients of $f$ up to order $s^{r}-1$ since $k$ varies between $0$ and $s_1^r-1$. 
The first row matrix appearing on the last right hand side is the row number $k$ of the Toeplitz matrix $Y_2^r$
(here $Y_2^r$ is not a square matrix since $s^{r_1}\neq 0$, see \eqref{defY2r})
whereas the second row matrix is the same row of the block $K_{\widetilde\zeta}(P_1^{(r)},P_2^{(r)})$
in $K_{\widetilde\zeta}(P_1,P_2)$ corresponding to $p_1=p_2=p^{(r)}$, see Notation \ref{notatMatrix1}, and then we have:
\begin{align*}
\Phi^{p^{(r)}}_{1} M^{y_2} f
= &
\col_{k=0,\ldots,s^r_1-1}\left( \frac{d^k }{d t(u) ^k} \frac{1}{k!} M^{y_2} f(u) \restrict{t(u)=0}\right)
\\ = &
Y_2^r \Phi^{p^{(r)}}_{12}(f) + i {\mathbf K}_{\widetilde\zeta}(P_1^{(r)},P_2^{(r)}) {\mathbf A}_{y_2,r} \Phi_2^{p^{(r)}}(f).
\end{align*}
Similarly, we have
$$
\Phi_2^{p^{(r)}}(M^{y_1}g) =
Y_1^r \Phi_{12}^{p^{(r)}}(g) + i {\mathbf K}_{\widetilde\zeta}(P_2^{(r)},P_1^{(r)}) 
{\mathbf A}_{y_1,r} \Phi_1^{p^{(r)}}(f).
$$
If $p^{(r)}$ is a real pole, we substitute the last two equalities on the right hand side of the last line of \eqref{mymyColl},
repeat the same argument as in \eqref{eqCollMultE}--\eqref{eqK0} to see that two of the terms vanish,
and for the remaining two terms we have up multiplication by $-(-1)^{\mu_{i_r}+\nu_{i_r}}$
(we expand naturally $\Phi^{p^{(r)}}_{1}$ and $\Phi^{p^{(r)}}_{2}$ to $\Phi^{p^{(r)}}_{12}$)
\begin{align*}
&
\Phi^{p^{(r)}}_{1}(g)^*
{\mathbf A}_{y_1,r}
Y_2^r \Phi^{p^{(r)}}_{12}(f)
+
\Phi^{p^{(r)}}_{12}(g)^*
(Y_1^r)^*
{\mathbf A}_{y_2,r}^{*}
\Phi^{p^{(r)}}_{2}(f)
\\ = &
\Phi^{p^{(r)}}_{12}(g)^*
\left(
\begin{pmatrix}
\mathbf{P}_1^{r} Y_2^r \\
0_{s^r_2\times s^r}
\end{pmatrix}
+
\begin{pmatrix}
(Y_1^r)^T \mathbf{P}_2^{r}  & 
0_{s^r\times s^r_1}
\end{pmatrix}
\right)
\Phi^{p^{(r)}}_{12}(f)
.
\end{align*}
The sum in parentheses is exactly the matrix $A_{y_1y_2,r}$, see \eqref{eqSigma12} above,
so that we obtain the diagonal block corresponding to $p^{(r)}$ in $\sigma_{y_1y_2}$ of the colligation condition for $M^{y_1y_2}$.
The calculation for the case of a pair complex conjugate poles is similar.
\end{pf}

To conclude this section, we present the proof of Theorem \ref{vesselAreEq}, but first we present a short technical lemma.

\begin{lem}
\label{eqHankelInv}
Let $G$ be an $s\times s$ upper-skew-triangular Hankel matrix (a Hankel matrix with zero entries below the main skew-diagonal)
$$G_{i,\ell} = g_{i+\ell+1}, \quad
0 \leq \ell, i <s, \quad
i + \ell < s,$$ 
where $g_{s}\neq 0$.
Then, the inverse of $G$ 
is a lower-skew-triangular Hankel matrix which is given by
$$ 
\left(G^{-1}\right)_{i,\ell} 
=  
\restr{ \frac{d^{i+\ell+1}}{dx^{i+\ell+1}} 
\frac{1}{(i+\ell+1)! g(x) } 
}{x = 0} 
\quad  
0 \leq \ell, i <s, 
\quad 
\ell+ i \geq s -1,
$$
where $g(x) = \sum_{i=1}^{s} g_{i} x^{-i}+g_1(x)$ for any $g_1(x)$ analytic at zero.
\end{lem}
\begin{pf}
Let us denote the lower-skew-triangular Hankel matrix given in the lemma by $F$.
Then for $\alpha \leq \beta$ (otherwise, the entry is zero) we have
\begin{align}
\nonumber
\left(G F \right)_{\alpha,\beta}
= &
\sum_{i = s-\beta-1}^{s-\alpha-1} 
G_{\alpha,i} F_{i,\beta}
=
\sum_{i = 0}^{\beta - \alpha} 
G_{\alpha,s-\beta+i-1} F_{s-\beta+i-1,\beta}
\\ = &
\nonumber
\sum_{i = 0}^{\beta - \alpha} 
g_{s-(\beta-\alpha -i)}
\restr{ \frac{d^{s+i}}{dx^{s+i}} \frac{1}{(s+i)! g(x)}}{x=0},
\end{align}
which is the $(\beta-\alpha)$-th Laurent coefficient of 
the multiplication of $g(x)$ by $\frac{1}{g(x)}$ and hence equal to $\delta_{\beta,\alpha}$.
\end{pf}
\begin{pf}[of Theorem \ref{vesselAreEq}]
To prove the validity of the right most equality in \eqref{eqVesselEq}, 
we first define the unitary mapping $U$ by \eqref{ModelMap1}. 
Therefore, we have
\[
U \, \colon \, h 
\mapsto  
\varphi_h(p) = {\bf \widetilde{u}}_l^{\times} (p) (\xi \sigma)\Phi (\xi A - \xi y(p) )^{-1}h,
\]
where $\xi \sigma$, $\xi A$ and $\xi y(p)$ stand for $\xi_1 \sigma_1+\xi_2\sigma_2$, $\xi_1A_1+\xi_2A_2$ and $\xi_1y_1(p) +\xi_2y_2(p)$, respectively.
We fix an arbitrary element $h\in H$, 
we consider a local coordinate $t_j(u)$ centered at $p^{(j)} \in \mathcal P (y_1 ) \cup \mathcal P (y_2 )$ 
and we fix $(\xi_1,\xi_2)$ such that $\xi y(p)$ has a pole of order $\hat{s}^{j}$ at $p^{(j)}$ for all $j$ 
(see Notation \ref{notatMatrix1}; notice that this is equivalent to $\xi\sigma$ being invertible).
Then the operator valued function $(\xi A - \xi y(p) )^{-1} $ has a zero of order $\hat{s}^{j}$ at $p^{(j)}$
with the $\hat{s}^j,\ldots, 2 \hat{s}^j - 1$ Taylor coefficients equal to those of $\frac{-1}{\xi y(p)}I_{H}$.
We recall that ${\bf \widetilde{u}}_l^{\times} (p)$ has a pole at $(p^{(j)})$ of order $\hat{s}^{j}$.
It follows that $\varphi_h(p)$ is analytic at the poles of $\xi y$.
Let $\Phi_{\xi y}\varphi_h(p)$ be the vector of the Taylor coefficients of $\varphi_h(p)$ at the poles of $\xi y(p)$ as before \eqref{defPhiYOp}.
Then, applying Lemma \ref{eqHankelInv} to $\xi y$
at each of the poles $p^{(j)}$, we have that
\begin{equation}
\label{eq123asd}
\Phi_{\xi y} \varphi_h(p)
= 
-\mathbf{A}_{\xi y}^{-1} \widetilde{U}_l^\times (\xi \sigma) \Phi h
.
\end{equation}
Here ${\mathbf A}_{\xi y}$ is the block diagonal matrix \eqref{lblAmatrix} with upper-skew-triangular Hankel blocks consisting of the coefficients
of the singular parts of $\xi y(p)$ and $\widetilde{U}^\times_l$ denotes the $n \times n$ matrix with rows the coefficients 
of the principal parts of $\widetilde{\mathbf u}^\times_l(p)$, that is,
\begin{equation}\label{eq123asd12}
\widetilde{U}_l^\times 
=  
\col_{j=1,\ldots,\abs{\mathcal{P} (y_1) \cup \mathcal{P}(y_2)}} \col_{k=1,\ldots,\hat{s}_j}({\bf \widetilde{u}}_l^{\times})_{(j,-k)}.
\end{equation}
where $(\cdot)_{(j,-k)}$ denotes the $-k$ Laurent coefficient at $p^{(j)}$.
It follows from \eqref{eq123asd} that $\Phi h = \Phi_{\mathrm{Mod}} (\varphi_h)$, where 
\begin{equation}
\label{eq123asdasdas}
\Phi_{\mathrm{Mod}} = -(\xi\sigma)^{-1} (\widetilde{U}_l^\times)^{-1} \mathbf{A}_{\xi y} \Phi_{\xi y}.
\end{equation}
In particular, when $y_1$ and $y_2$ have simple poles and one considers
the canonical determinantal representation constructed explicitly in \cite{MR1704479}, 
see Steps \ref{step1} \& \ref{step2} of the proof of Theorem \ref{preTh} in Section \ref{secProofStruct},
we have $\widetilde U^\times_l = \Per{y_1,y_2}$ (see \eqref{normSectionCK}),
$\xi\sigma = - \Per{y_1,y_2} {\mathbf A}_{\xi y}$ and hence $\Phi_{\text{Mod}} = \Phi_{\xi y}$; 
for the general case, see Remark \ref{rkUs123123} below.
We note that $\Phi_{\mathrm{Mod}} $ in \eqref{eq123asdasdas} does not depend on 
the selection of $(\xi_1,\xi_2)$ since \eqref{ModelMap0}
does not depend on $(\xi_1,\xi_2)$, see Lemma \ref{lemmModelSpaceMapProp}.
\smallskip

To verify the validity of the two remaining identities in \eqref{eqVesselEq},
we illustrate and continue the calculation for the operator $A_1$ and the direction $\xi=(1,0)$
(we can always assume by a linear change of variables that $\sigma_1$ is invertible and $s^j_1 = \hat s^j$ for all $j$):
\begin{align}
\nonumber
\varphi_{A_1 h}(p   ) 
= &
{\bf \widetilde{u}}_l^{\times} (p)
\sigma_1
\Phi 
( A_1 - y_1(p) )^{-1} A_1 \, h
\\ = &
{\bf \widetilde{u}}_l^{\times} (p)
\sigma_1 \Phi 
( A_1 - y_1(p) )^{-1} 
\left( A_1 - y_1(p) +y_1(p) \right)  \, h
\nonumber
\\ = &
y_1(p)
\varphi_h(p) 
+
{\bf \widetilde{u}}_l^{\times} (p) \sigma_1 \Phi h.
\label{eq123asd1}
\end{align}
Substituting \eqref{eq123asd} where $\xi=(1,0)$ shows that the second summand on the right hand side of \eqref{eq123asd1} equals 
$- \widetilde{\mathbf u}^\times_l(p) (\widetilde{U}^\times_l)^{-1} {\mathbf A}_{y_1} \Phi_{y_1} \varphi_h$.
However it follows from the definition of $\widetilde{U}^\times_l$ that $\widetilde{\mathbf u}^\times_l(p)$ and 
$-i {\mathbf K}_{\widetilde\zeta}(p,P_1) \widetilde U^\times_l$ have the same poles and the same principal parts.
Since these are row vectors with entries sections of $L_{\widetilde\zeta} \otimes \Delta$
and since $L_{\widetilde\zeta} \otimes \Delta$ has no global holomorphic sections it follows that
\begin{equation} \label{eQWeNeedItCfBelow}
\widetilde{\mathbf u}^\times_l(p) = -i {\mathbf K}_{\widetilde\zeta}(p,P_1) \widetilde{U}^\times_l.
\end{equation}
Therefore the second summand on the right hand side of \eqref{eq123asd1} equals 
$i {\mathbf K}_{\widetilde{\zeta}}(p,P_1) {\mathbf A}_{y_1} \Phi_{y_1} \varphi_h$, and this leads to
\begin{equation*}
\varphi_{A_1 h}(p) = y_1(p) \varphi_h(p) + i \mathbf{K}_{\widetilde{\zeta}}(p,P_1) \mathbf{A}_{y_1} \Phi_{y_1}  \varphi_h.
\end{equation*}
This is, see \eqref{mymyColl2}, the definition of the model operator $M^{y_1} \varphi_h$ corresponding to $y_1$. 
\smallskip

Finally, it remains to show that the Hilbert space $\mathcal{H}(T) = \{ \varphi_h \, : \, h \in H \}$
with the inner product inherited form $H$, 
is a reproducing kernel Hilbert space with the reproducing kernel \eqref{eqKernPos}.

We use the properties of the characteristic functions of the vessel $\mathcal{V}$. 
Since $\widetilde{\mathbf u}^\times_l(p) = \widetilde{\mathbf u}^\times(\tauBa{p})^*$, 
\eqref{st} and \eqref{eqJcf}, together with $T(p)T(\tauBa{p})^*=1$ and \eqref{eqWinA}, imply that
$$ T(p) {\mathbf u}^\times_l(p) = \widetilde{\mathbf u}^\times_l(p) \xi\sigma W(\xi_1,\xi_2,\xi y(p)) (\xi\sigma)^{-1}.$$
Together with an identity which is presented in the next section, see \eqref{eqsec5A} below,
and \eqref{eqMatNormSecRelation} applied to both the input and the output determinantal representations, this yields the following
\begin{align*}
K_T(p,q) = & T(p) K_{\zeta}(p,q) T(q)^* - K_{\widetilde{\zeta}}(p,q)
\\ = &
\frac
{ T(p) {\bf u}^{\times} _l (p) ( \xi \sigma ) {\bf u}^{\times} _l (q)^* T(q)^* }
{-i  (\xi y(p) - \xi \overline{y(q)}) }
-
\frac{\widetilde{\bf u}^{\times} _l (p)(\xi\sigma){\widetilde{\bf u}}^{\times} _l (q)^* } 
{-i  (\xi y(p) - \xi \overline{y(q)}) }
\\ = &
\widetilde{\bf u}^{\times} _l (p)
\frac{\xi\sigma W(\xi_1,\xi_2,\xi y(p)) (\xi\sigma)^{-1} W(\xi_1,\xi_2,\xi y(q))^* \xi\sigma - \xi\sigma}{-i  (\xi y(p) - \xi \overline{y(q)}) }
\widetilde{\bf u}^{\times} _l (q)^*
\\ = &
\widetilde{\bf u}^{\times} _l (p)
(\xi\sigma)
\Phi
(\xi A - \xi y(p))^{-1} 
(\xi A - \xi y(q))^{-*}
\Phi^*
(\xi \sigma)^*
\widetilde{\bf u}^{\times} _l (q)^*
.
\end{align*}
Thus, we may conclude that
\begin{align*}
\innerProductTri{\varphi_h(p)}{ K_T(p,q)}{\mathcal{H}(T)}
= &
\innerProductTri{h}{
(\xi A - \xi y(q))^{-*}
\Phi^*
(\xi \sigma)^*
\widetilde{\bf u}^{\times} _l (q)^*
}{H}
\\ = &  
\widetilde{\bf u}^{\times} _l (q)
(\xi \sigma)
\Phi
(\xi A - \xi y(q))^{-1}
h
\\ = &
\varphi_h(q).
\end{align*}
\end{pf}

\begin{Rk}
\label{rkUs123123}
Let us define analogously to \eqref{eq123asd12}
$$\widetilde U^\times = \row_{j=1,\ldots,\abs{\mathcal{P} (y_1) \cup \mathcal{P}(y_2)}} \row_{k=1,\ldots,\hat{s}_j} (\widetilde{\mathbf u}^\times)_{(j,-k)},$$
then it follows as in \eqref{eQWeNeedItCfBelow} that
$$ \widetilde{\mathbf u}^\times(p) = i U^\times {\mathbf K}_{\widetilde\zeta}(P_1,p),$$
where the column vector ${\mathbf K}_{\widetilde\zeta}(P_1,p)$  is defined similarly to Notation \eqref{notatMatrix1}.
Plugging \eqref{eQWeNeedItCfBelow} and the last equality into the first identity of \eqref{eqMatNormSecRelation},
and comparing to \eqref{eqCollFor}, which can be rewritten as
${\mathbf K}_{\widetilde\zeta}(p,P_1) {\mathbf A}_{\xi y} {\mathbf K}_{\widetilde\zeta}(P_1,\tauBa{q}) = i(\xi y(p) - \xi \overline{y(q)}) K_{\widetilde\zeta}(p,q)$,
yields then $-\widetilde{U}^\times_l \xi\sigma \widetilde{U}^\times = {\mathbf A}_{\xi y}$.
It follows that $(\xi\sigma)^{-1} (\widetilde{U}^\times_l)^{-1} {\mathbf A}_{\xi y} = -\widetilde{U}^\times$ 
and then finally it follows from \eqref{eq123asdasdas} that 
$\Phi_{\mathrm{Mod}} = (\widetilde{U}^\times) \Phi_{\xi y} $.
\end{Rk}

\begin{Rk}
\label{rkUs2}
Using the notations introduced above, my may generalized Remark \ref{rk26D} for the general case of not necessarily simple poles.
In general, the joint characteristic function $S$ satisfies the property that $S(p)-I$ has zero at $p^{(j)}$ of order $\hat{s}^{j}$
(this is referred to as equal to identity at infinity, see \cite{KLMV} and the remark in \cite[p. 280]{MR97m:30051}).
Then, we consider a block diagonal matrix where the $j$-th block is a lower triangular Toeplitz matrix
with the $(\alpha,\beta)$-entry equal to the Taylor coefficient $T_{j,\alpha-\beta}$ of $T$ at $p^{(j)}$ ($1 \leq \beta \leq \alpha \leq \hat s^j$).
Denoting this matrix by $T_{\infty}$ and analyzing the power series expansion of \eqref{st} at each of the poles $p^{(j)}$, 
it follows that the equalities in \eqref{eq26A} and \eqref{eq26B} still hold (see \cite[Proposition 11.2.2]{KLMV}, for a particular case).
\end{Rk}

\section{Proof of the structure theorem}
\label{secProofStruct}

As mentioned in Section \ref{secMainThm}, our strategy to prove Theorem \ref{preTh} is to embed  the operators $M^{y_1}$ and $M^{y_2}$
in a commutative two-operator vessel of the form
\begin{equation}
\label{vessel123}
( \, M^{y_1} \, ,  \,  M^{y_2} \, ; \, \mathcal X \, , \,  \Phi \, , \, E \, ; \,  \sigma_1 \, , \,  \sigma_2 \, , \, \gamma  \, , \, \widetilde{\gamma} \, ),
\end{equation}
where (see Notation \ref{notatMatrix1}) $n = \Card{\mathcal{P} (y_1) \cup \mathcal{P}(y_2)}$, $E=\mathbb{C}^n$, and where
$\Phi$ is the evaluation operator from $\mathcal{X}$ to $E$ at the poles of $y_1$ and $y_2$, namely,
\begin{equation}
\label{phi123}
f \mapsto 
 \begin{pmatrix}
  f(p^{(1)}) \\
  \vdots \\
  f(p^{(n)})
 \end{pmatrix}.
\end{equation}
The discriminant curve of the vessel \eqref{vessel123} will turn out to be 
the image $C$ of the birational embedding $p \mapsto (y_1(p),y_2(p))$ of $X$ into ${\mathbb P}^2$
(see Notation \ref{not4A} \ref{not4A_3}).
\smallskip

\begin{pf}[of Theorem \ref{preTh}]
The proof consists of seven steps.
We start by recalling the relation between the colligation condition and the structure identity.
	
\begin{Step}
\label{step1}
The colligation conditions \eqref{collCond} of the collection \eqref{vessel123} are equivalent to 
the structure identities (condition \ref{assump3}, Theorem \ref{preTh}), 
where $\sigma_k$ for $k=1,2$ are given in \eqref{sigma12} below.
\end{Step}

This result is just a special case of Proposition \ref{myGeneral1} which yields
the colligation condition \eqref{collCond} for the collection \eqref{vessel123}
\begin{equation}
\label{eqCollCondAA}
\frac{1}{i} 
\innerProductTri{\left( M^{y_k}-(M^{y_k})^* \right)F}{G}{\mathcal X}
= 
\innerProductTri{\sigma_k \Phi (F)}{\Phi (G)}{E}
\qquad
k=1,2.
\end{equation}
Here $\Phi$ is the evaluation operator \eqref{phi123} at the $n$ poles of $y_1$ and $y_2$ and $\sigma_k$ is given by
\begin{equation} \label{sigma12}
\sigma_k = \Per{y_1,y_2} \diag_{j=1,\ldots,n}(c_k^j),
\end{equation}
where $\text{Per}(y_1,y_2)$ is a matrix of size $n \times n$ 
constructed as in \eqref{eqPerY} but for the union of poles of $y_1$ and $y_2$
and where $c^j_k$ is minus the residue of $y_k$ at $p^{(j)}$
and is set to be zero whenever $p^{(j)}$ is not a pole of $y_k$.
The matrices $\sigma_1$ and $\sigma_2$ are selfadjoint, since $y_1$ and $y_2$ are real.

\begin{Step}
\label{step2}
Let $\widetilde{\gamma}$ be defined by $\widetilde{\gamma} = \Per{y_1,y_2} \widetilde{\gamma}^\prime$ where
\begin{equation}
\label{tildeGamma}
\widetilde{\gamma}^{\prime}_{j,l} =
\begin{cases}
c_{1}^{l} h_{2}^{l} - c_{2}^{l} h_{1}^{l} , 
& 
p^{(j)} = p^{(l)}
\\
\bracketsA{ c_{1}^{j} c_{2}^{l}  - c_{2}^{j} c_{1}^{l}  }
\frac
{\vartheta[\widetilde{\zeta}]( p^{(l)}- {p^{(j)}})}
{\vartheta[\widetilde{\zeta}] (0)E(p^{(l)},{p^{(j)}})},
& 
{\text otherwise},
\end{cases}
\end{equation}
where $c_{k}^{j}$ and $h_{k}^{j}$, for a $p^{(j)} \in \mathcal{P}(y_k)$, 
are defined by the expansion of $y_k$ at $p^{(j)}$ 
using a local coordinate $t_j$ centered at $p^{(j)}$
\[
y_k(u) = - \frac{c^{m}_{k}}{t_j(u)}+h^{m}_{k} + O(|t_j|),
\]
otherwise ($p^{(j)} \notin \mathcal{P}(y_k)$), we set $c_{k}^{j}=0$ and $h_{k}^{j} = y_k(p^{(j)})$.
Then, the output vessel condition
\begin{equation}
\label{eq4_3}
\sigma_1\Phi M^{y_2} - \sigma_2\Phi M^{y_1}  = \widetilde{\gamma}\Phi
\end{equation}
holds for $\sigma_1$ and $\sigma_2$ as defined in \eqref{sigma12}.
\end{Step}

The matrix $\widetilde{\gamma}$ is given explicitly by
\begin{equation}\label{eqTildeGamma2}
\widetilde{\gamma}_{j,k} =
\begin{cases}
(-1)^{\mu_{i_j}+\nu_{i_j}}
\left(c_{1}^{k} h_{2}^{k}- c_{2}^{k} h_{1}^{k}\right)
& 
\tauBa{p^{(j)}} = p^{(k)}, \quad j=k,
\\
c_{1}^{k} h_{2}^{k} - c_{2}^{k} h_{1}^{k} & 
\tauBa{p^{(j)}} = p^{(k)}, \quad j \neq k,
\\
\bracketsA{\overline{c_{1}^{j}} c_{2}^{k}  - \overline{c_{2}^{j}} c_{1}^{k} }
\frac
{\vartheta[\widetilde{\zeta}]( p^{(k)}- \tauBa{p^{(j)}})}
{\vartheta[\widetilde{\zeta}] (0)E(p^{(k)},\tauBa{p^{(j)}})}
,
& 
{\text otherwise},
\end{cases}
\end{equation}
where a real pole $p^{(j)}$ belongs to $X_{i_j}$.
$\widetilde{\gamma}$ is selfadjoint since $\zeta \in T_\nu$ and the meromorphic functions $y_1$ and $y_2$ are real.
This completes the proof of Step \ref{step2}.

\begin{pf}[of Step \ref{step2}]
We note that $\mathbf{A}_k$, defined in the proof of Proposition \ref{myGeneral1}, 
is a diagonal matrix (since we are now in a simple poles setting).
Furthermore \eqref{eq4_3} becomes
\begin{equation}
\label{eq4_3A}
\mathbf{A}_2 \Phi M^{y_1} - \mathbf{A}_1 \Phi M^{y_2} = \Per{y_1,y_2} \widetilde{\gamma}^\prime \Phi.
\end{equation}
Below, for the sake of simplicity, 
we begin by calculating the $j$-th entry of the left hand side of \eqref{eq4_3A}
for the case where $p^{(j)} \in \mathcal{P}(y_2) \cap \mathcal{P}(y_1)^c$, we use $h_1^j=y_1(p^{(j)})$ and $c_1^j=0$ to conclude:
\begin{align}
\nonumber
[LHS]_j
= &
\left(c_1^j y_2(p^{(j)}) - c_2^j y_1(p^{(j)})\right) f(p^{(j)}) -
\sum_{ p^{(l)} \in \mathcal{P}(y_1)}{c_1^l f(p^{(l)})
\frac
{\vartheta[\widetilde{\zeta}](p^{(l)}-p^{(j)})}
{\vartheta[\widetilde{\zeta}] (0)E(p^{(l)},p^{(j)})}
}
\\ = &
-
c_2^j
h_1^j 
f(p^{(j)})
-
c_2^j
\sum_{ p^{(l)} \in \mathcal{P}(y_1)}
c_1^l f(p^{(l)})
\frac
{\vartheta[\widetilde{\zeta}](p^{(l)}-p^{(j)})}
{\vartheta[\widetilde{\zeta}] (0)E(p^{(l)},p^{(j)})}
\label{step2_2B}
.
\end{align}
For the case where $p^{(j)} \in \mathcal{P}(y_2) \cap \mathcal{P}(y_1) $,
we use the local coordinate $t_j(u)$ and the properties of the Cauchy kernel to get:
\begin{align*}
[\mathbf{A}_2 \Phi M^{y_1}f]_j
= & -
c_2^j
 f(p^{(j)}) 
\restr{\left(\frac{-c^{j}_{1}}{t_j(u)}+h^{j}_{1} + o(|t_j|) + \frac{c_1^j}{t_j(u)}\right)}{t_j(u)=0} - \\ &
c_2^j \sum_{ p^{(j)} \neq p^{(l)} \in \mathcal{P}(y_1)  } c_1^l f(p^{(l)})
\frac{\vartheta[\widetilde{\zeta}](p^{(l)}-p^{(j)})}{\vartheta[\widetilde{\zeta}] (0)E(p^{(l)},p^{(j)})},
\end{align*}
and similarly for $[\mathbf{A}_1 \Phi M^{y_2} f]_j$.
Therefore, the LHS of \eqref{eq4_3A}, in the case of a joint pole, is given by
\begin{align}
\nonumber
[LHS]_j
= 
\left(
c_1^j h_2^j 
-
c_2^j h_1^j 
\right) f(p^{(j)})
& -
c_2^j
\sum_{ \substack{p^{(l)} \in \mathcal{P}(y_1) \\  p^{(j)} \neq p^{(l)} } } 
c_1^l f(p^{(l)})
\frac
{\vartheta[\widetilde{\zeta}](p^{(l)}-p^{(j)})}
{\vartheta[\widetilde{\zeta}] (0)E(p^{(l)},p^{(j)})}
\\ & +
c_1^j
\sum_{ \substack{p^{(l)} \in \mathcal{P}(y_2) \\  p^{(j)} \neq p^{(l)} } } 
c_2^l f(p^{(l)})
\frac
{\vartheta[\widetilde{\zeta}](p^{(l)}-p^{(j)})}
{\vartheta[\widetilde{\zeta}] (0)E(p^{(l)},p^{(j)})}
\label{step2_1}
.
\end{align}
One may note, using that $c_k^j = 0$ when $p^{(j)} \notin \mathcal{P}(y_k)$, that both \eqref{step2_2B} (and the similar expression for the case where 
$p^{(j)} \in {\mathcal P}(y_1) \cap {\mathcal P}(y_2)^c$) and \eqref{step2_1} can be written as
\begin{align}
\nonumber
[LHS]_j
= & f(p^{(j)})\left(c_1^j h_2^j - c_2^j h_1^j  \right) - \\ &
\sum_{ \substack{p^{(l)} \in \mathcal{P}(y_1) \cup \mathcal{P}(y_2) \\  p^{(j)} \neq p^{(l)} } } 
f(p^{(l)})\left(c_2^j c_1^l - c_1^j c_2^l\right) 
\frac{\vartheta[\widetilde{\zeta}](p^{(l)}-p^{(j)})}{\vartheta[\widetilde{\zeta}] (0)E(p^{(l)},p^{(j)})}.
\label{step2_VV}
\end{align}
The right hand side of \eqref{step2_VV} coincides exactly with $[\widetilde\gamma' \Phi f]_j$, where $\widetilde\gamma'$ is given by \eqref{tildeGamma}.
This completes the proof of Step \ref{step2}.
\end{pf}

Before moving to prove Step \ref{steppp}, we make several remarks.
At this stage, we have constructed, using the Cauchy kernels, the output determinantal representation of the collection \eqref{vessel123}
from the line bundle $L_{\wt{\z}}\otimes \Delta$ and the pair of meromorphic functions that gives the birational
embedding $p \mapsto (y_1(p),y_2(p))$ of $X$ into $\mathbb{P}^2$.
This canonical determinantal representation was introduced in \cite{MR1704479} (see also \cite{shamovich2014livsic}),
up to left multiplication by $\Per{y_1,y_2}$.
It is shown there that it is indeed a determinantal representation of $C$ corresponding to $\widetilde{\zeta}$,
i.e., its kernel bundle is isomorphic up to a twist (see the discussion preceding \eqref{eqIsoExplicit}) 
to $L_{\widetilde\zeta} \otimes \Delta$ (the canonical determinantal representation is always maximal and
fully saturated, the proof is given in \cite[Theorem 5.1]{MR1704479} under the assumption
that the singularities of $C$ are ordinary multiple points, e.g., nodes).
\smallskip

The associated left normalized section at $v$, is given by a vector of Cauchy kernels at $v$ and $\pi^{-1}({\infty})$,
namely, it is a vector of form 
\begin{equation}
\label{normSectionCK}
[\widetilde{\bf u}^{\times}_l(v)]_j = 
\frac{\vartheta[\widetilde{\zeta}](\tauBa{p^{(j)}}-v)}{\vartheta[\widetilde{\zeta}](0)E(\tauBa{p^{(j)}},v)}
\end{equation}
(see \cite{MR1704479} and note that here $p^{(j)}$ is evaluated under $\tauBa{\cdot}$ due to the multiplication by $\Per{y_1,y_2}$).
\smallskip

The next result follows from the Generalized Cayley--Hamilton Theorem (Theorem \ref{polyVanish})
since Step \ref{step4} below implies that the vessel ${\mathcal V}$ is irreducible. 
We use a different approach based on Theorem \ref{rkModelOpAlg}.
\begin{lem}
The joint spectrum of $M^{y_1}$ and $M^{y_2}$ lies in $C_0$.
\end{lem}
\begin{pf}
Let us assume that $(\lambda_1,\lambda_2)\in \Spec{M^{y_1},M^{y_2}}$.
Then, using the spectral mapping theorem, we have 
\[
p(\lambda_1,\lambda_2) \in \Spec{p(M^{y_1},M^{y_2})},
\]
where $p$ is the discriminant polynomial of the vessel. 
Furthermore, using Theorem \ref{rkModelOpAlg}, we conclude that
\[
p(M^{y_1},M^{y_2}) = M^{p(y_1,y_2)}  = 0.
\]
Thus, $p(\lambda_1,\lambda_2)=0$ and $(\lambda_1,\lambda_2)\in C_0$ follows.
\end{pf}

\begin{Step}
\label{steppp}
We define $\gamma$ by
\begin{equation}
\label{input}
\gamma_{j,k} = \Psi_{j,k}+
\begin{cases}
(-1)^{\mu_{i_j}+\nu_{i_j}} \left(c_{1}^{k} h_{2}^{k}- c_{2}^{k} h_{1}^{k}\right)
& 
\tauBa{p^{(j)}} = p^{(k)}, \quad j=k,
\\
c_{1}^{k} h_{2}^{k} - c_{2}^{k} h_{1}^{k}
& 
\tauBa{p^{(j)}} = p^{(k)}, \quad j \neq k,
\\
\bracketsA{\overline{c_{1}^{j}} c_{2}^{k}  - \overline{c_{2}^{j}} c_{1}^{k}}
\frac{\vartheta[\widetilde{\zeta}]( p^{(k)}- \tauBa{p^{(j)}})}{\vartheta[\widetilde{\zeta}] (0)E(p^{(k)},\tauBa{p^{(j)}})},
& 
{\text otherwise},
\end{cases}
\end{equation}
where
\begin{align*}
\Psi = & \Per{y_1,y_2} 
\left[i(c^i_1 \overline{c^j_2} - c^i_2 \overline{c^j_1}) K_{{\mathcal X}}(p^{(i)},p^{(j)})\right]_{i,j=1,\ldots,n} \Per{y_1,y_2}  
\\ = &
\left[ i(\overline{c^i_1} c^j_2 - \overline{c^i_2} c^j_1) K_{{\mathcal X}}(\tau(p^{(i)}),\tau(p^{(j)})) \right]_{i,j=1,\ldots,n}.
\end{align*}
Then, the input vessel condition \eqref{VesselCond1} and the linkage condition \eqref{VesselCond3} hold.
\end{Step}
Since $\mathcal{X}$ is a reproducing kernel space, it follows that
\begin{equation}
\label{kernelEval}
\Phi\Phi^* = \left[ K _ {\mathcal{X}}(p^{(i)},p^{(j)}) \right]_{i,j=1,\ldots,n}.
\end{equation}
Then, Equation \eqref{input} is derived by substituting
\eqref{sigma12}, \eqref{eqTildeGamma2} and \eqref{kernelEval} in the linkage condition \eqref{VesselCond3}.
\begin{Step}
\label{step4}
The mapping
\begin{equation}
\label{ModelMap}
h 
\mapsto
\widetilde{\bf u}^{\times}_l(v) 
(\xi_1 \sigma_1 + \xi_2 \sigma_2)
\Phi(\xi_1 A_1 + \xi_2 A_2 - \xi_1 y_1(v) -\xi_2 y_2(v))^{-1}h,
\end{equation}
where $h\in\mathcal{X}$,
is the identity,
in the sense that the two sections of $L_{\wt{\z}}\otimes \Delta$
coincide in $\Omega \cap \left(X \setminus \pi^{-1}(\Spec{M^{y_1},M^{y_2}})\right)$.
\end{Step}

We prove that the section on the right hand side coincide with $h$ in the neighborhood of the poles of $y_1$ and $y_2$.
The statement then follows by analytic continuation.
\smallskip

Let us recall that the mapping \eqref{ModelMap} is independent of the choice of the direction $(\xi_1,\xi_2)$ 
(see \cite[Section 3]{MR1634421} and Proposition \ref{lemmModelSpaceMapProp} above for more details).
Hence, for the sake of simplicity, we illustrate the calculation in the direction $\xi=(1,0)$.
Furthermore, we take $v\in \Omega$ such that $y_1(v)$ does not belong to the spectrum of $M^{y_1}$
and $y_1^{-1}(y_1(v))\subset \Omega$. 
A direct computation leads to:
\begin{align}
\nonumber
h \mapsto &
\widetilde{\bf u}^{\times}_l(v) (1 \cdot \sigma_1 + 0 \cdot \sigma_2)\Phi(M^{y_1} - y_1(v) )^{-1}h
\\= & \nonumber
\widetilde{\bf u}^{\times}_l(v) \sigma_1 \Phi R^{y_1}_{y_1(v)}(h)
\\\nonumber= &
\widetilde{\bf u}^{\times}_l(v) \sigma_1 \Phi
\left(
\frac{h(\cdot)}{y_1(\cdot) - y_1(v)} - \sum_{j=1}^{n} \frac{h(u^{(j)})}{dy_1(u^{(j)})}
\frac{\vartheta[\widetilde{\zeta}](u^{(j)}-\cdot)}{\vartheta[\widetilde{\zeta}] (0)E(u^{(j)},\cdot)}
\right)
\end{align}
and so
\begin{align}
\nonumber
h 
\mapsto &
\widetilde{\bf u}^{\times}_l(v) \sigma_1
\left[
\restr{
\frac{h(u)}{y_1(u) - y_1(u)} -
\sum_{j=1}^{n} \frac{h(u^{(j)})}{dy_1(u^{(j)})}
\frac{\vartheta[\widetilde{\zeta}](u^{(j)}-u)}{\vartheta[\widetilde{\zeta}] (0)E(u^{(j)},u)}
}{u=p^{(l)}}
\right]_{l=1}^n
\\ = &
- \widetilde{\bf u}^{\times}_l(v) \sigma_1 
\left[
\sum_{j=1}^{n} \frac{h(u^{(j)})}{dy_1(u^{(j)})}
\frac{\vartheta[\widetilde{\zeta}](u^{(j)}-p^{(l)})}{\vartheta[\widetilde{\zeta}] (0)E(u^{(j)},p^{(l)})}
\right]_{l=1}^n,
\label{eqMapToMS}
\end{align}
where $ \left(  u^{(j)}  \right)_{j=1} ^ {n}$ are the points in $\Omega$ such that 
$y_1(u^{(j)}) = y_1(v)$.
Note that using the $2 \times 2$ diagonal-block structure of $\sigma_1$ (defined in \eqref{sigma12}), we have
\begin{align}
\label{eqMapToMS123}
\sigma_1 & 
\left[
\sum_{j=1}^{n} \frac{h(u^{(j)})}{dy_1(u^{(j)})}
\frac{\vartheta[\widetilde{\zeta}](u^{(j)}-p^{(l)})}{\vartheta[\widetilde{\zeta}] (0)E(u^{(j)},p^{(l)})}
\right]_{l=1}^n
\\ = &
\left[
\overline{c_1^l}\sum_{j=1}^{n} \frac{h(u^{(j)})}{dy_1(u^{(j)})}
\frac{\vartheta[\widetilde{\zeta}](u^{(j)}-\tauBa{p^{(l)}})}{\vartheta[\widetilde{\zeta}] (0)E(u^{(j)},\tauBa{p^{(l)}})}
\right]_{l=1}^{n}.
\nonumber
\end{align}
Note that using the fact that $\sigma$ as given in \eqref{sigma12}
is the product of the signed permutation matrix $\Per{y_1,y_2}$ and of a diagonal matrix, we have
\begin{align}
h \mapsto
- \sum_{j=1}^{n}\frac{h(u^{(j)})}{dy_1(u^{(j)})} \sum_{l=1}^{n}
\frac{\vartheta[\widetilde{\zeta}](\tauBa{p^{(l)}} - v)}{\vartheta[\widetilde{\zeta}](0) E(\tauBa{p^{(l)}},v)}
\overline{c_1^l}
\frac{\vartheta[\widetilde{\zeta}](u^{(j)}-\tauBa{p^{(l)}})}{\vartheta[\widetilde{\zeta}] (0)E(u^{(j)},\tauBa{p^{(l)}})}.
\label{eq12344}
\end{align}
We use a version of the collection formula as presented in \eqref{CollFormLimitVer} to simplify the second summation in \eqref{eq12344},
which is equal to $-dy_1(v)$ whenever $u^{(j)} = v$ and zero otherwise. Hence one may conclude that
\begin{align*}
h \mapsto & - \sum_{j=1}^{n}h(u^{(j)})(- \delta_{u^{(j)},v} )= h(v),
\end{align*}
where $\delta_{u^{(j)},v}$ denotes the Kronecker delta, that is, 
$\delta_{u^{(j)},v} = 1$ if and only if $v = u^{(j)}$.

\begin{Rk}
\label{rkOmegaSym}
Note that the section defined by the mapping \eqref{ModelMap}, 
can be extended analytically (with bounded point evaluations) to 
$\Omega \cup (X \setminus \pi ^{-1} (\Spec{M^{y_1},M^{y_2}}))$.
Since the operators $M^{y_1},M^{y_2}$ have finite dimensional non Hermitian parts,
the spectrum of $M^{y_1}$ and $M^{y_2}$ contains real points and isolated non-real points.
Therefore, we can assume without loss of generality that also $X \setminus \Omega$ contains real points and isolated non-real points.
In particular, $\Omega \cap \tauBa{\Omega}$ is open, connected and symmetric,
and therefore, without loss of generality, we may assume that $\Omega$ is symmetric.
\end{Rk}

\begin{Step}
\label{step6A}
The reproducing kernel of $\mathcal{X}$ restricted to $\Omega \setminus \pi^{-1}(\Spec{M^{y_1},M^{y_2}})$
is given in terms of the complete characteristic function $W$ \eqref{CCF} of the vessel \eqref{vessel123} by
\begin{align}
\label{eqRkS}
K_{\mathcal X}(p,q) = & 
- \frac{\widetilde{\bf u}^{\times} _l (p)(\xi\sigma)\widetilde{\bf u}^{\times} _l (q)^* }{-i  (\xi y(p) - \xi \overline{y(q)}) } +
\\ & \nonumber
\frac
{\widetilde{\bf u}^{\times} _l (p) ( \xi \sigma ) W(\xi_1,\xi_2,\xi y(p))( \xi \sigma )^{-1} 
W(\xi_1,\xi_2,\xi y(q))^* (\xi\sigma)\widetilde{\bf u}^{\times} _l (q)^*} 
{-i (\xi  y(p) - \xi  \overline{y(q)}) }.
\end{align}
\end{Step}
$\mathcal X$, by assumption, is a reproducing kernel Hilbert space and in particular
\begin{equation}
\label{qwe123A}
h(p) =\innerProductReg{h}{K_{\mathcal X}(\cdot,p)}_{\mathcal X}.
\end{equation}
On the other hand, the mapping to the model space is, by Step \ref{step4}, the identity.
Thus, we have
\begin{equation}
\label{qwe123B}
\widetilde{\bf u}^{\times} _l (p) (\xi_1 \sigma_1 + \xi_2 \sigma_2) \Phi 
(\xi_1A_1+\xi_2A_2 - \xi_1 y_1(p) -\xi_2y_2(p))^{-1}h 
=  h(p).
\end{equation}
Combining equations \eqref{qwe123A} and \eqref{qwe123B}, one may conclude that the reproducing kernel 
can be expressed explicitly in terms of the model space mapping by:
\begin{align}
\label{eqsec5B}
K_{\mathcal X}(p,q) = & \widetilde{\bf u}^{\times} _l (p) (\xi \sigma) \Phi(\xi A - \xi y(p))^{-1}\times 
\\ & \nonumber \hspace{21mm}
(\xi A - \xi y(q))^{-*}\Phi^* (\xi \sigma)^* \widetilde{\bf u}^{\times} _l (q)^*.
\end{align}
A classical computation in the single-operator colligation setting 
(a similar computation can be found in \cite[Chapter 10]{KLMV}) 
yields the following relation:
\begin{align}
\nonumber
(\xi\sigma)
\Phi
(\xi A - \xi y(p))^{-1}
&
(\xi A - \xi y(q))^{-*}
\Phi^*
(\xi \sigma)^*
\\ & =
\frac{\xi\sigma W(\xi_1,\xi_2,\xi y(p)) (\xi\sigma)^{-1}W(\xi_1,\xi_2,\xi y(q))^* \xi\sigma - \xi\sigma}{-i(\xi y(p) - \xi \overline{y(q)})},
\label{eqsec5A}
\end{align}
where $\xi y = \xi_1 y_1 + \xi_2 y_2$, $\xi \sigma = \xi_1 \sigma_1 + \xi_2 \sigma_2$
and $W(\xi_1,\xi_2,\cdot)$ is the complete characteristic function \eqref{CCF}.
Finally, we substitute \eqref{eqsec5A} in \eqref{eqsec5B} to get \eqref{eqRkS}.
\smallskip

A corollary from Step \ref{step6A} is given below and is used later in this proof.

\begin{Cy}\label{cySpecC}
Let $p_0$ be a point on $X$ such that
all the pre-images of $\xi_1 y_1(p_0)+\xi_2 y_2(p_0)$ with respect to 
$\xi_1 y_1+ \xi_2 y_2 $ belong to $\Omega$ for some $(\xi_1,\xi_2)$.
Then $(y_1(p_0),y_2(p_0))$ does not belong to the joint spectrum of $M^{y_1}$ and $M^{y_2}$.
\end{Cy}
\begin{pf}
Let us choose $\alpha$ such that
the pre-images with respect to $\xi y$ of $\alpha$, assumed to be distinct and denoted by $p_1,\ldots,p_n$, belong to $\Omega$.
Then $R^{(\xi y)}_\alpha f$ is a well defined section of $L_{\widetilde{\zeta}}\otimes \Delta$ on $\Omega$
for any $f$ in $\mathcal X$, but this does not mean that $\mathcal X$ is invariant under $R^{(\xi y)}_\alpha$.
As a consequence, it does not a-priori imply that $\alpha$ does not belong to the spectrum of $M^{(\xi y)}$.
\smallskip

We use \eqref{eqRkS} to show that $W(\xi_1,\xi_2,z)$ can be extended analytically to a neighborhood of $\xi y(p_0)$.
First, we take $\beta$ in a neighborhood of infinity and, as a consequence, 
the pre-images $q_1,\ldots,q_n$ of $\beta$ under $\xi y$ belong to $\Omega$. 
Notice that the rows ${\mathbf u}^\times_l(q_j)$ for $j=1,\ldots,n$ 
are linearly independent and form an invertible matrix (follows by \eqref{CollFormLimitVer} and \eqref{CollFormLimitVer3BLA}).
Furthermore, we assume without loss of generality that $\xi\sigma$ invertible.
Hence, using \eqref{eqRkS}, 
we can analytically express $F(p)\defEq \widetilde{\bf u}^{\times} _l (p) W(\xi_1,\xi_2,\xi y(p))$
in terms of $K_{\mathcal X}(p,q_i)$, $W(\xi_1,\xi_2,\xi y(q_i))$ and $\widetilde{\bf u}^{\times} _l (q_i)$ where $i= 1,\ldots,n$
as well as ${\mathbf u}^\times_l(p)$ and $\xi y(p) - \beta$.
It follows that $F(p)$ is analytic on $\Omega$ except possibly for poles at the poles of $\xi y$.
\smallskip

We choose $\alpha$ to be an element in a punctured neighborhood of $\alpha_0=\xi y(p_0)$.
We may assume that $\alpha$ has $n$ distinct pre-images $p_1,\ldots,p_n$ with respect to $\xi y$
which all belong to $\Omega$. We then consider \eqref{eqRkS} and define
the matrix $\widetilde{\mathcal{U}}^{\times}_{l}(\alpha)$, with rows the values of
$\widetilde{\bf u}^\times_l$ at the pre-images of $\alpha$
(similarly, we define the matrix $\widetilde{\mathcal{U}}^{\times}(\alpha)$
with columns the values of $\widetilde{\bf u}^\times$ at the pre-images of $\alpha$).
The matrix $\widetilde{\mathcal{U}}^{\times}_{l}(\alpha)$ 
is invertible and the inverse is given by (see \eqref{eqMatNormSecRelation})
$
\widetilde{\mathcal{U}}^{\times}_{l}(\alpha_0)^{-1} 
= 
(\xi \sigma) 
\widetilde{\mathcal{U}}^{\times}(\alpha) \DiagTwo{j}{\frac{1}{dy(p_j)}}$.
Then 
\begin{equation}
\label{eqFp}
W(\xi_1, \xi_2 , \alpha) = 
(\xi \sigma) \widetilde{\mathcal{U}}^{\times}(\alpha)\DiagTwo{j=1,\ldots n}{\frac{1}{\xi dy(p_j)}}{\col_{j=1,\ldots n}F(p_j)}.
\end{equation}

If $(y_1(p_0),y_2(p_0))$ is not a singular point of the algebraic curve $C$, we can, without loss of generality, 
take $\xi_1$ and $\xi_2$ such that $ \alpha_0$ has $n$ distinct pre-images. 
By continuity, the same is then true for $\alpha$ in a neighborhood of $\alpha_0$.
It then follows from \eqref{eqFp} that $W(\xi_1,\xi_2,z)$ can be extended analytically to a full neighborhood of $\alpha_0$.
\smallskip

We turn now to the case where $(y_1(p_0),y_2(p_0))$ is singular. 
We use a local coordinate $t$ centered at $p_0$ so that $\xi y(p)=t(p)^r$ where $r$ is the ramification index of $\xi y$ at $p_0$.
Then \eqref{eqFp} can be rewritten in the form
\begin{equation}
\label{eqFpA}
W(\xi_1, \xi_2 , \alpha) = 
(\xi \sigma) \sum_{i=1}^{n}{\frac{ \widetilde{\mathcal{U}}^{\times}(p_i)F(p_i)}{\xi dy(p_i)}}.
\end{equation}
We may assume that as $\alpha$ goes to $\alpha_0$, the points $p_1,\ldots,p_r$ go to $p_0$ whereas $p_{r+1},\ldots,p_n$ remain distinct.
Then the summation in \eqref{eqFpA} from ${r+1}$ to $n$ extends to an analytic function on a full neighborhood of $\alpha_0$.
For the sum of the first $r$ terms $p_{1},\ldots,p_r$, correspond to $t, \e t,\ldots, \e^{r-1} t$, where $\e$ is a primitive $r$-th root of unity.
All the elements in the matrix $\phi(p) = \widetilde{\bf u}^{\times}(p)F(p)$ are analytic near $p_0$,
so that, in terms of the local coordinate $t$, we may write $\phi(t)=\sum_{k=0}^{\infty} \phi_k  t^k$.
Then the sum of the first $r$ terms becomes
\begin{align*}
\sum_{j=0}^{r-1} \frac{\phi(\e^j t)}{y'(\e^j t)}
= &
\sum_{j=0}^{r-1} \frac{\sum_{k=0}^{\infty} \phi_k (\e^j t)^k}{r(\e^j t)^{r-1}}
=
\frac{1}{r t^{r-1}} \sum_{k=0}^{\infty} \phi_k t^k \sum_{j=0}^{r-1}  (\e^j)^{-(r-1)}  (\e^j)^k
\\ = &
\frac{1}{r} \sum_{k=0}^{\infty} \phi_k t^{k-r+1} \sum_{j=0}^{r-1}  (\e^{k+1})^j
.
\end{align*}
The inner summation vanishes as long as $\e^{k+1} \neq 1$, 
as it is the sum of $k+1$-st powers of all $r$-th roots of unity
In particular, the inner summation vanishes whenever $k < r-1$
and thus the negative-index coefficients are zero.
It follows that $W(\xi_1,\xi_2,z)$ can be extended to 
$\alpha_0$ and hence to a (full) neighborhood of $\alpha_0$.
\smallskip

In both case, the singular and the non-singular cases,  $W(\xi_1,\xi_2,z)$ can be extended analytically to $\xi y(p_0)$.
It is well known, see for instance \cite{MR20:7221}, that if the characteristic function of an irreducible colligation 
can be extended analytically to a (full) neighborhood of $\alpha_0$, then $\alpha_0$ lies outside the spectrum of the operator.
By Step \ref{step4}, the mapping \eqref{ModelMap} is injective and hence by Proposition \ref{lemmModelSpaceMapProp123} the vessel $\mathcal V$ is irreducible,
it follows that the single-operator colligation derived from $\mathcal V$ in the direction $(\xi_1,\xi_2)$ is also irreducible.
Hence $\xi y(p_0)$ does not belong to the spectrum of $\xi_1M^{y_1}+\xi_2M^{y_2}$ and therefore $(y_1(p_0),y_2(p_0))$ does not belong 
to the joint spectrum of $M^{y_1}$ and $M^{y_2}$.
\end{pf}

\begin{Step}
\label{stepMax}
The input and the output determinantal representations 
of the vessel \eqref{vessel123} are fully saturated.
\end{Step}

For the definitions of maximality and fully saturated,
we refer to Section \ref{sub22} above, see also \cite{MR97m:30051,MR1634421}.
The output determinantal representation of \eqref{vessel123} is maximal and fully saturated by construction.
We proceed to prove maximality and full saturation of the input determinantal representation. First, we present several preliminary results.
\smallskip

\begin{lem}
\label{gExists}
Let $X$ be a compact Riemann surface of genus $g$ which is the normalization of 
a curve $C \subset \mathbb P ^2$ with the bi-rational embedding $\pi: X \to \mathbb P^2$.
Let $\Omega$ be an open subset of $X$ containing the pre-images of the singular points of $C$ and let $p \in \Omega$.
Then there exists a polynomial $h \in \mathbb C [z_1,z_2]$ such that $h(\pi(p)) = 0$ and all the zeros of $h \circ \pi$ on $X$ are in $\Omega$.
\end{lem}
\begin{pf}
We denote by $\mu_k$ the Abel-Jacobi mapping from $X^{(k)}$ to $J(X)$, 
sending an effective divisors of degree $k$ to the Jacobian
and we continue to denote by $n$ the degree of $C$.
Let us assume, without loss of generality, that the base point of the Abel-Jacobi mapping $p_0$ belongs to $\Omega$ and,
furthermore, we choose $p_0$ which does not belong to the finite set of the Weierstrass points.
\smallskip

Let $k>0$ be the degree of a polynomial $h \in \mathbb C [z_1,z_2]$ and
let $D$ be the effective divisor of $h$.
We know that $D$ is linear equivalent to $k \cdot L$, where $L$ is a divisor of straight line, that is,
\begin{equation}
\label{gDivisor}
D \equiv k \, \cdot \, L.
\end{equation}
Conversely, let $D$ be an effective divisor of degree $k \, n$ satisfying  \eqref{gDivisor} 
and containing the divisor of singularities, see the discussion preceeding \eqref{eqOmegaDef}.
Then, for sufficiently large $k$, $D$ is the divisor of some homogeneous polynomial of degree $k$
(by the completeness of the linear system of the adjoint curves, see e.g. \cite{fulton}). 
\smallskip

Let $D$ be a divisor that contains $p$ and the singular points, that is $D \geq p +D_{\rm sing}$.
We define an effective divisor $D^\prime$ by
\[
D^\prime = D - p - D_{\rm sing}.
\] 
By assumption $\Omega$ contains the pre-images of the singularities of $C$ and $p\in \Omega$ and then
it follows that $D$ is supported in $\Omega$ if and only if $D^\prime$ is supported in $\Omega$.
\smallskip

Recall that by the Abel-Jacobi theorem, \eqref{gDivisor} holds if and only if
\begin{equation}
\label{gDivisor12}
\mu_{k\,n}(D) = k \, \cdot \, \mu_{n} (L),
\end{equation}
and hence, using \eqref{gDivisor12}, we have
\begin{align*}
\mu_{Deg(D^\prime)}(D^\prime) 
= 
k \, \mu_n(L) - \mu_{(1+\deg (D_{\rm sing}))}(p + D_{\rm sing}).
\end{align*}

We turn to consider the divisor $g \cdot p_0$. 
This divisor, since $p_0$ is the base point of the Abel-Jacobi mapping, is zero under the mapping $\mu_{g}(\cdot)$.
Furthermore, since $p_0$ is not a Weierstrass point, $\mu_g(\cdot)$ is invertible at $g \cdot p_0$.
Therefore, using that $p_0 \in \Omega$, $\mu_g(\Omega^{(g)})$ contains an open ball around zero
(where $\Omega^{(\ell)}$ stands for the set of effective divisors of degree $\ell$ with support in $\Omega$).
\smallskip

For any $j>0$ we have $ j \mu_g(\Omega^{(g)}) \subseteq  \mu_{j \cdot g}(\Omega^{(jg)}) $
and hence for a sufficiently large $j_0$, the equality $\mu_{j_0 \cdot g}(\Omega^{(j_0 g)}) = J(X)$ holds.
As a consequence, $\mu_{l}(\Omega^{(l)}) = J(X)$ for $l \geq j_0 g$.
It remains to choose $k$ such that $k n - 1 - \deg \left(D_{\rm sing}\right) \geq j_0 g$.
\end{pf}

The following result appeared in \cite[Section 2]{MR97m:30051}) without a proof.

\begin{lem}
\label{lemMaxIO}
Let us assume that the singular points of the discriminant curve $C$ 
lie outside the joint spectrum of $\left( A_1 , A_2 \right)$.
Then the following statements hold:
\begin{enumerate}
\item If either the input or the output determinantal representations of a vessel is maximal then so is the other.
\item If either the input or the output determinantal representations of a vessel is fully saturated then so is the other.
\end{enumerate}
\end{lem}
\begin{pf}
Assume that the output determinantal representation is maximal and all affine singular points $\lambda$ 
lie outside the joint spectrum of $\left( A_1 , A_2 \right)$.
Then, $S(\lambda)$, for any $\lambda$ outside the spectrum of $\left( A_1 , A_2 \right)$, 
is a linear mapping from $\E(\lambda)$ to $\widetilde{\E}(\lambda)$.
Similarly, the JCF of the adjoint vessel $\widetilde{S}(\lambda)$ (see \cite{MR2043236}), 
for any $\lambda$ outside the spectrum of $\left( A_1 , A_2 \right)$,
is a linear mapping from $\widetilde{\E}_l(\lambda)$ to $\E_l(\lambda)$.
Hence, $\E(\lambda)$ and $\widetilde{\E}(\lambda)$ share the same dimension for all affine singular points $\lambda \in C$.
For the points at infinity, the fibers of the input and output determinantal representations simply coincide,
since the input and output determinantal representations share the same matrices $\sigma_1$ and $\sigma_2$. 
In particular, they have the same dimension.
\smallskip

Let us now assume that the output determinantal representation is fully saturated.
Then, for a singular affine point $\lambda_0=(\lambda^0_1,\lambda^0_2) \in C_0$, we choose $\xi_1$ and $\xi_2$ 
such that $\xi\sigma=\xi_1\sigma_1+\xi_2\sigma_2$ is invertible and 
such that
$\xi_1 \lambda^0_1 + \xi_2 \lambda^0_2$ does not belong to spectrum of $\xi_1 A_1 + \xi_2 A_2$,
by continuity this is also true in a neighborhood of $\lambda^0$.
Then, one may consider the identity (see for instance \cite{KLMV})
\begin{align} \nonumber
\left(
\lambda_1\sigma_2 - \lambda_2\sigma_1 + \widetilde{\gamma}
\right)
&
W(\xi_1,\xi_2, {\xi_1\lambda_1 + \xi_2 \lambda_2} ) 
= \\ &
\label{eqBLa5353}
\widetilde{W}(\xi_1,\xi_2, \xi_1\lambda_1 + \xi_2 \lambda_2) 
\left(
\lambda_1\sigma_2 - \lambda_2\sigma_1 + \gamma
\right)
\end{align}
where
\[
\widetilde{W}(\xi_1,\xi_2,z) = I - i (\xi_1 \sigma _1 + \xi_2 \sigma_2) \Phi(\xi_1 A_1 + \xi_2 A_2 -zI)^{-1} \Phi ^*,
\]
and hence, using \eqref{eqWinA}, we have:
\begin{align*}
\left(
\lambda_1\sigma_2 - \lambda_2\sigma_1 + \widetilde{\gamma}
\right)
= 
&
\widetilde{W}(\xi_1,\xi_2, \xi_1\lambda_1 + \xi_2 \lambda_2 ) 
\left(
\lambda_1\sigma_2 - \lambda_2\sigma_1 + \gamma
\right) \times
\\ & (\xi\sigma) W(\xi_1,\xi_2, \xi_1 \overline{\lambda_1} + \xi_2 \overline{\lambda_2} )^* (\xi\sigma)^{-1}.
\end{align*}
The adjoint operator is anti-multiplicative, that is, it satisfies $\operatorname{adj} (AB) = \operatorname{adj} (B) \operatorname{adj} (A)$, and therefore:
\begin{equation}
\label{eqBLa535A}
\widetilde{V}(\lambda)
=
\operatorname{adj}  \widetilde{W}(\xi_1,\xi_2, \xi_1\lambda_1 + \xi_2 \lambda_2 )  \, 
V(\lambda) \, 
\operatorname{adj}  (\xi\sigma)W(\xi_1,\xi_2, \xi_1 \overline{\lambda_1} + \xi_2 \overline{\lambda_2} )^*(\xi\sigma)^{-1}.
\end{equation}
Since $\widetilde{W}(\xi_1,\xi_2, \xi_1\lambda_1 + \xi_2 \lambda_2 )$ and 
$W(\xi_1,\xi_2, \xi_1 \overline{\lambda_1} + \xi_2 \overline{\lambda_2} )^*$
are by assumption analytic in a neighborhood of $\lambda_0$,
it follows by \eqref{eqBLa535A} that if all entries of $V(\lambda)$ vanish to a certain order, so do all the entries of $\widetilde{V}(\lambda)$
and hence if the input determinantal representation is fully saturated,
so is the output one.
Similarly, the second direction follows by multiplying \eqref{eqBLa5353} on the left by the inverse of
$\widetilde{W}(\xi_1,\xi_2,\xi_1\lambda_1 + \xi_2\lambda_2)$.

To understand the behavior at infinity, we consider the homogeneous version of \eqref{eqBLa5353}, that is
\begin{align*}
\left(
\nu_1\sigma_2 - \nu_2\sigma_1 + \nu_0\widetilde{\gamma}
\right)
&
W(\nu_0\xi_1,\nu_0\xi_2, {\xi_1\nu_1 + \xi_2 \nu_2} ) 
= \\ &
\widetilde{W}(\nu_0\xi_1,\nu_0\xi_2, \xi_1\nu_1 + \xi_2 \nu_2) 
\left(
\nu_1\sigma_2 - \nu_2\sigma_1 + \nu_0\gamma
\right).
\end{align*}
Then, we use a different affine chart and repeat the argument above.
\end{pf}

Using the previous results, we may conclude the main argument of this step.

\begin{pf}[of Step \ref{stepMax}]
Using Lemma \ref{lemMaxIO}, it is enough to show that for any $p_0 \in \Omega$ the point $(y_1(p_0),y_2(p_0))$ 
lies outside the joint spectrum of $M^{y_1}$ and $M^{y_2}$.
\smallskip

Using Lemma \ref{gExists} and assuming $p_0\in\Omega$, there exists a polynomial $g$ satisfying $g(\pi(p_0))=0$ 
such that the entire fiber of the meromorphic function $g(y_1,y_2)$ above $0$ belongs to $\Omega$.
We define a new meromorphic function $w_1(p) = g(y_1(p),y_2(p))$.
We assume without lose of generality that $w_1$ is real, otherwise we replace $w_1(p)$ by $w_1(p)\overline{w_1(\tauBa{p})}$
(it is well-defined since we may assume that $\Omega$ is symmetric, see Remark \ref{rkOmegaSym}).
Then we define $w_2(p) = h(y_1(p),y_2(p))$ for some two-variable polynomial $h$ with real coefficients
such that $w_1$ and $w_2$ generate $\mathcal{M} (X)$.
\smallskip

We now take $\alpha_1$ and $\alpha_2$ in a neighborhood of infinity and such that 
$w^\prime_1 \defEq \frac{1}{w_1(z)- \alpha_1}$ and $w^\prime_2 \defEq \frac{1}{w_2(z)- \alpha_2}$
have simple poles.
Then all the pre-images of $\frac{-1}{\alpha_1}$ under $w^\prime_1$ lie in $\Omega$.
Since the structure identity holds for $y_1$ and $y_2$,
then by Proposition \ref{myGeneral1} the colligation condition hold for $M^{y_1}$ and $M^{y_2}$.
Therefore, by Theorem \ref{lemCollCodI}, $M^{w_1}$ satisfies the colligation condition
and then again by Proposition \ref{myGeneral1} the structure identity holds for $w_1$.
We then take $\beta'$ in a neighborhood of infinity and let $\beta = \frac{1+\beta'\alpha_1}{\beta'}$ 
in a neighbourhood of $\alpha_1$.
Then one can show that
$\frac{1}{w_1'-\beta'} = -\frac{1}{\beta'}\left(
1+\frac{1}{\beta'}\frac{1}{w_1-\beta}\right)$ and therefore
$R^{w_1'}_{\beta'} = -\frac{1}{\beta'}\left(
I+\frac{1}{\beta'}R^{w_1}_\beta\right)$.
It follows by a straightforward verification 
that if the structure identity holds for $w_1$ then it also true for $w^\prime_1$ (and similarly for $w^\prime_2$).
Hence, repeating Steps (\ref{step1}-\ref{step6A}) but now with $M^{w^\prime_1}$ and $M^{w^\prime_2}$ implies that the collection
\[
( \, M^{w^\prime_1} \, ,  \,  M^{w^\prime_2} \, ; \, \mathcal X \, , \,  \Phi_{(w^\prime_1,w^\prime_2)} \, , \, E_{(w^\prime_1,w^\prime_2)} \, ; \,  
\sigma_{w^\prime_1} \, , \,  \sigma_{w^\prime_2} \, , \, \gamma_{(w^\prime_1,w^\prime_2)}  \, , \, \widetilde{\gamma}_{(w^\prime_1,w^\prime_2)} \, ),
\]
is again a commutative two-operator vessel.
Since all the pre-images of $w_1'(p_0)=\frac{-1}{\alpha_1}$ 
under $w_1'$ lie in $\Omega$,
Corollary 7.3 implies that $(w_1'(p_0),w_2'(p_0))$  
does not belong to the joint spectrum of $M^{w^\prime_1}$ and $M^{w^\prime_2}$.
\smallskip

Assume that the statement is not true, that is, $(y_1(p_0),y_1(p_0))$ belongs to the joint the spectrum of $M^{y_1}$ and $M^{y_2}$.
Then by the spectral mapping theorem for a pair of commuting operators,
$(w_1(p_0),w_1(p_0))$ belongs to the joint of spectrum $M^{w_1}$ and $M^{w_2}$ (see \cite{MR1320541} and Theorem \ref{rkModelOpAlg}).
Furthermore, using the spectral mapping theorem for rational function \cite{MR42:3603},
$(w^\prime_1(p_0),w^\prime_1(p_0))$ belongs to the joint spectrum of $M^{w^\prime_1}$ and $M^{w^\prime_2}$, a contradiction.
\end{pf}

\begin{Step}
The reproducing kernel of $\mathcal{X}$ is equal to
\begin{align*}
K_{\mathcal X}(p,q) = T(p) K_\zeta(p,q) T(q)^* - K_{\tilde\zeta}(p,q),
\end{align*}
where $T$ is the normalized joint characteristic function
associated to the vessel \eqref{vessel123}.
\end{Step}
Using Step \ref{stepMax}, the input and output determinantal representations are maximal and hence 
we turn to examine the NJCF (the normalized joint characteristic function, which is related to the JCF by \eqref{st}).
Using \eqref{eqRkS} and repeating in the reverse order the calculation at 
the end of the proof of Theorem \ref{vesselAreEq}
yields the required result.
\end{pf}

We now turn to present the proof of Theorem \ref{MainTh}.\\

\begin{pf}[of Theorem \ref{MainTh}]
Let $y_1(\cdot)$ and $y_2(\cdot)$ be two meromorphic functions,
not necessarily with simple poles, generating $\mathcal{M} (X)$.
The corresponding bi-rational embedding of $X$ onto a curve $C$
defined by the closure in $\mathbb P ^2$ of the curve $C_0$ is given by
\begin{align*}
\pi \colon X & \to  C_0 \subseteq {\mathbb C}^2
\\
p & \mapsto (y_1(p),y_2(p)).
\end{align*}
Let us consider the pair of meromorphic functions defined by
\[
\widetilde{y}_k (p) \defEq \frac{1}{y_k(p) - \alpha_k^0}
\qquad
\alpha_k^0 \in \mathbb R
\qquad {\rm for} \qquad k=1,2.
\]
Since $y_k$ has a finite number of critical points,
we can choose $\alpha_k^0$ in a neighborhood of $\beta_k^0$ such that the pre-images 
of $\alpha_1^0$ and $\alpha_2^0$ are non critical points of $y_1$ and $y_2$ and
the functions $\widetilde{y}_1 (p) $ and $\widetilde{y}_2 (p)$ have only simple poles.
\smallskip

Note that $\mathcal X$ is invariant under $R^{y_k}_{\alpha_k^0}$ implies 
$\mathcal X$ is invariant under $M^{\wt{y}_k}$.
Furthermore, the same verification as in the proof of Step \ref{stepMax} of Theorem \ref{preTh} above
shows that the validity of the structure identity for $y_k$ 
(for some $\alpha_k,\beta_k$ in a neighborhood of $\alpha_k^0$)
implies that the structure identity holds for $\widetilde{y}_k$
(for some $\widetilde{\alpha}_k,\widetilde{\beta}_k$ in a neighborhood of infinity).
\smallskip

The functions $\widetilde{y}_1 (p)$ and $\widetilde{y}_2 (p)$ define a new birational 
embedding of $X$ into $\mathbb P ^2$, but onto a different curve, denoted by $\widetilde C$, 
defined by the closure of $\widetilde C_0$:
\begin{align*}
\widetilde \pi \colon X & \to \widetilde C_0 \subseteq {\mathbb C}^2
\\
x & \mapsto ((y_1(p)-\alpha_1^0)^{-1},(y_2(p)-\alpha_2^0)^{-1})
.
\end{align*}
The new critical points of $\widetilde\pi$ 
(except the critical points inherited from $\pi$) 
can be only the points above the poles of 
$(y_1(p)-\alpha_1^0)^{-1}$ and $(y_2(p)-\alpha_2^0)^{-1}$.
However, by assumption, the pre-images of $\alpha_1^0$ and $\alpha_2^0$ are regular.
Thus
$(\widetilde \pi)^{-1} \widetilde C{\text{\rm sing}} \subseteq (\pi)^{-1} C_{\text{\rm sing}} $ 
and the regularity assumption of $(y_1(p),y_2(p))$ 
implies regularity of $(\widetilde{y}_1(p),\widetilde{y}_2(p))$.
\smallskip

We now can apply Theorem \ref{preTh}. We embed the operators $M^{\widetilde{y}_1}$ and $M^{\widetilde{y}_2}$ into a commutative vessel of the form
\begin{equation*}
\mathcal{V}^{\prime} =
( \, M^{(y_1-\alpha_1^0)^{-1}} \, , \, M^{(y_2-\alpha_2^0)^{-1}} \, ;  \, \mathcal X \, , \,  \Phi \, , \, E \, ; \, \sigma_1 \, , \, \sigma_2 \, , \, \gamma, \, \widetilde{\gamma} \, ).
\end{equation*}
Then the reproducing kernel is of the form
\begin{equation*}
K(p,q)
=
T^{\prime}(p)
K_{\zeta}(p,q)
T^{\prime}(q)^*
-
K_{\widetilde{\zeta}}(p,q)
,
\end{equation*}
where $T^{\prime}$ is the normalized joint characteristic function of the vessel $\mathcal{V}^{\prime}$.
This completes the proof of part \ref{stateA}.
\smallskip

We now turn to prove part \ref{stateB}.
Let $y_1 = y$ be a real meromorphic function,
we choose $y_2$, another real meromorphic function, 
such that all the poles of $y_2$ are contained in 
$\Omega$ and $y_1$ and $y_2$ generate $\mathcal{M} (X)$.
The pair of meromorphic functions defines 
a birational embedding onto a curve $C$ in $\P^2$.
\smallskip

Then we apply the realization theorem for the NJCF $T$ (i.e. Theorem \ref{realizationTh})
and thus we have a commutative two-operator vessel $\mathcal V$ with 
input and output determinantal representations corresponding to $\z$ and $\wt{\z}$ and 
with the normalized joint characteristic function $T$.
The associated model vessel is given by
\[
\mathcal{V}_{T} 
= 
( \, M^{y_1} \, , \, M^{y_2}  \, ; \,  \mathcal{H}(T) \, , \, \Phi_{\rm mod} \, , \, \mathbb C ^n \, ;
\,  \sigma_1 \, , \, \sigma_2 \, , \, \gamma \, , \,  \widetilde{\gamma} \, )
.
\]
Then, using Theorem \ref{vesselAreEq},
$\mathcal{V}_{T}$ is an irreducible commutative two-operator vessel which is unitary equivalent, 
on its principal subspace, to $\mathcal{V}$ and with the normalized joint characteristic function $T$.
\smallskip

Furthermore, $\mathcal X$ is invariant under $M^{y_1}$ and hence also $R_\alpha^{y_1}$-invariant 
for $\alpha$ in the neighborhood of infinity.
To complete the proof, we mention that by Proposition \ref{myGeneral1},
the structure identity is equivalent to the colligation condition. 
Hence the colligation condition for $M^{y_1}$ in $\mathcal{V}_{T}$ 
implies that the structure identity for $y_1(\cdot) = y(\cdot)$ holds.
\end{pf}

We note that in the proof of Theorem \ref{MainTh}, if $y$ has only simple poles then 
we can choose $y_2$ with simple poles.
Furthermore, we can take the output determinantal representation to be the canonical determinantal
corresponding to $\wt{\z}$.
Then the input determinantal representation is the canonical determinantal representation corresponding to $\z$
multiplied the by the values of $T$ at infinity, see Remark \ref{rk26D},
and $\Phi_{\rm mod}$ is just the evaluation operator at the poles of $y_1$ and $y_2$.


\ifdefined \elsarticle
    \section{\refname}    
\fi

\def\cfgrv#1{\ifmmode\setbox7\hbox{$\accent"5E#1$}\else
  \setbox7\hbox{\accent"5E#1}\penalty 10000\relax\fi\raise 1\ht7
  \hbox{\lower1.05ex\hbox to 1\wd7{\hss\accent"12\hss}}\penalty 10000
  \hskip-1\wd7\penalty 10000\box7} \def\cprime{$'$} \def\cprime{$'$}
  \def\cprime{$'$} \def\lfhook#1{\setbox0=\hbox{#1}{\ooalign{\hidewidth
  \lower1.5ex\hbox{'}\hidewidth\crcr\unhbox0}}} \def\cprime{$'$}
  \def\cprime{$'$} \def\cprime{$'$} \def\cprime{$'$} \def\cprime{$'$}
  \def\cprime{$'$}


\end{document}